\documentclass[11pt]{article}
  
\title{Enumeration and investigation of acute 0/1-simplices modulo the action of the hyperoctahedral group}
\date{\today} 
\author{Jan Brandts and Apo Cihangir} 
       
\usepackage{a4wide,authblk,epsfig,amsmath,amssymb,makeidx,graphics,color,tikz}
\usepackage{caption}
\begin{document}              
         
\newtheorem{Th}{Theorem}[section]       
\newtheorem{Le}[Th]{Lemma}     
\newtheorem{Co}[Th]{Corollary}     
\newtheorem{Pro}[Th]{Proposition}               
\newtheorem{Ex}[Th]{Example}   
\newtheorem{Def}[Th]{Definition}       
\newtheorem{rem}[Th]{Remark}             
\newcommand{\be}{\begin{equation}}         
\newcommand{\ee}{\end{equation}}        
\newcommand{\RR}{\mathbb{R}}       
\newcommand{\half}{\frac{1}{2}}  
\newcommand{\hdrie}{\hspace{3mm}} 
\newcommand{\und}{\hdrie\mbox{\rm and }\hdrie} 
\newcommand{\sth}{\hdrie | \hdrie}
    
\def\zeros {{\bf 0}} 
\def\ones{{\bf 1}}
\def\CC {\mathcal{C}}
\def\Bc {\mathcal{B}}
\def\PP {\mathcal{P}}
\def\Bn {\mathcal{B}_n} 
\def\bb {\mathbb{b}}
\def\GG {\mathcal{G}}
\def\OO {\mathcal{O}} 
\def\II {\mathbb{I}} 
\def\BB {\mathbb{B}}  
\def\MM {\mathcal{M}}
\def\EE {\mathcal{E}}
\def\PP {\mathcal{P}}
\def\NN {\mathcal{N}}
\def\SS {\mathcal{S}}
\def\ZZ {\mathcal{Z}}
\def\TT {\mathcal{T}}
\def\AA {\mathcal{A}}
\def\csim {\overset{c}{\sim}} 
\def\psim {\overset{p}{\sim}} 
\def\zsim {\overset{{\small 0/1}}{\sim}} 
\def\sort {{\rm sort}}
\def\Eqv {\hdrie\Leftrightarrow\hdrie} 
\def\Bnk {\BB^{n\times k}}
\def\Bnn {\BB^{n\times n}}
\def\Mnn {\MM^{n\times n}}
\def\Mnk {\MM^{n\times k}}  
\def\mod {\,{\rm mod }\,}
\def\conv {\mbox{\rm conv}}
\newcommand{\supp}{{\rm supp}} 
\newcommand{\ol}{\overline} 
\newcommand{\vep}{\varepsilon}
\newcommand{\zok}{0/1-$k$} 
\newcommand{\zo}{0/1} 
\newcommand{\nmo}{n\!-\!1}
\newcommand{\npo}{n\!+\!1}
\newcommand{\npt}{n\!+\!2}
\maketitle  
   
\begin{abstract}
The convex hull of~$n+1$ affinely independent vertices of the unit~$n$-cube~$I^n$ is called a {\em 0/1-simplex}. \index{0/1-simplex}
It is {\em nonobtuse} if none its dihedral angles is obtuse, and {\em acute} \index{acute simplex} if additionally none of them is right. In 
terms of linear algebra, acute 0/1-simplices in~$I^n$ can be described by nonsingular 0/1-matrices~$P$ of 
size~$n\times n$ whose Gramians~$G=P^\top P$ have an inverse that is strictly diagonally dominant, with negative 
off-diagonal entries \cite{BrKoKr,BrKoKrSo}.

\smallskip

The first part of this paper deals with giving a detailed description of how to efficiently compute, by means of a 
computer program, a representative from each orbit of an acute 0/1-simplex under the action of the 
{\em hyperoctahedral group}~$\Bn$ \cite{GeKi} \index{hyperoctahedral group} of symmetries of~$I^n$. A side product of the investigations is 
a simple code that computes the {\em cycle index} of~$\Bn$, which can in explicit form only be found in the 
literature \cite{ChGu2} for~$n\leq 6$. Using the computed cycle indices for~
$\mathcal{B}_3,\dots,\mathcal{B}_{11}$ in combination with P\'olya's theory of enumeration shows that 
acute 0/1-simplices are extremely rare among all 0/1-simplices. 

\smallskip

In the second part of the paper, we study the 0/1-matrices that represent the acute 0/1-simplices that were 
generated by our code from a mathematical perspective. One of the patterns observed in the data involves 
{\em unreduced upper Hessenberg} 0/1-matrices of size~$n\times n$, block-partitioned according to certain 
{\em integer compositions} of~$n$. These patterns will be fully explained using a so-called 
{\em One Neighbor Theorem} \cite{BrDiHaKr}. Additionally, we are able to prove that the volumes of 
the corresponding acute simplices are in one-to-one correspondence with the part of 
{\em Kepler's Tree of Fractions} \cite{AiDuFi,KimMo} \index{Kepler's Tree of Fractions} that enumerates~$\mathbb{Q}\cap(0,1)$. Another 
key ingredient in the proofs is the fact that the Gramians of the unreduced upper Hessenberg matrices involved 
are {\em strictly ultrametric} \cite{DeMaMa,NabVar} matrices. 
\end{abstract}
\section{Introduction} \label{Eintro}
A 0/1-simplex is an~$n$-dimensional 0/1-polytope \cite{KaZi} with~$n+1$ vertices. Equivalently, it is the 
convex hull of~$n+1$ of the~$2^n$ elements of the set~$\BB^n$ of vertices of the unit~$n$-cube~$I^n$ 
whenever this hull has dimension~$n$. To support the mathematical studies of 0/1-simplices, and in particular 
of those whose dihedral angles are all {\em nonobtuse} or even {\em acute} \cite{BrKoKrSo}, we investigate how 
to enumerate such 0/1-simplices modulo the action of the hyperoctahedral group~$\Bn$ of symmetries of~$I^n$ 
by means of a computer program. The motivation to generate such computational data was quite appropriately 
phrased by G\"unther Ziegler in Chapter 1 of {\em Lectures on 0/1-Polytopes} \cite{KaZi}, as 
``{\em Low-dimensional intuition does not work!}\,''.

\smallskip

This statement expresses the fact that although it is tempting to formulate conjectures on~$n$-dimensional 
0/1-polytopes and related 0/1-matrices based on computational data obtained for a few small values of~$n$, 
these conjectures often fail to be true. Finding out that a conjecture is false using general mathematical arguments 
may be much harder than generating the necessary computational data for large enough~$n$ to {\em disprove} it, 
not in the least because the tendency towards a conjecture is rather to believe its validity and aim to prove it. This 
is why we concentrate on the enumeration problem for acute 0/1-simplices. Using the date produced by the 
enumeration, we will also formulate and prove some mathematical results on certain classes of 0/1-matrices. We 
will summarize the most important of these results in Section~\ref{sect-E1.2}. First, in Section~\ref{sect-E1.1}, we 
give two examples that illustrate Ziegler's claim above, also based on the computational data. 
\begin{rem}{\rm Especially the larger 0/1-matrices in this paper we will often display as a picture of an 
array with black and white squares representing its ones and zeros, respectively}.
\end{rem}
\subsection{Two examples that illustrate Ziegler's claim above}\label{sect-E1.1}
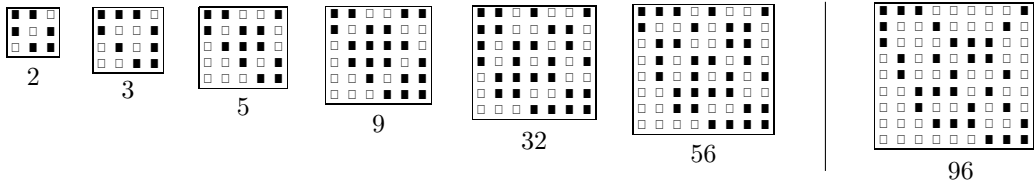
\begin{figure}[!]
\begin{center}  
\begin{tikzpicture}[scale=0.9, every node/.style={scale=0.9}]
\node[xscale=0.4, yscale=0.5] at (0,0.65) {$\begin{array}{|ccc|}
\hline
\blacksquare&\blacksquare&\square\\
\blacksquare&\square&\blacksquare\\
\square&\blacksquare&\blacksquare\\\hline
\end{array}$};
\node at (0,0) {$2$};
\node[xscale=0.4, yscale=0.5] at (1.4,0.53) {$\begin{array}{|cccc|}
\hline
\blacksquare&\blacksquare&\blacksquare&\square\\
\blacksquare&\square&\square&\blacksquare\\
\square&\blacksquare&\square&\blacksquare\\
\square&\square&\blacksquare&\blacksquare\\\hline
\end{array}$};
\node at (1.4,-0.2) {$3$};
\node[xscale=0.4, yscale=0.5] at (3.1,0.42) {$\begin{array}{|ccccc|}
\hline
\blacksquare&\blacksquare&\square&\square&\blacksquare\\
\blacksquare&\square&\blacksquare&\blacksquare&\square\\
\square&\blacksquare&\blacksquare&\blacksquare&\square\\
\square&\square&\blacksquare&\square&\blacksquare\\
\square&\square&\square&\blacksquare&\blacksquare\\\hline
\end{array}$};
\node at (3.1,-0.45) {$5$};
\node[xscale=0.4, yscale=0.5] at (5.1,0.31) {$\begin{array}{|rrrrrr|}
\hline
\blacksquare&\blacksquare&\square&\square&\blacksquare&\blacksquare\\
\blacksquare&\square&\blacksquare&\blacksquare&\square&\square\\
\square&\blacksquare&\blacksquare&\blacksquare&\blacksquare&\square\\
\square&\blacksquare&\blacksquare&\blacksquare&\square&\blacksquare\\
\square&\square&\blacksquare&\square&\blacksquare&\blacksquare\\
\square&\square&\square&\blacksquare&\blacksquare&\blacksquare\\
\hline
\end{array}$};
\node at (5.1,-0.7) {$9$};
\node[xscale=0.4, yscale=0.5] at (7.4,0.2) {$\begin{array}{|ccccccc|}
\hline
\blacksquare&\blacksquare&\square&\blacksquare&\square&\square&\blacksquare\\
\blacksquare&\blacksquare&\square&\square&\blacksquare&\blacksquare&\square\\
\blacksquare&\square&\blacksquare&\blacksquare&\square&\blacksquare&\square\\
\blacksquare&\square&\blacksquare&\square&\blacksquare&\square&\blacksquare\\
\square&\blacksquare&\blacksquare&\blacksquare&\blacksquare&\square&\square\\
\square&\blacksquare&\blacksquare&\square&\square&\blacksquare&\blacksquare\\
\square&\square&\square&\blacksquare&\blacksquare&\blacksquare&\blacksquare\\
\hline \end{array}$};
\node at (7.4,-0.95) {$32$};
\node[xscale=0.4, yscale=0.5] at (9.9,0.1) {$\begin{array}{|cccccccc|}
\hline
\blacksquare&\blacksquare&\blacksquare&\square&\blacksquare&\square&\square&\blacksquare\\
\blacksquare&\square&\square&\blacksquare&\square&\blacksquare&\blacksquare&\square\\
\square&\blacksquare&\blacksquare&\square&\square&\blacksquare&\blacksquare&\square\\
\square&\blacksquare&\square&\blacksquare&\blacksquare&\square&\blacksquare&\square\\
\square&\blacksquare&\square&\blacksquare&\square&\blacksquare&\square&\blacksquare\\
\square&\square&\blacksquare&\blacksquare&\blacksquare&\blacksquare&\square&\square\\
\square&\square&\blacksquare&\blacksquare&\square&\square&\blacksquare&\blacksquare\\
\square&\square&\square&\square&\blacksquare&\blacksquare&\blacksquare&\blacksquare\\
\hline\end{array}$};
\node at (9.9,-1.15) {$56$};
\draw (11.7,-1.4)--(11.7,1.1);
\node[xscale=0.4, yscale=0.5] at (13.6,0) {$\begin{array}{|ccccccccc|}
\hline
\blacksquare&\blacksquare&\blacksquare&\square&\square&\square&\square&\square&\blacksquare\\
\blacksquare&\square&\square&\blacksquare&\square&\square&\square&\blacksquare&\square\\
\blacksquare&\square&\square&\square&\blacksquare&\blacksquare&\blacksquare&\square&\square\\
\square&\blacksquare&\square&\blacksquare&\square&\blacksquare&\blacksquare&\square&\square\\
\square&\blacksquare&\square&\square&\blacksquare&\square&\square&\blacksquare&\square\\
\square&\square&\blacksquare&\blacksquare&\blacksquare&\square&\blacksquare&\square&\square\\
\square&\square&\blacksquare&\square&\square&\blacksquare&\square&\blacksquare&\square\\
\square&\square&\square&\blacksquare&\blacksquare&\blacksquare&\square&\square&\blacksquare\\
\square&\square&\square&\square&\square&\square&\blacksquare&\blacksquare&\blacksquare\\
\hline\end{array}$};
\node at (13.7,-1.4) {$96$};
\end{tikzpicture}    
\end{center} 
\caption{\small{0/1-matrices with maximal determinant that represent {\em acute} 0/1-simplices.}}
\label{figureE1}
\end{figure}

A first example of a statement that is valid in~$I^n$ for~$n\leq 8$ but that does not hold in~$I^9$ is the following. 
In Figure~\ref{figureE1} we display 0/1-matrices~$P$ having maximal absolute value of the determinant, when ranging over 
all those~$n\times n$ matrices whose~$n$ columns together with the origin are the vertices of a so-called 
{\em acute} 0/1-simplex. See Definition~\ref{Edef-acute} for a linear algebraic characterization acute 
0/1-simplices and the 0/1-matrices associated with them. 
 \\[2mm] 
For~$n\leq 8$, these values turn out to be even maximal when ranging over {\em all} 0/1-matrices of size~
$n\times n$. However, the maximal determinant over all~$9\times 9$ 0/1-matrices is~$144$ and not~$96$.

\smallskip

This, of course,  disproves the conjecture that maximal determinants of 0/1-matrices are attained by 
0/1-matrices that represent acute 0/1-simplices. Notwithstanding, the Hadamard maximal determinant 
conjecture \cite{Hada} is equivalent \cite{Gr} with the existence of a {\em regular} simplex in~$I^n$ for dimensions 
$n$ whose remainder after devision by~$4$ equals~$3$. Regular simplices have acute dihedral angles, and indeed, 
the~$3\times 3$ and the~$7\times 7$ matrix in Figure~\ref{figureE1} are so-called {\em Hadamard matrices}. This motivates 
a further study of acute 0/1-simplices and their determinants, as the set of acute 0/1-simplices is a small 
and highly structured set in which the Hadamard matrices figure as the most structured ones. It thus puts the 
Hadamard matrices in a wider context in which, as far as we know, they have not yet been studied.

\smallskip

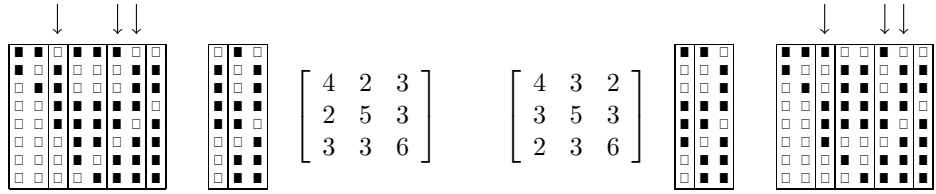
\begin{figure}[h]
\begin{center}  
\begin{tikzpicture}
\draw[->] (-0.4,1.5)--(-0.4,1.1);
\draw[->] (0.4,1.5)--(0.4,1.1);
\draw[->] (0.65,1.5)--(0.65,1.1);
\node[xscale=0.4, yscale=0.5] at (0,0) {$\begin{array}{|cc|c|cc|cc|c|}
\hline
\blacksquare&\blacksquare&\square&\blacksquare&\blacksquare&\blacksquare&\square&\square\\
\blacksquare&\square&\blacksquare&\square&\square&\square&\blacksquare&\blacksquare\\
\square&\blacksquare&\blacksquare&\square&\square&\square&\blacksquare&\blacksquare\\
\square&\square&\blacksquare&\blacksquare&\blacksquare&\blacksquare&\blacksquare&\square\\
\square&\square&\blacksquare&\blacksquare&\blacksquare&\blacksquare&\square&\blacksquare\\
\square&\square&\square&\blacksquare&\blacksquare&\square&\blacksquare&\blacksquare\\
\square&\square&\square&\blacksquare&\square&\blacksquare&\blacksquare&\blacksquare\\
\square&\square&\square&\square&\blacksquare&\blacksquare&\blacksquare&\blacksquare\\
\hline\end{array}$};
\node[xscale=0.4, yscale=0.5] at (2,0) {$\begin{array}{|c|cc|}
\hline
\square&\blacksquare&\square\\
\blacksquare&\square&\blacksquare\\
\blacksquare&\square&\blacksquare\\
\blacksquare&\blacksquare&\blacksquare\\
\blacksquare&\blacksquare&\square\\
\square&\square&\blacksquare\\
\square&\blacksquare&\blacksquare\\
\square&\blacksquare&\blacksquare\\
\hline\end{array}$};
\node[scale=0.9] at (3.7,0) {$\left[\begin{array}{ccc} 4 & 2 & 3 \\ 2 & 5 & 3\\ 3 & 3 & 6\end{array}\right]$};
\node[scale=0.9] at (6.5,0) {$\left[\begin{array}{ccc} 4 & 3 & 2 \\ 3 & 5 & 3\\ 2 & 3 & 6\end{array}\right]$};
\node[xscale=0.4, yscale=0.5] at (8.2,0) {$\begin{array}{|c|cc|}
\hline
\blacksquare&\blacksquare&\square\\
\square&\square&\blacksquare\\
\square&\square&\blacksquare\\
\blacksquare&\blacksquare&\blacksquare\\
\blacksquare&\blacksquare&\square\\
\blacksquare&\square&\blacksquare\\
\square&\blacksquare&\blacksquare\\
\square&\blacksquare&\blacksquare\\
\hline\end{array}$};
\draw[->] (9.8,1.5)--(9.8,1.1);
\draw[->] (10.6,1.5)--(10.6,1.1);
\draw[->] (10.85,1.5)--(10.85,1.1);
\node[xscale=0.4, yscale=0.5] at (10.2,0) {$\begin{array}{|cc|c|cc|cc|c|}
\hline
\blacksquare&\blacksquare&\blacksquare&\square&\square&\blacksquare&\square&\square\\
\blacksquare&\square&\square&\blacksquare&\blacksquare&\square&\blacksquare&\blacksquare\\
\square&\blacksquare&\square&\blacksquare&\blacksquare&\square&\blacksquare&\blacksquare\\
\square&\square&\blacksquare&\blacksquare&\blacksquare&\blacksquare&\blacksquare&\square\\
\square&\square&\blacksquare&\blacksquare&\blacksquare&\blacksquare&\square&\blacksquare\\
\square&\square&\blacksquare&\square&\square&\square&\blacksquare&\blacksquare\\
\square&\square&\square&\blacksquare&\square&\blacksquare&\blacksquare&\blacksquare\\
\square&\square&\square&\square&\blacksquare&\blacksquare&\blacksquare&\blacksquare\\
\hline\end{array}$};
\end{tikzpicture} 
\end{center} 
\caption{\small{Examples of 0/1-matrices representing acute 0/1-simplices whose Gramians are }
{\em not} ultrametric. Columns 3,6,7 do not satisfy the inequalities in (\ref{Eult}).}
\label{figureE2}
\end{figure}

As a second example in which low dimensional intuition does not work, here is a statement that holds in~$I^n$ for 
all~$n\leq 7$: the Gramian~$G=P^\top P$ of any 0/1-matrix~$P$ whose columns together with the origin are 
the vertices of an acute simplex is a {\em strictly ultrametric} matrix \cite{DeMaMa,NabVar}. A strictly ultrametric 
matrix is a highly structured {\em positive} matrix in the sense that all its~$3\times 3$ principal submatrices are, 
modulo simultaneous permutation of rows and columns, of the form
\be\label{Eult} \left[\begin{array}{rrr} d & b & a\\b & c & a\\a & a & f\end{array}\right], \hdrie \mbox{\rm with}\hdrie a\leq b < c \leq d \und a< f.\ee 
Even though the columns of the two~$8\times 8$ matrices displayed in Figure~\ref{figureE2} together with the origin are indeed 
vertices of acute 0/1-simplices in~$I^8$, their Gramians are however {\em not} strictly ultrametric. In both 
matrices, the inner products between columns~$3$,~$7$, and~$8$ do not satisfy the relations in (\ref{Eult}), not even 
after simultaneous row- and column permutations.

\smallskip

We will return to strictly ultrametric and related matrices in Section~\ref{ESect-6}, because they will turn out to be 
a  powerful tool to prove our main results. See also \cite{BrCi} for a detailed account on the geometric properties 
of the special type of simplices whose Gramians are ultrametric.
\subsection{Main results obtained from analyzing the generated data}\label{sect-E1.2}
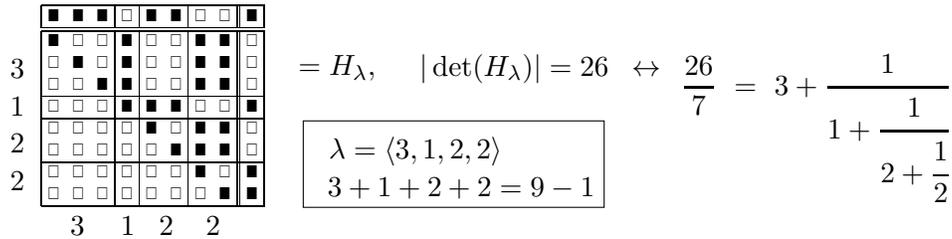
\begin{figure}[h]
\begin{center}  
\begin{tikzpicture}
\node[xscale=0.5, yscale=0.6] at (0,0) {$\begin{array}{|ccc|c|cc|cc||c|}
\hline
\blacksquare&\blacksquare&\blacksquare&\square&\blacksquare&\blacksquare&\square&\square&\blacksquare\\\hline\hline
\blacksquare&\square&\square&\blacksquare&\square&\square&\blacksquare&\blacksquare&\square\\
\square&\blacksquare&\square&\blacksquare&\square&\square&\blacksquare&\blacksquare&\square\\
\square&\square&\blacksquare&\blacksquare&\square&\square&\blacksquare&\blacksquare&\square\\\hline
\square&\square&\square&\blacksquare&\blacksquare&\blacksquare&\square&\square&\blacksquare\\\hline
\square&\square&\square&\square&\blacksquare&\square&\blacksquare&\blacksquare&\square\\
\square&\square&\square&\square&\square&\blacksquare&\blacksquare&\blacksquare&\square\\\hline
\square&\square&\square&\square&\square&\square&\blacksquare&\square&\blacksquare\\
\square&\square&\square&\square&\square&\square&\square&\blacksquare&\blacksquare\\
\hline\end{array}$};
\node at (3.5,-0.6) {$\lambda=\langle 3,1,2,2\rangle$};
\node at (4.1,-1.1) {$3+1+2+2=9-1$};
\draw (2,-1.35)--(6,-1.35)--(6,-0.2)--(2,-0.2)--cycle;
\node at (2.45,0.5) {$=H_\lambda,$};
\node at (5.15,0.5) {$|\det(H_\lambda)| = 26\,\,\,\leftrightarrow$};
\node at (8.8,-0.252)  {$\dfrac{26}{7} \,\,=\,\, 3 + \cfrac{1}{1 + \cfrac{1}{2 + \cfrac{1}{2} } }$};
\node at (-1.8,0.5) {$3$};
\node at (-1.8,0) {$1$};
\node at (-1.8,-0.5) {$2$};
\node at (-1.8,-1) {$2$};
\node at (-1,-1.6) {$3$};
\node at (-0.33,-1.6) {$1$};
\node at (0.18,-1.6) {$2$};
\node at (0.8,-1.6) {$2$};
\end{tikzpicture}
\end{center} 
\caption{\small{The matrix~$H_\lambda$ for a composition~$\lambda=\langle\lambda_1,\dots,\lambda_k\rangle$ of~$n-1$ 
with first and last parts larger than one, and its determinant as numerator of 
the continued fraction~$[\lambda_1;\lambda_2,\dots,\lambda_k]$.}}
\label{figureE3}
\end{figure}

A positive result in this context is as follows. Let~$n\geq 3$, and let the ordered tupel 
$\lambda=\langle \lambda_1,\dots,\lambda_k\rangle$ be a {\em composition} of the integer~$n-1$ whose first 
and last part are at least~$2$. Associate with~$\lambda$ the~$n\times n$ matrix~$H_\lambda$ as is done for 
the example~$\lambda=\langle 3,1,2,2\rangle$ in Figure~\ref{figureE3}. 

\smallskip

The matrix~$H_\lambda$ is constructed as an {\em unreduced upper Hessenberg} \index{unreduced upper Hessenberg} matrix with identity matrices 
$I_j$ of size~$\lambda_j\times\lambda_j$ covering the lower co-diagonal from top left to bottom right. The 
matrices~$I_1,\dots,I_k$ define a {\em checkerboard pattern} in~$H_\lambda$ above~$I_1,\dots,I_k$, with 
blocks containing either only ones, or only zero entries, where the blocks directly bordering~$I_j$ and~$I_{j+1}$ 
contain only ones. This uniquely defines~$H_\lambda$ in terms of the composition~$\lambda$.

\smallskip

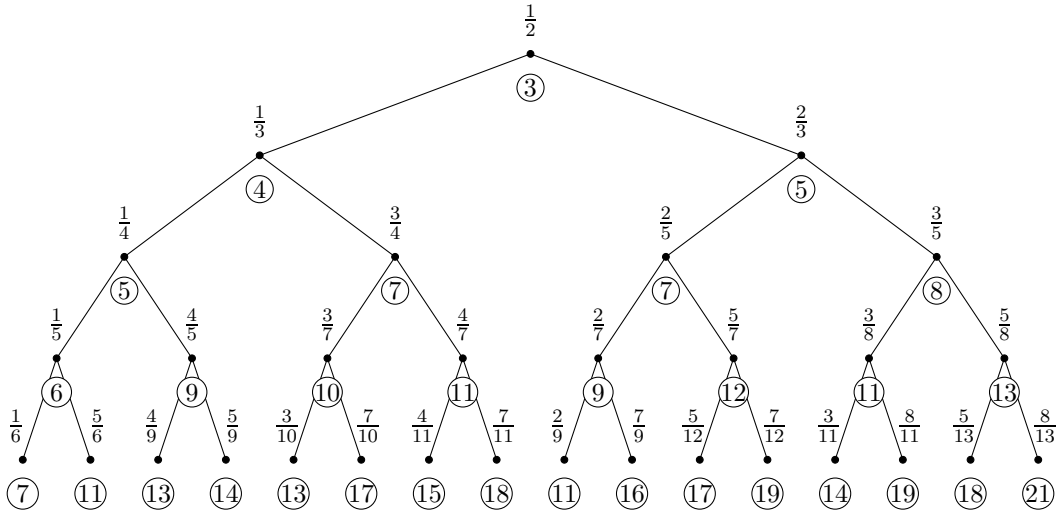
\begin{figure}[h]
\begin{center}  
\begin{tikzpicture}[scale=0.9, every node/.style={scale=0.9}]
\draw[fill=black] (0,0) circle [radius=0.05];
\draw[fill=black] (1,0) circle [radius=0.05];
\draw[fill=black] (2,0) circle [radius=0.05];
\draw[fill=black] (3,0) circle [radius=0.05];
\draw[fill=black] (4,0) circle [radius=0.05];
\draw[fill=black] (5,0) circle [radius=0.05];
\draw[fill=black] (6,0) circle [radius=0.05];
\draw[fill=black] (7,0) circle [radius=0.05];
\draw[fill=black] (8,0) circle [radius=0.05];
\draw[fill=black] (9,0) circle [radius=0.05];
\draw[fill=black] (10,0) circle [radius=0.05];
\draw[fill=black] (11,0) circle [radius=0.05];
\draw[fill=black] (12,0) circle [radius=0.05];
\draw[fill=black] (13,0) circle [radius=0.05];
\draw[fill=black] (14,0) circle [radius=0.05];
\draw[fill=black] (15,0) circle [radius=0.05];
\draw[fill=black] (0.5,1.5) circle [radius=0.05];
\draw[fill=black] (2.5,1.5) circle [radius=0.05];
\draw[fill=black] (4.5,1.5) circle [radius=0.05];
\draw[fill=black] (6.5,1.5) circle [radius=0.05];
\draw[fill=black] (8.5,1.5) circle [radius=0.05];
\draw[fill=black] (10.5,1.5) circle [radius=0.05];
\draw[fill=black] (12.5,1.5) circle [radius=0.05];
\draw[fill=black] (14.5,1.5) circle [radius=0.05];
\draw[fill=black] (1.5,3) circle [radius=0.05];
\draw[fill=black] (5.5,3) circle [radius=0.05];
\draw[fill=black] (9.5,3) circle [radius=0.05];
\draw[fill=black] (13.5,3) circle [radius=0.05];
\draw[fill=black] (3.5,4.5) circle [radius=0.05];
\draw[fill=black] (11.5,4.5) circle [radius=0.05];
\draw[fill=black] (7.5,6) circle [radius=0.05];
\draw (7.5,6)--(3.5,4.5);
\draw (7.5,6)--(11.5,4.5);
\node at (7.5,6.5) {$\frac{1}{2}$};
\draw (3.5,4.5)--(1.5,3);
\draw (11.5,4.5)--(9.5,3);
\node at (3.5,5) {$\frac{1}{3}$};
\node at (11.5,5) {$\frac{2}{3}$};
\draw (3.5,4.5)--(5.5,3) ;
\draw (11.5,4.5)--(13.5,3);
\draw (1.5,3)--(0.5,1.5);
\draw (1.5,3)--(2.5,1.5);
\draw (9.5,3)--(8.5,1.5);
\draw (9.5,3)--(10.5,1.5);
\draw (5.5,3)--(4.5,1.5);
\draw (5.5,3)--(6.5,1.5);
\draw (13.5,3)--(12.5,1.5);
\draw (13.5,3)--(14.5,1.5);
\node at (1.5,3.5) {$\frac{1}{4}$};
\node at (5.5,3.5) {$\frac{3}{4}$};
\node at (9.5,3.5) {$\frac{2}{5}$};
\node at (13.5,3.5) {$\frac{3}{5}$};
\draw (0.5,1.5)--(0,0);
\draw (0.5,1.5)--(1,0);
\draw (2.5,1.5)--(2,0);
\draw (2.5,1.5)--(3,0);
\draw (4.5,1.5)--(4,0);
\draw (4.5,1.5)--(5,0);
\draw (6.5,1.5)--(6,0);
\draw (6.5,1.5)--(7,0);
\draw (8.5,1.5)--(8,0);
\draw (8.5,1.5)--(9,0);
\draw (10.5,1.5)--(10,0);
\draw (10.5,1.5)--(11,0);
\draw (12.5,1.5)--(12,0);
\draw (12.5,1.5)--(13,0);
\draw (14.5,1.5)--(14,0);
\draw (14.5,1.5)--(15,0);
\node at (0.5,2) {$\frac{1}{5}$};
\node at (2.5,2) {$\frac{4}{5}$};
\node at (4.5,2) {$\frac{3}{7}$};
\node at (6.5,2) {$\frac{4}{7}$};
\node at (8.5,2) {$\frac{2}{7}$};
\node at (10.5,2) {$\frac{5}{7}$};
\node at (12.5,2) {$\frac{3}{8}$};
\node at (14.5,2) {$\frac{5}{8}$};
\node at (-0.1,0.5) {$\frac{1}{6}$};
\node at (1.1,0.5) {$\frac{5}{6}$};
\node at (1.9,0.5) {$\frac{4}{9}$};
\node at (3.1,0.5) {$\frac{5}{9}$};
\node at (3.9,0.5) {$\frac{3}{10}$};
\node at (5.1,0.5) {$\frac{7}{10}$};
\node at (5.9,0.5) {$\frac{4}{11}$};
\node at (7.1,0.5) {$\frac{7}{11}$};
\node at (7.9,0.5) {$\frac{2}{9}$};
\node at (9.1,0.5) {$\frac{7}{9}$};
\node at (9.9,0.5) {$\frac{5}{12}$};
\node at (11.1,0.5) {$\frac{7}{12}$};
\node at (11.9,0.5) {$\frac{3}{11}$};
\node at (13.1,0.5) {$\frac{8}{11}$};
\node at (13.9,0.5) {$\frac{5}{13}$};
\node at (15.1,0.5) {$\frac{8}{13}$};
\node at (7.5,5.5) {$3$};
\draw (7.5,5.5) circle [radius=0.2];
\node at (3.5,4) {$4$};
\node at (11.5,4) {$5$};
\draw (3.5,4) circle [radius=0.2];
\draw (11.5,4) circle [radius=0.2];
\node at (1.5,2.5) {$5$};
\node at (5.5,2.5) {$7$};
\node at (9.5,2.5) {$7$};
\node at (13.5,2.5) {$8$};
\draw (1.5,2.5) circle [radius=0.2];
\draw (5.5,2.5) circle [radius=0.2];
\draw (9.5,2.5) circle [radius=0.2];
\draw (13.5,2.5) circle [radius=0.2];
\draw[fill=white] (0.5,1) circle [radius=0.22];
\draw[fill=white] (2.5,1) circle [radius=0.22];
\draw[fill=white] (4.5,1) circle [radius=0.22];
\draw[fill=white] (6.5,1) circle [radius=0.22];
\draw[fill=white] (8.5,1) circle [radius=0.22];
\draw[fill=white] (10.5,1) circle [radius=0.22];
\draw[fill=white] (12.5,1) circle [radius=0.22];
\draw[fill=white] (14.5,1) circle [radius=0.22];
\node at (0.5,1) {$6$};
\node at (2.5,1) {$9$};
\node at (4.5,1) {$10$};
\node at (6.5,1) {$11$};
\node at (8.5,1) {$9$};
\node at (10.5,1) {$12$};
\node at (12.5,1) {$11$};
\node at (14.5,1) {$13$};
\node at (0,-0.5) {$7$};
\node at (1,-0.5) {$11$};
\node at (2,-0.5) {$13$};
\node at (3,-0.5) {$14$};
\node at (4,-0.5) {$13$};
\node at (5,-0.5) {$17$};
\node at (6,-0.5) {$15$};
\node at (7,-0.5) {$18$};
\node at (8,-0.5) {$11$};
\node at (9,-0.5) {$16$};
\node at (10,-0.5) {$17$};
\node at (11,-0.5) {$19$};
\node at (12,-0.5) {$14$};
\node at (13,-0.5) {$19$};
\node at (14,-0.5) {$18$};
\node at (15,-0.5) {$21$};
\draw (0,-0.5) circle [radius=0.23];
\draw (1,-0.5) circle [radius=0.23];
\draw (2,-0.5) circle [radius=0.23];
\draw (3,-0.5) circle [radius=0.23];
\draw (4,-0.5) circle [radius=0.23];
\draw (5,-0.5) circle [radius=0.23];
\draw (6,-0.5) circle [radius=0.23];
\draw (7,-0.5) circle [radius=0.23];
\draw (8,-0.5) circle [radius=0.23];
\draw (9,-0.5) circle [radius=0.23];
\draw (10,-0.5) circle [radius=0.23];
\draw (11,-0.5) circle [radius=0.23];
\draw (12,-0.5) circle [radius=0.23];
\draw (13,-0.5) circle [radius=0.23];
\draw (14,-0.5) circle [radius=0.23];
\draw (15,-0.5) circle [radius=0.23];
\end{tikzpicture}
\end{center} 
\caption{\small{Part of Kepler's {\em Tree of Fractions} and absolute determinants of the matrices~$H_\lambda$.}}
\label{figureE4}
\end{figure}

In Section~\ref{Esect-7.2} we will prove the following results in this context.
\begin{Th}\label{Eth-main} Let~$H$ be an~$n\times n$ unreduced upper Hessenberg 0/1-matrix whose columns 
and the origin are the~$n+1$ vertices of a simplex~$S\subset I^n$ with acute dihedral angles only. Then, 
possibly after exchanging its first two rows and/or last two columns,~$H$ is equal to the matrix~$H_\lambda$ 
for some composition~$\lambda=\langle\lambda_1,\dots,\lambda_k\rangle$ of~$n-1$ with first and last parts 
larger than one. Moreover, 
\be  |\det(H_\lambda)| = f_k, \hdrie \mbox{\rm where }\hdrie \frac{f_k}{g_k} = \lambda_1 + \cfrac{1}{\lambda_2 + \cfrac{1}{\ddots + \cfrac{1}{\lambda_k} } }, \hdrie \gcd(f_k,g_k)=1.\ee
Conversely, each such matrix~$H_\lambda$ has the property that its Gramian is strictly ultrametric, which implies 
that its columns together with the origin are the vertices of an acute 0/1-simplex.
\end{Th}
As a corollary of this theorem, all attainable absolute values of the determinant function on the set of all 
unreduced~$n\times n$ upper Hessenberg 0/1-matrices~$H$ for which~$(H^\top H)^{-1}$ is a diagonally 
dominant Stieltjes matrix with negative off-diagonal entries, can be explicitly read from a part of 
{\em Kepler's Tree of Fractions} \cite{KimMo}. This part is depicted in Figure~\ref{figureE4}. It has the fraction~$\frac{1}{2}$ 
as root. The children of a vertex~$\frac{p}{q}$ are~$\frac{p}{p+q}$ and~$\frac{q}{p+q}$. Transversing the 
tree level by level corresponds to an enumeration of all the rationals~$\mathbb{Q}\cap(0,1)$. 
\\
The circled integers displayed in Figure~\ref{figureE4} below each vertex equal the sum of numerator and denominator of 
the fraction belonging to that vertex. At level~$k$ these integers correspond to the absolute values of 
the determinants of each of the~$2^k$ matrices~$H_\lambda$ of size~$(k+4)\times (k+4)$.

\smallskip

Observe that the determinants in the rightmost branch in the tree equal the Fibonacci numbers, which were 
proved in \cite{Ching} to be the maximal value of the determinant function over {\em all}~$n\times n$ 
Hessenberg 0/1-matrices. We can now conclude that this maximum is (also) attained by matrices 
representing {\em acute} simplices. More generally we see that any branch of the tree that, starting at a 
given vertex~$p/q$ corresponding to the determinantal value~$p+q$, extends only to the right, yields a 
Fibonacci-type sequence~$d_r(j)$,
\be d_r(j+2) = d_r(j)+d_r(j+1) \hdrie\mbox{\rm with} \hdrie d_r(-1)=p,\hdrie d_r(0)=q \und d_r(1)=p+q,\ee
whereas any branch from a vertex~$p/q$ that extends only to the left, yields a family of acute 0/1-simplices 
with determinants increasing linearly as 
\be d_\ell(j)=jp+q. \ee
The corresponding matrices~$H_\lambda$ in this latter case have integer compositions of which the last 
part increases by one when the size of~$H_\lambda$ increases by one while all the other parts of~$\lambda$ 
remain the same. The existence of such families with linearly increasing determinants was first observed 
in \cite{BrHoKuSt}. In Section~\ref{ESect-6} we give a full explanation of their structure.
 %%%%%%
 %%%%%
 %%%%%
 %%%%%%%%%
\subsection{Outline}  
Our aim is to give a self-contained account of all necessary ingredients. For this, we first recall in 
Section~\ref{ESect-1} the group structure of~$\Bn$ and the permutation subgroup of~$S_{2^n}$ it induces on the 
set~$\BB^n$ of 0/1-vectors of length~$n$. These induced permutations were studied by Harrison and High, 
who derived a formula in \cite{HaHi} for the corresponding {\em cycle index polynomial}~$Z_n$ of~$\Bn$. 
This formula was later claimed to be simplified by Chen in \cite{Che}, who also studied the induced permutations 
of the {\em edges} of~$I^n$. Unfortunately, in view of the standard counting paradigm of P\'olya \cite{Polya}, 
neither formula allows a straightforward evaluation, modulo the induced action of~$\Bn$ on~$\BB^n$, of the 
number~$\vep_n^k$ of 0/1-polytopes in~$I^n$ with~$k$ vertices. Therefore, here we will aim for a more 
pragmatic approach, also motivated  by the fact that for~$n>6$ we failed to find explicit expressions for~$Z_n$ 
in the literature. First, in Section~\ref{ESect-1}, we give transparent algorithmic descriptions of how to compute 
$Z_n$ by means of a simple computer code. This code yields~$Z_n$ as a table of coefficients and exponents 
of monomials in a minimal effort: the table for~$n=9$ in Section~\ref{ESect-7} was, for instance, produced on a 
simple laptop within half a second. As a next step, in Section~\ref{ESect-2} we explain how to compute, modulo 
the action of~$\Bn$, the numbers~$\vep_n^k$ of 0/1-polytopes with~$k$ vertices by applying P\'olya's theory 
to the specific situation at hand. Again, the emphasis is to show how to algorithmically obtain the 
{\em concrete values} of~$\vep_n^k$ by means of a computer code, using the tables for the cycle indices of 
$\Bn$. As we will be interested in 0/1-simplices, we pay special attention to the values~$k\leq n+1$. 
In Section~\ref{ESect-7} we present a selection of the numbers produced by the algorithms.

\smallskip

In Section~\ref{ESect-3} we change our perspective from 0/1-polytopes and two-colorings to 
0/1-{\em matrices}. A 0/1-polytope~$c$ with~$k$ vertices can trivially be represented by a 0/1-matrix of 
size~$n\times k$ whose columns are the vertices of~$c$. Although convenient, this unfortunately introduces 
another non-trivial redundancy, as there are~$k!$ matrices having this property. Consequently, we investigate 
how to establish whether two given 0/1-matrices represent 0/1-polytopes in the same orbit under~$\Bn$. 
From all 0/1-matrices representing all the 0/1-polytopes in the same orbit under~$\Bn$, we select one 
designated matrix, the {\em minimal matrix representation}~$P^\ast$, and study its properties. As a first 
application, this concept enables us to enumerate all 0/1-triangles in~$I^n$ modulo the action of~$\Bn$: we 
give the minimal matrix representation of each of the~$\vep_n^2$ distinct orbits of 0/1-triangles under 
$\Bn$ using~$\mathcal{O}(1)$ arithmetic operations per triangle. We do the same for the subset of {\em acute} 
0/1-triangles. Basically, we parametrize both sets with the points with integer coordinates in a three-dimensional 
polyhedron, which in both cases turns out to be a simple tetrahedron. We also derive an explicit formula for their 
cardinalities by counting the integer points in the respective tetrahedra. In theory, the same can be done for 
$k$-simplices. This however results in enumerating and counting the points with integer coordinates in a polytope 
of dimension~$2^k-1$ constrained by at most~$(k+1)!$ inequalities. Although the enumeration would still cost 
$\mathcal{O}(1)$ per~$k$-simplex independent of~$n$, the dependence on~$k$ makes such enumeration impractical.

\smallskip
 
As a consequence of this intractability, in Section~\ref{ESect-5} we use the assistance of the computer to extend 
the minimal matrix representations of the acute 0/1-{\em triangles} from Section~\ref{ESect-3} into minimal 
matrix representations of  acute 0/1-{\em tetrahedra}, and similarly further into minimal matrix representations 
of all acute 0/1-{\em simplices} with~$n+1$ vertices. Since acute simplices have acute {\em facets} \cite{Fie}, 
each minimal matrix representation of an acute 0/1-tetrahedron equals a minimal matrix representation~$T$ of 
an acute 0/1-triangle {\em with one additional column~$t$ appended}. Hence, in theory, one could append one 
by one all feasible columns~$t$ to~$T$ such that~$[T|t]$ represents an acute tetrahedron, discard the ones that 
do not yield a minimal matrix representation, and continue to add more columns. Unfortunately, the verification 
of {\em minimality} is computationally much more expensive than verifying {\em acuteness}. It may thus be 
much quicker to find out if~$[T|t]$ can be extended to the desired number of columns, then to find out if it is 
minimal. This saves computational effort if it {\em cannot} be acutely extended. 
\begin{table}[h]
\begin{center}
$\begin{array}{|r||r|r|r|r|r|r|r|r|r|r|r|r|r|r|}
\hline
n &1 & 2 & 3 & 4 & 5 & 6 & 7 & 8 & 9 & 10 & 11 \\
 \hline\hline
a(n) &1 & 0 & 1 & 1 & 2 & 6 & 13 & 29 & 67 & 162 & 392\\ 
  \hline
s(n) &1 & 1 & 6 & 27 & 472 & 19735 & 2773763 & 1245930065 & 1.8e12 & 8.7e15& 1.3e20\\
 \hline
\end{array}$
\end{center}
\caption{\small{The number~$a(n)$ of {\em acute} 0/1~$n$-simplices in~$I^n$ related to their 
{\em total} number~$s(n)$. All cardinalities are modulo the action of~$\Bn$.}}
\label{Etable1}
\end{table}\\
According to the data from Table~\ref{Etable1} (see also Section~\ref{ESect-7}), acute~$n$-simplices in~$I^n$ are extremely rare. 
Thus it seems likely that simply extending a minimal matrix representation until no acute extensions are available 
anymore is quicker than discarding the matrix representations that are not minimal. However, it turns out that 
the amount of data in the intermediate phases becomes unacceptably large. Thus, the challenge to make 
our algorithms as efficiently as possible is therefore nontrivial, and involves the well-known struggle between 
time and memory requirements. It requires a subtle balance between spending time in computing minimal matrix 
representations, and allowing the data to take more and more memory space. Along the way, and also for 
the purpose of mathematical analysis, we introduce the sets of {\em candidate acute extensions}~$\CC^n(S)$ 
and of {\em acute extensions}~$\AA^n(S)$ of a given acute simplex~$S\subset I^n$. Using the theory of symmetric 
inverse M-matrices (also called Stieltjes matrices) \cite{Joh,JoSm}, it is possible to derive relations between 
the members of these sets that make their computation in many cases much less expensive than at first sight. 
In Section~\ref{Esect-6.2} we display the minimal matrix representations of all the acute 0/1-simplices in~$I^n$ 
for~$n\in\{3,4,\dots,9\}$.

\smallskip

In Section~\ref{ESect-6} we analyze these minimal matrix representations. Proofs of the results given in 
Section~\ref{sect-E1.2} will be based on the so-called {\em One Neighbor Theorem} for acute 0/1-simplices 
\cite{BrDiHaKr}. This theorem states that the set~$\CC^n(S)$ of candidate acute extentions of a simplex 
$S\subset I^n$ with~$n$ vertices consists of at most two antipodal points. As a consequence, only points in higher 
dimensional cubes that project orthogonally on this antipodal pair can be added to form acute 0/1-simplices with 
 more vertices. If one demands their matrix representation to be unreduced upper Hessenberg, this limits the 
 possible options even further. What results is a complete description of the corresponding simplices in 
 Section~\ref{Esect-7.2} together with the values of the determinants of their matrix representations in terms of 
 continued fractions. Apart from the One Neighbor Theorem, we also use the properties of {\em strictly ultrametric 
 matrices} \cite{NabVar} to prove acuteness of the simplices involved.
 
\section{The hyperoctahedral group~$\Bn$}\label{ESect-1}
Write~$I^n=[0,1]^n$ for the unit~$n$-cube and~$\BB^n=\{0,1\}^n$ for the set of its vertices. Let~$\Bn$ be the set 
of all rigid transformations~$h:I^n\rightarrow I^n$. Endowing~$\Bn$ with the usual composition of map as 
multiplication rule, it becomes the {\em hyperoctahedral group} of~$n$-cube symmetries, with the {\em dihedral} 
group~$\mathcal{B}_2$ and the {\em octahedral} group~$\mathcal{B}_3$ as well-known instances. Each~$h\in\Bn$ 
bijectively maps~$k$-facets of~$I^n$ to~$k$-facets and thus induces a permutation of these~$k$-facets; in 
particular, it permutes~$\BB^n$. To describe this latter permutation, we choose the following bijection~$\beta$ as 
{\em numbering} of~$\BB^n$. It interprets the 0/1-vector~$v\in\BB^n$ as a binary number. 
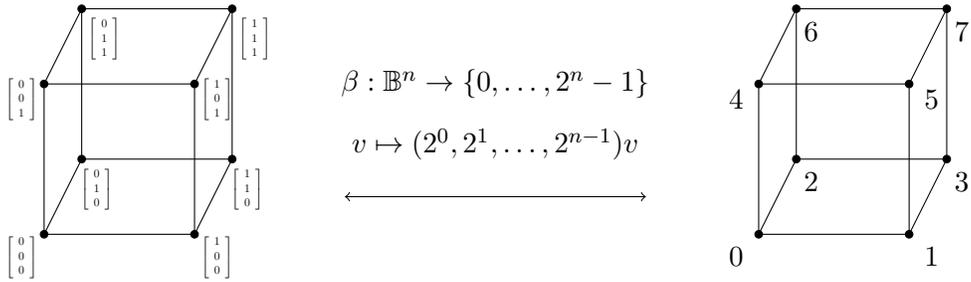
\begin{figure}[h]
\begin{center}  
\begin{tikzpicture}[scale=1]
\begin{scope}
\draw (0,0)--(2,0)--(2,2)--(0,2)--cycle;
\draw (0.5,1)--(2.5,1)--(2.5,3)--(0.5,3)--cycle;
\draw (0,0)--(0.5,1);
\draw (2,0)--(2.5,1);
\draw (0,2)--(0.5,3);
\draw (2,2)--(2.5,3);
\draw[fill=black] (0,0) circle [radius=0.05];
\draw[fill=black] (0.5,1) circle [radius=0.05];
\draw[fill=black] (2.5,1) circle [radius=0.05];
\draw[fill=black] (0.5,3) circle [radius=0.05];
\draw[fill=black] (2,0) circle [radius=0.05];
\draw[fill=black] (0,2) circle [radius=0.05];
\draw[fill=black] (2,2) circle [radius=0.05];
\draw[fill=black] (2.5,3) circle [radius=0.05];
\node at (6,2) {$\beta: \mathbb{B}^n\rightarrow \{0,\dots,2^n-1\}$};
\node at (6,1.2) {$v \mapsto (2^0,2^1,\dots,2^{n-1}) v$};
\draw[<->] (4,0.5)--(8,0.5);
\node[scale=0.4] at (-0.3,-0.3) {$\left[\begin{array}{c}0\\0\\0\end{array}\right]$};
\node[scale=0.4] at (2.3,-0.3) {$\left[\begin{array}{c}1\\0\\0\end{array}\right]$};
\node[scale=0.4] at (0.7,0.6) {$\left[\begin{array}{c}0\\1\\0\end{array}\right]$};
\node[scale=0.4] at (2.7,0.6) {$\left[\begin{array}{c}1\\1\\0\end{array}\right]$};
\node[scale=0.4] at (-0.3,1.8) {$\left[\begin{array}{c}0\\0\\1\end{array}\right]$};
\node[scale=0.4] at (2.3,1.8) {$\left[\begin{array}{c}1\\0\\1\end{array}\right]$};
\node[scale=0.4] at (0.8,2.6) {$\left[\begin{array}{c}0\\1\\1\end{array}\right]$};
\node[scale=0.4] at (2.8,2.6) {$\left[\begin{array}{c}1\\1\\1\end{array}\right]$};
\end{scope}
\begin{scope}[shift={(9.5,0)}]
\draw (0,0)--(2,0)--(2,2)--(0,2)--cycle;
\draw (0.5,1)--(2.5,1)--(2.5,3)--(0.5,3)--cycle;
\draw (0,0)--(0.5,1);
\draw (2,0)--(2.5,1);
\draw (0,2)--(0.5,3);
\draw (2,2)--(2.5,3);
\draw[fill=black] (0,0) circle [radius=0.05];
\draw[fill=black] (0.5,1) circle [radius=0.05];
\draw[fill=black] (2.5,1) circle [radius=0.05];
\draw[fill=black] (0.5,3) circle [radius=0.05];
\draw[fill=black] (2,0) circle [radius=0.05];
\draw[fill=black] (0,2) circle [radius=0.05];
\draw[fill=black] (2,2) circle [radius=0.05];
\draw[fill=black] (2.5,3) circle [radius=0.05];
\node at (-0.3,-0.3) {$0$};
\node at (2.3,-0.3) {$1$};
\node at (0.7,0.7) {$2$};
\node at (2.7,0.7) {$3$};
\node at (-0.3,1.8) {$4$};
\node at (2.3,1.8) {$5$};
\node at (0.7,2.7) {$6$};
\node at (2.7,2.7) {$7$};
\end{scope}
\end{tikzpicture}
\end{center} 
\caption{\small{Binary-to-decimal numbering~$\beta$ of the vertex set~$\BB^n$ of~$I^n$, depicted for~$n=3$.}}
\label{figureE5}
\end{figure}

\smallskip

Each~$h\in\Bn$ induces a permutation~$\pi_h\in S_{2^n}$ of the numbers~$0,\dots,2^n-1$ via~$\beta$ by 
\be\label{Ep-1} \pi_h: \hdrie  \{0,\dots,2^n-1\} \rightarrow  \{0,\dots,2^n-1\}: \hdrie k\mapsto (\beta\circ h \circ \beta^{-1})(k). \ee
In Section~\ref{ESect-2} we will count the number of {\em 0/1-polytopes} modulo~$n$-cube symmetries using 
P\'olya's Theorem \cite{Polya}. For this we need to know {\em how many} permutations~$\pi_h$ 
{\em of which cycle type} are induced in~$S_{2^n}$ when~$h$ ranges over~$\Bn$.  Recalling that any permutation 
$\pi$ of~$\ell$ objects can be written as the product of {\em disjoint cycles}, we can define the {\em cycle type} of 
$\pi$.   
\begin{Def}[Cycle type]\label{Ecycle-type} {\rm If~$\pi\in S_\ell$ has~$t_i$ cycles of length~$i$ in its cycle \index{cycle type}
factorization, then the vector 
\be t(\pi) = (t_1,\dots,t_\ell), \hdrie\mbox{ with}  \hdrie t_1\cdot 1+\dots+t_\ell\cdot \ell = \ell, \ee 
is an integer partition of~$\ell$ called the {\em cycle type} of~$\pi$}. 
\end{Def}
From basic algebra we know that the cycle types of two permutations in~$S_\ell$ coincide if and only if they are 
{\em conjugate}.     
\begin{Le}\label{Elem-5} Let~$\pi,\tilde{\pi}\in S_\ell$, then~$ t(\tilde{\pi})=t(\pi)$ if and only if 
$\tilde{\pi} = \sigma^{-1} \circ \pi \circ \sigma$ for some~$\sigma\in S_\ell$. 
\end{Le}
An important consequence is the following corollary, whose formulation uses  (\ref{Ep-1}). 
\begin{Co}\label{ECo-0} If two elements~$h,\tilde{h}$ are conjugate in~$\Bn$ then~$t(\pi_h)=t(\pi_{\tilde{h}})$ in 
$S_{2^n}$. 
\end{Co} 
{\bf Proof. } The mapping~$\Bn\rightarrow S_{2^n}: f\mapsto \pi_f$ is an injective homomorphism. Therefore, if 
$\tilde h = g\circ h \circ g^{-1}$ then~$\pi_{\tilde{h}} = \pi_g\circ \pi_h \circ \pi_g^{-1}$. Lemma~\ref{Elem-5} now  
proves the statement. \hfill~$\Box$
\begin{rem}\label{Erem-1} {\rm The table in  (\ref{Ep-36}) constitutes an example of the fact that~$t(\pi_g)=t(\pi_h)$ 
while~$g$ and~$h$ are not conjugate in~$\Bn$, hence the converse implication in Corollary~\ref{ECo-0} does not hold}. 
\end{rem} 
Corollary~\ref{ECo-0} shows that counting how many permutations of which type are induced by the elements of 
$\Bn$ reduces to the following two tasks, \\[2mm]
$\bullet$ find the cycle type of~$\pi_h$ of a single element~$h$ from each conjugacy class of~$\Bn$, \\[2mm]
$\bullet$ count the number of elements in each conjugacy class of~$\Bn$.\\[2mm]
Before performing these tasks in Section~\ref{ESect-1.5}, we recall some basic facts about~$\Bn$. We identify two 
subgroups~$\Bn^c$ and~$\Bn^p$ of~$\Bn$ and show that~$\Bn = \Bn^p \times \Bn^c$. This enables us to associate 
with each~$h\in\Bn$ a so-called {\em signed permutation}. The corresponding {\em signed cycle type} of such a 
signed permutation will then be used to describe and count the conjugacy classes of~$\Bn$, and consequently, the 
number and cycle type of their induced permutations in~$S_{2^n}$.
 %%%%%%%%
 %%%%%%%%%
 %%%%%%%%%%
\subsection{The subgroups~$\Bn^c$ \& $\Bn^p$: complementations \& permutations}\label{ESect-1.1} \index{coordinate complementation} \index{coordinate permutation}
Let the~$n$-tupel~$\langle e_1,\dots, e_n\rangle$ be the standard basis for~$\RR^n$. For~$j\in\{1,\dots,n\}$, let 
$c_j: I^n\rightarrow I^n:x\mapsto e_j + x-2e_je_j^\top x$ be the reflection in the affine hyperplane~$2x_j=1$. 
The set~$\{c_1,\dots,c_n\}$ generates a subgroup~$\Bn^c$ of~$\Bn$. Note that~$c_i\circ c_j = c_j\circ c_i$ and 
$c_j^2=id$. Thus, the mapping
\be\label{Ep-2} \BB^n \rightarrow \Bn^c: \hdrie w \mapsto c_w =  c_1^{w_1}\circ c_2^{w_2} \circ \dots \circ c_n^{w_n} \ee   
is a bijection, showing that~$|\Bn^c|=2^n$. One can verify that~$c_w(v) = \mbox{\rm xor}(w,v) = (w+v)\mod 2$, 
where xor is the logical {\em exclusive or} operation performed entry-wise on the pair~$w,v\in\BB^n$.

\smallskip

Next, for each~$j\in\{2,\dots,n\}$, consider the reflection~$s_{j}:I^n\rightarrow I^n:x\mapsto x-(e_1-e_j)(e_1-e_j)^\top x$ 
in the hyperplane~$x_1=x_j$. The set~$\{s_2,\dots,s_n\}$ generates a subgroup~$\Bn^p$ of~$\Bn$. The action of~$s_{j}$ on 
$v\in\BB^n$ interchanges the first and~$j$th entry of~$v$. Since {\em each} permutation of~$n$ objects is a product of 
transpositions with the first object \cite{Ehr}, we conclude that~$|\Bn^p|=n!$. 

\smallskip

For each permutation~$u=[u(1),\dots,u(n)]\in S_n$, we write~$p_u$ for the element from~$\Bn^p$ defined by its action on~$\BB^n$ as
\be\label{Epuv} (p_u)(v) = v\circ u = (v_{u(1)},\dots,v_{u(n)})^\top. \ee 
 
\begin{Def}[Coordinate complementation and permutation] {\rm An element~$c_w\in\Bn^c$ will be called a 
{\em coordinate complementation} and an element~$p_u\in\Bn^p$ a {\em coordinate permutation}}.
\end{Def}
\begin{figure}[b]
\begin{center}  
\begin{tikzpicture}[scale=1, every node/.style={scale=1}]
\begin{scope}
\draw[fill=gray!20!white] (1,0)--(1.5,0.8)--(1.5,2.8)--(1,2)--cycle;
\draw (0,0)--(2,0)--(2,2)--(0,2)--cycle;
\draw (0.5,0.8)--(2.5,0.8)--(2.5,2.8)--(0.5,2.8)--cycle;
\draw (0,0)--(0.5,0.8);
\draw (2,0)--(2.5,0.8);
\draw (0,2)--(0.5,2.8);
\draw (2,2)--(2.5,2.8);
\node at (1,-0.3) {$c_1$};
\end{scope}
\begin{scope}[shift={(3,0)}]
\draw[fill=gray!20!white] (0.25,0.4)--(2.25,0.4)--(2.25,2.4)--(0.25,2.4)--cycle;
\draw (0,0)--(2,0)--(2,2)--(0,2)--cycle;
\draw (0.5,0.8)--(2.5,0.8)--(2.5,2.8)--(0.5,2.8)--cycle;
\draw (0,0)--(0.5,0.8);
\draw (2,0)--(2.5,0.8);
\draw (0,2)--(0.5,2.8);
\draw (2,2)--(2.5,2.8); 
\node at (1,-0.3) {$c_2$};
\end{scope}
\begin{scope}[shift={(6,0)}]
\draw[fill=gray!20!white] (0,1)--(2,1)--(2.5,1.8)--(0.5,1.8)--cycle;
\draw (0,0)--(2,0)--(2,2)--(0,2)--cycle;
\draw (0.5,0.8)--(2.5,0.8)--(2.5,2.8)--(0.5,2.8)--cycle;
\draw (0,0)--(0.5,0.8);
\draw (2,0)--(2.5,0.8);
\draw (0,2)--(0.5,2.8);
\draw (2,2)--(2.5,2.8); 
\node at (1,-0.3) {$c_3$};
\end{scope}
\begin{scope}[shift={(9,0)}]
\draw[fill=gray!20!white] (0,0)--(2.5,0.8)--(2.5,2.8)--(0,2)--cycle;
\draw (0,0)--(2,0)--(2,2)--(0,2)--cycle;
\draw (0.5,0.8)--(2.5,0.8)--(2.5,2.8)--(0.5,2.8)--cycle;
\draw (0,0)--(0.5,0.8);
\draw (2,0)--(2.5,0.8);
\draw (0,2)--(0.5,2.8);
\draw (2,2)--(2.5,2.8); 
\node at (1,-0.3) {$s_2$};
\end{scope}
\begin{scope}[shift={(12,0)}]
\draw[fill=gray!20!white] (0,0)--(2,0)--(2.5,2.8)--(0.5,2.8)--cycle;
\draw (0,0)--(2,0)--(2,2)--(0,2)--cycle;
\draw (0.5,0.8)--(2.5,0.8)--(2.5,2.8)--(0.5,2.8)--cycle;
\draw (0,0)--(0.5,0.8);
\draw (2,0)--(2.5,0.8);
\draw (0,2)--(0.5,2.8);
\draw (2,2)--(2.5,2.8); 
\node at (1,-0.3) {$s_3$};
\end{scope}
\end{tikzpicture}
\end{center} 
\caption{\small{Generators of~$\mathcal{B}_3$: complementations~$c_1,c_2,c_3$ and the permutations~$s_2$ and~$s_3$.}}
\label{figureE6}
\end{figure}
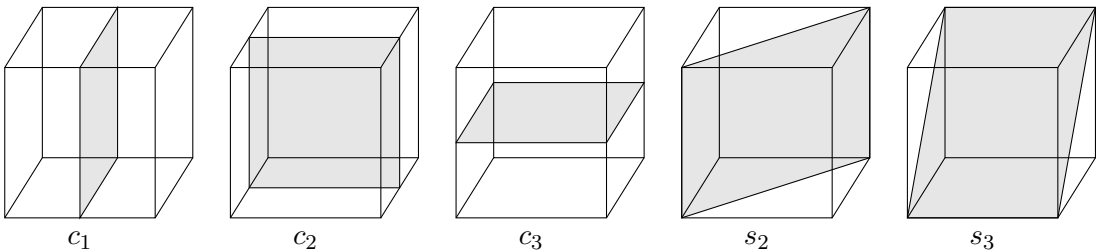
\begin{Ex} Consider the group~$\mathcal{B}_3$ of unit cube symmetries. It contains a subgroup 
$\mathcal{B}_3^c$ of order~$8=2^3$ with generators~$c_1,c_2,c_3$, the reflections in the planes~$2x_j=1$, 
and a subgroup $\mathcal{B}_3^p$ of order~$6=3!$ with generators~$s_2$ and~$s_3$, the reflectors in the planes~$x_1=x_2$ and~$x_1=x_3$. 

\smallskip

To illustrate the actions of elements from~$\mathcal{B}_3^c$ and~$\mathcal{B}_3^p$, let for instance 
$w=(0,1,1)^\top$, then~$c_w\in\Bn^c$ acts on~$I^3$ as depicted in the right part of Figure~\ref{figureE7}. Note that 
$c_w=c_1\circ c_2=c_2\circ c_1$. Also given is its induced permutation~$\pi_{c_w}\in S_8$. Next, given the 
permutation~$u=[3\,\, 1\,\, 2]\in S_3$, the action of~$p_u\in\Bn^p$ on~$I^3$ is depicted on the left, also 
together with its induced permutation~$\pi_{p_u}\in S_8$. Observe that~$p_u=s_3\circ s_2$, but that the 
product~$s_3\circ s_2$ does not equal~$p_u$.~\hfill~$\diamondsuit$
\begin{figure}[h]
\begin{center}
\begin{tikzpicture}[scale=1, every node/.style={scale=1}]
\begin{scope}
\draw (0,0)--(2,0)--(2,2)--(0,2)--cycle;
\draw (0.5,1)--(2.5,1)--(2.5,3)--(0.5,3)--cycle;
\draw (0,0)--(0.5,1);
\draw (2,0)--(2.5,1);
\draw (0,2)--(0.5,3);
\draw (2,2)--(2.5,3);
\draw[fill=black] (0,0) circle [radius=0.05];
\draw[fill=black] (0.5,1) circle [radius=0.05];
\draw[fill=black] (2.5,1) circle [radius=0.05];
\draw[fill=black] (0.5,3) circle [radius=0.05];
\draw[fill=black] (2,0) circle [radius=0.05];
\draw[fill=black] (0,2) circle [radius=0.05];
\draw[fill=black] (2,2) circle [radius=0.05];
\draw[fill=black] (2.5,3) circle [radius=0.05];
\node[scale=0.4] at (-0.3,-0.3) {$\left[\begin{array}{c}0\\0\\0\end{array}\right]$};
\node[scale=0.4] at (2.3,-0.3) {$\left[\begin{array}{c}0\\1\\0\end{array}\right]$};
\node[scale=0.4] at (0.7,0.6) {$\left[\begin{array}{c}0\\0\\1\end{array}\right]$};
\node[scale=0.4] at (2.7,0.6) {$\left[\begin{array}{c}0\\1\\1\end{array}\right]$};
\node[scale=0.4] at (-0.3,1.8) {$\left[\begin{array}{c}1\\0\\0\end{array}\right]$};
\node[scale=0.4] at (2.3,1.8) {$\left[\begin{array}{c}1\\1\\0\end{array}\right]$};
\node[scale=0.4] at (0.8,2.6) {$\left[\begin{array}{c}1\\0\\1\end{array}\right]$};
\node[scale=0.4] at (2.8,2.6) {$\left[\begin{array}{c}1\\1\\1\end{array}\right]$};
\node[scale=0.9] at (1,-1.5) {[0\,\,2\,\,4\,\,6\,\,1\,\,3\,\,5\,\,7]};
\end{scope}
\begin{scope}[shift={(6,0)}]
\draw (0,0)--(2,0)--(2,2)--(0,2)--cycle;
\draw (0.5,1)--(2.5,1)--(2.5,3)--(0.5,3)--cycle;
\draw (0,0)--(0.5,1);
\draw (2,0)--(2.5,1);
\draw (0,2)--(0.5,3);
\draw (2,2)--(2.5,3);
\draw[fill=black] (0,0) circle [radius=0.05];
\draw[fill=black] (0.5,1) circle [radius=0.05];
\draw[fill=black] (2.5,1) circle [radius=0.05];
\draw[fill=black] (0.5,3) circle [radius=0.05];
\draw[fill=black] (2,0) circle [radius=0.05];
\draw[fill=black] (0,2) circle [radius=0.05];
\draw[fill=black] (2,2) circle [radius=0.05];
\draw[fill=black] (2.5,3) circle [radius=0.05];
\node[scale=0.4] at (-0.3,-0.3) {$\left[\begin{array}{c}0\\0\\0\end{array}\right]$};
\node[scale=0.4] at (2.3,-0.3) {$\left[\begin{array}{c}1\\0\\0\end{array}\right]$};
\node[scale=0.4] at (0.7,0.6) {$\left[\begin{array}{c}0\\1\\0\end{array}\right]$};
\node[scale=0.4] at (2.7,0.6) {$\left[\begin{array}{c}1\\1\\0\end{array}\right]$};
\node[scale=0.4] at (-0.3,1.8) {$\left[\begin{array}{c}0\\0\\1\end{array}\right]$};
\node[scale=0.4] at (2.3,1.8) {$\left[\begin{array}{c}1\\0\\1\end{array}\right]$};
\node[scale=0.4] at (0.8,2.6) {$\left[\begin{array}{c}0\\1\\1\end{array}\right]$};
\node[scale=0.4] at (2.8,2.6) {$\left[\begin{array}{c}1\\1\\1\end{array}\right]$};
\draw[->] (3,1.5)--(5,1.5);
\node[scale=0.9] at (4,1.7) {$c_w$};
\draw[->] (-1,1.5)--(-3,1.5);
\node[scale=0.9] at (-2,1.7) {$p_u$};
\draw[->] (3,-1.5)--(5,-1.5);
\node[scale=0.9] at (4,-1.3) {$\pi_{c_w}$};
\draw[->] (-1,-1.5)--(-3,-1.5);
\node[scale=0.9] at (-2,-1.3) {$\pi_{p_u}$};
\node[scale=0.9] at (1,-1.5) {[0\,\,1\,\,2\,\,3\,\,4\,\,5\,\,6\,\,7]};
\begin{scope}[scale=0.7,shift={(4.5,-1.2)}]
\draw[fill=gray!20!white] (0.25,0.4)--(2.25,0.4)--(2.25,2.4)--(0.25,2.4)--cycle;
\draw[fill=gray!20!white] (0,1)--(2,1)--(2.5,1.8)--(0.5,1.8)--cycle;
\draw (0.25,0.4)--(2.25,0.4)--(2.25,2.4)--(0.25,2.4)--cycle;
\draw (0,0)--(2,0)--(2,2)--(0,2)--cycle;
\draw (0.5,0.8)--(2.5,0.8)--(2.5,2.8)--(0.5,2.8)--cycle;
\draw (0,0)--(0.5,0.8);
\draw (2,0)--(2.5,0.8);
\draw (0,2)--(0.5,2.8);
\draw (2,2)--(2.5,2.8);
\end{scope}
\begin{scope}[scale=0.7,shift={(-4,-1.2)}]
\draw[fill=gray!20!white] (0,0)--(2,0)--(2.5,2.8)--(0.5,2.8)--cycle;
\draw[fill=gray!20!white] (0,0)--(2.5,0.8)--(2.5,2.8)--(0,2)--cycle;
\draw (0,0)--(2,0)--(2.5,2.8)--(0.5,2.8)--cycle;
\draw (0,0)--(2,0)--(2,2)--(0,2)--cycle;
\draw (0.5,0.8)--(2.5,0.8)--(2.5,2.8)--(0.5,2.8)--cycle;
\draw (0,0)--(0.5,0.8);
\draw (2,0)--(2.5,0.8);
\draw (0,2)--(0.5,2.8);
\draw (2,2)--(2.5,2.8); 
\end{scope}
\end{scope}
\begin{scope}[shift={(12,0)}]
\draw (0,0)--(2,0)--(2,2)--(0,2)--cycle;
\draw (0.5,1)--(2.5,1)--(2.5,3)--(0.5,3)--cycle;
\draw (0,0)--(0.5,1);
\draw (2,0)--(2.5,1);
\draw (0,2)--(0.5,3);
\draw (2,2)--(2.5,3);
\draw[fill=black] (0,0) circle [radius=0.05];
\draw[fill=black] (0.5,1) circle [radius=0.05];
\draw[fill=black] (2.5,1) circle [radius=0.05];
\draw[fill=black] (0.5,3) circle [radius=0.05];
\draw[fill=black] (2,0) circle [radius=0.05];
\draw[fill=black] (0,2) circle [radius=0.05];
\draw[fill=black] (2,2) circle [radius=0.05];
\draw[fill=black] (2.5,3) circle [radius=0.05];
\node[scale=0.4] at (-0.3,-0.3) {$\left[\begin{array}{c}0\\1\\1\end{array}\right]$};
\node[scale=0.4] at (2.3,-0.3) {$\left[\begin{array}{c}1\\1\\1\end{array}\right]$};
\node[scale=0.4] at (0.7,0.6) {$\left[\begin{array}{c}0\\0\\1\end{array}\right]$};
\node[scale=0.4] at (2.7,0.6) {$\left[\begin{array}{c}1\\0\\1\end{array}\right]$};
\node[scale=0.4] at (-0.3,1.8) {$\left[\begin{array}{c}0\\1\\0\end{array}\right]$};
\node[scale=0.4] at (2.3,1.8) {$\left[\begin{array}{c}1\\1\\0\end{array}\right]$};
\node[scale=0.4] at (0.8,2.6) {$\left[\begin{array}{c}0\\0\\0\end{array}\right]$};
\node[scale=0.4] at (2.8,2.6) {$\left[\begin{array}{c}1\\0\\0\end{array}\right]$};
\node[scale=0.9] at (1,-1.5) {[6\,\,7\,\,4\,\,5\,\,2\,\,3\,\,0\,\,1]};
\end{scope}
\end{tikzpicture}
\end{center} 
\caption{\small{A complementation~$c_w$, a permutation~$p_u$, and their induced permutations in~$S_8$.}}
\label{figureE7}
\end{figure}
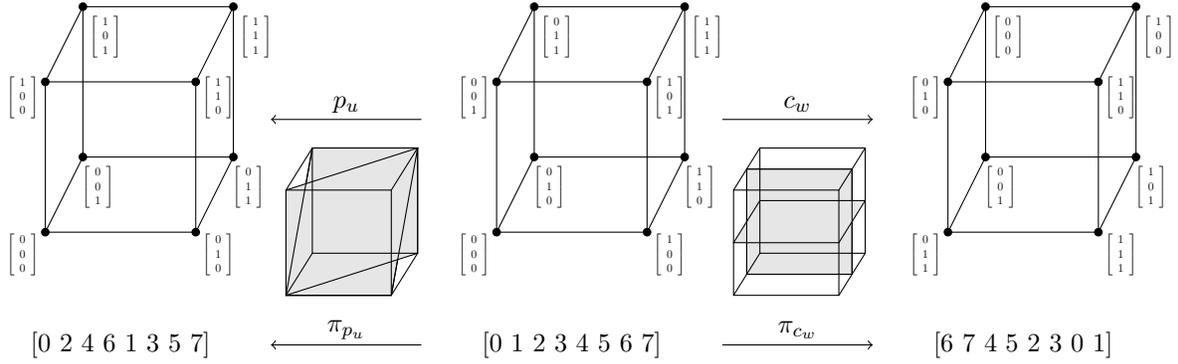
\end{Ex}
%%%%%%%%
%%%%%%%%%%
%%%%%%%%
%%%%%%%%%
\subsection{$\Bn$ and the group of signed permutations of~$n$ objects}\label{ESect-1.3}
An~$n$-cube symmetry~$h\in\Bn$ is a rigid transformation and thus an affine isometry. As 
such, it is uniquely determined by the combination of both the items (1) and (2), being\\[2mm]
$(1)$ the vertex~$v\in\BB^n$ that is mapped to the origin by~$h$,\\[2mm]
$(2)$ how the~$n$ vertices of~$I^n$ at distance one from~$v$ are mapped to~$e_1,\dots,e_n$.\\[2mm]
Note that~$c_v$ is the unique element from~$\Bn^c$ with~$c_v(v)=0$. Also note that~$p_u(0)=0$ 
and~$p(e)=e$ for all~$p_u\in\Bn^p$, where~$e=e_1+\dots+e_n$ is the {\em all-ones vector}. Thus, we have that
\be \left\{ h\in\Bn\sth h(v)=0\right\} \hdrie = \hdrie\left\{ h = p_u\circ c_v \sth p_u\in\Bn^p\right\}.\ee
Also observe that each~$p_u\in\Bn^p$ corresponds to a unique permutation of the basis vectors.
\begin{Co} For each~$h\in\Bn$ there exist unique pair~$p_u\in\Bn^p$ and~$c_v\in\Bn^c$ such that
\be\label{Ecomm} p_u\circ c_v = h =  c_{p_u(v)}\circ p_u.\ee
\end{Co}
One of the consequences of the uniqueness is that the order~$|\Bn|$ of~$\Bn$ equals~$n!2^n$. Another 
consequence is that we can now identify with each~$h\in\Bn$ a so-called {\em signed permutation}.
\begin{Def}[Signed permutation] We will index~$h\in\Bn$ as~$h_w$, where the vector~$w$, called \index{signed permutation} 
a {\em signed permutation}, has entries given by
\be w_j = u_j \hdrie\mbox{\rm if~$v_j=0$} \und w_j = \overline{u_j} \hdrie\mbox{\rm if~$v_j=1$}, \ee
where~$u$ and~$v$ are the indexes of the unique~$p_u\in\Bn^p$ and~$c_v\in\Bn^c$ such that~$h=p_u\circ c_v$.
\end{Def} 
The set of all signed permutations of~$n$ objects is obviously isomorphic to~$\Bn$. 
 %%%%%%%%%%
 %%%%%%%%%%
 %%%%%%%%%%%%%
 %%%%%%%%%%%%
\subsection{Conjugacy classes and signed cycle types in~$\Bn$}\label{ESect-1.5}
We now introduce the {\em signed cycle type} of a signed permutation. It will have 
the property that two elements in~$\Bn$ are conjugate if and only if they have the same signed cycle type. 
\begin{Def}[Signed cycle type] \index{signed cycle type}Let~$h=p_u\circ c_v$ with~$p_u\in\Bn^p$ and~$c_v\in\Bn^c$. For each cycle 
$\gamma$ in the decomposition of the permutation~$u$ into disjoint cycles, set
\be \chi(\gamma) = \sum_{j\in\gamma} v_j. \ee
Let~$u_+$ be the product of the cycles~$\gamma$ of~$u$ for which~$\chi(\gamma)$ is even, and~$u_-$ such that~$u=u_+\circ u_-$. Then the~$2\times n$ array
\be t_\pm(h) = \left\{ \begin{array}{cc} t(u_+)\\ t(u_-)\end{array}, \right. \ee
where~$t$ is the cycle type from Definition~\ref{Ecycle-type}, is called the {\em signed cycle type} of~$h$.
\end{Def}
Note that the signed cycle types of the elements~$h\in\Bn$ are in one-to-one correspondence with the 
{\em double partitions} \cite{GeKi} of~$n$, which are ordered pairs of partitions of~$k$ and~$\ell$ with~$k+\ell=n$.
\begin{Ex} Let~$h_w\in\Bc_{10}$ be the ten-cube symmetry indexed by the {\em signed permutation}
\be w=[2\,\,\ol{6}\,\,7\,\,\ol{4}\,\,8\,\,\ol{1}\,\,9\,\,5\,\,\ol{10}\,\,3].\ee
Then~$h_w=p_u\circ c_v$, with~$v=(0,1,0,1,0,1,0,0,1,0)^\top$ and~$u=[2\,\,6\,\,7\,\,4\,\,8\,\,1\,\,9\,\,5\,\,10\,\,3]$. 
The latter can be written as a product of cycles as~$u = (1\,\,2\,\,6)(3\,\,7\,\,9\,10)(4)(5\,\,8)$. 
This~results in~$\chi(1\,\,2\,\,6)=v_1+v_2+v_6 = 2$ and 
similarly,~$\chi(5\,\,8)=0$, and~$\chi(3\,\,7\,\,9\,10)=\chi(4)=1$. Therefore,~$u_+=(1\,\,2\,\,6)(5\,\,8)$ and 
$u_-=(3\,\,7\,\,9\,10)(4)$ and thus,
\be t_\pm(h) = \left\{ \begin{array}{cc}(0,1,1,0,0,0,0,0,0,0)\\ (1,0,0,1,0,0,0,0,0,0)\end{array}\right., \ee
is the {\em signed cycle type} of~$h$, corresponding to the partitions~$2+3$ and~$1+4$ of~$5$.~\hfill~$\diamondsuit$
\end{Ex}
In view of Corolary~\ref{ECo-0}, we will now state one of the main results in this section.
\begin{Th} Two elements~$g,h\in\Bn$ are conjugate in~$\Bn$ if and only~$t_\pm(g)=t_\pm(h)$.\end{Th} 
Thus, we have been successful in our aim to characterize the conjugacy classes of~$\Bn$. 
Now we will consider the question of counting the number 
of elements of each of these classes. Firstly, Definition~\ref{Ecycle-type} implies 
that~$S_\ell$ has~$p(\ell)$ conjugacy classes, where~$p(\ell)$ stands for 
the number of integer partitions of~$\ell$. The sizes of these classes are well known.
\begin{Pro} The size of each conjugacy class of~$S_\ell$, being the number 
of~$\sigma\in S_\ell$ such that~$t(\sigma)= (t_1,\dots,t_\ell)$ for a given 
cycle type~$(t_1,\dots,t_\ell)$ equals
\be\label{Epolya-18} \left[\begin{array}{c} \ell \\ t\end{array}\right] = 
\frac{\ell!}{1^{t_1}\cdot\ldots\cdot\ell^{t_\ell} \cdot t_1!\cdot\ldots\cdot {t_\ell}!}.\ee
\end{Pro}
Similarly, writing~$\Delta(n)$ for the number of double partitions 
of~$n$, we have that~$\Delta(n)$ is the number of conjugacy classes of~$\Bn$ and that 
\be \Delta(n) = \sum_{k=0}^n p(k)p(n-k). \ee
For illustration, we list here the first few values of~$p$ and~$\Delta$ in Table~\ref{Epolya-19}. 
\begin{table}[h]
\begin{center}
$
\begin{array}{|r||r|r|r|r|r|r|r|r|r|r|r|r|r|r|}
\hline
n & 0 & 1 & 2 & 3 & 4 & 5 & 6 & 7 & 8 & 9 & 10 & 11 & 12 & 13\\
\hline
\hline
p(n) & 1 & 1 & 2 & 3 & 5 & 7 & 11 & 15 & 22 & 30 & 42 & 56 & 77 & 101\\
\hline
\Delta(n)  & 1 & 2 & 5 & 10 & 20 & 36 & 65 & 110 & 185 & 300 & 481 & 752 & 1165 & 1770\\
\hline
\end{array}
$
\end{center}
\caption{\small{Sequences A000041 and A000712 in the {\em Online Encyclopedia of Integer Sequences}.}}
\label{Epolya-19}
\end{table}
\begin{Pro} The number of elements in~$\Bn$ of the signed cycle type
\be\label{Eexpr} t_\pm(h) = \left\{ \begin{array}{cc} (t_1,\dots,t_n)\\ (s_1,\dots,s_n)\end{array}\right. 
\hdrie\mbox{\rm equals }\hdrie \binom{n}{k} 
\left[\begin{array}{c} k \\ t\end{array}\right] \left[\begin{array}{c} \ell \\ s \end{array}\right]2^{n-\sum(t_j+s_j)}, \ee
where~$k=t_1\cdot 1+\dots+t_n\cdot n$ and~$\ell=s_1\cdot 1+\dots+s_n\cdot n$ are the sums of the respective parts.
\end{Pro} 
{\bf Proof. } The only factor that needs explanation is the power of two. A cycle of length~$m$ can be given signs in 
$2^{m}$ ways,~$2^{m-1}$ of which resulting in an even number of signs and~$2^{m-1}$ of which in an odd number of signs.\hfill~$\Box$
%%%%%%%%%%%%%%
%%%%%%%%%%
%%%%%%%%%%%%
%%%%%%%%%%%%%
\subsection{An algorithm for the cycle index of $\Bn$}\label{ESect-1.4}
We are now able to answer the question how many permutations in~$S_{2^n}$ of which cycle type 
are induced by the~$2^nn!$ elements of~$\Bn$ by implementing the following algorithm.\\[2mm]
{\bf Algorithm 1}: Counting and tabulating the induced permutations of~$\Bn$.\\[2mm]
Let~$n\in\mathbb{N}$ be given.\\[2mm]
{\bf Step 1. } Generate the~$\Delta(n)$ double partitions~$(\tau_+,\tau_-)\,\vdash (k,\ell)$ of~$n$.\\[2mm]
{\bf Step 2. } For each such double partition, construct a single~$h\in\Bn$ with signed cycle type 
\[ t_{\pm}(h) = \left\{ \begin{array}{r} t_+(h) \\ t_-(h) \end{array} \right. = \left\{ \begin{array}{r} \tau_+ \\ \tau_- \end{array}, \right. \] 
and evaluate the expression in (\ref{Eexpr}) to count how many of them there are.\\[2mm]
{\bf Step 3. } Compute the type~$t(\pi_h)$ of the permutation~$\pi_h\in S_{2^n}$ induced by~$h$.\\[2mm]
{\bf Step 4. } Accumulate the result of Steps 2 and 3 over all double partitions in a table.
\begin{Ex} The conjugacy classes of~$\mathcal{B}_3$ are indexed by the ten double partitions of~$3$. 
Below we list these ten, and at their left we show how many elements of that type there are in~$\mathcal{B}_3$.
\[ 1:\left\{\begin{array}{cc}(3,0,0)\\ (0,0,0)\end{array}\right.\hdrie 6:
\left\{ \begin{array}{cc}(1,1,0)\\ (0,0,0)\end{array}\right. \hdrie 8:\left\{ \begin{array}{cc}(0,0,1)\\ (0,0,0)\end{array}\right. 
\hdrie 3:\left\{ \begin{array}{cc}(2,0,0)\\ (1,0,0)\end{array}\right. \hdrie 6:\left\{ \begin{array}{cc}(0,1,0)\\ (1,0,0)\end{array}\right.\]
\[ 1:\left\{ \begin{array}{cc}(0,0,0)\\ (3,0,0)\end{array}\right. \hdrie 6:
\left\{ \begin{array}{cc}(0,0,0)\\ (1,1,0)\end{array}\right. \hdrie 8:\left\{ \begin{array}{cc}(0,0,0)\\ (0,0,1)\end{array}\right. 
\hdrie 3:\left\{ \begin{array}{cc}(1,0,0)\\ (2,0,0)\end{array}\right. \hdrie 6:\left\{ \begin{array}{cc}(1,0,0)\\ (0,1,0)\end{array}\right.\]
Table~\ref{Econjclass} lists for each conjugacy class its cardinality, together with one element~$h\in\mathcal{B}_3$ 
from that class, and the cycle type of its induced permutation~$\pi_h$ in~$S_8$.
\begin{table}[h]
\be\label{Ep-36}\small \begin{array}{|c||c|c||c|c||c|} 
\hline
\# & h & \pi_h & \pi_h & h & \#\\
\hline
\hline 
1& (1)(2)(3) & (8,0,0,0,0,0,0,0) & (0,4,0,0,0,0,0,0)& (\ol{1})(\ol{2})(\ol{3}) & 1\\
\hline
6 & (1)(2\phantom{)(}3) & (4,2,0,0,0,0,0,0) & (0,0,0,2,0,0,0,0)& (\ol{1})(\ol{2}\phantom{)(}3) & 6\\
\hline
8 & (1\phantom{)(}2\phantom{)(}3) & (2,0,2,0,0,0,0,0) & (0,1,0,0,0,1,0,0)& (\ol{1}\phantom{)(}2\phantom{)(}3) & 8\\
\hline
3 & (1)(2)(\ol{3}) & (0,4,0,0,0,0,0,0) & (0,4,0,0,0,0,0,0) &(1)(\ol{2})(\ol{3}) & 3\\
\hline
6 & (1\phantom{)(}2)(\ol{3}) & (0,4,0,0,0,0,0,0) & (0,0,0,2,0,0,0,0) &(1)(\ol{2}\phantom{)(}3) & 6\\
\hline
\end{array}
\ee
\caption{\small{Cycle types of induced permutations and their cardinality.}}
\label{Econjclass}
\end{table}\\[2mm]
It also illustrates Remark~\ref{Erem-1}: elements from distinct conjugacy classes of~$\mathcal{B}_3$ may 
induce permutations in~$S_8$ the same cycle type. Table~\ref{EperminS8} groups them together.
\begin{table}[h]
\begin{center}
\begin{tabular}{|c|c|c|}
 \hline 
\# & $\pi_h$\\
\hline
\hline
1 & (8,0,0,0,0,0,0,0)\\
\hline
6 & (4,2,0,0,0,0,0,0)\\
\hline
13 & (0,4,0,0,0,0,0,0)\\
\hline
8 & (2,0,2,0,0,0,0,0)\\
\hline
12 & (0,0,0,2,0,0,0,0)\\
\hline
8 & (0,1,0,0,0,1,0,0)\\
\hline 
\end{tabular}
\end{center}
\caption{\small{Cycle index of~$\mathcal{B}_3$ in tabulated form.}}
\label{EperminS8}
\end{table}

\smallskip

Note that instead of computing the cycle types of the induced permutations of 
all~$2^nn!$ elements of~$\Bn$, we need to compute only~$\Delta(n)$ of them.~\hfill~$\diamondsuit$
\end{Ex}
The usual way in which Table~\ref{EperminS8} is expressed, is as a {\em cycle index} polynomial \cite{Bru,Polya}. 
\begin{Def}{\rm The {\em cycle index} \index{cycle index} of the induced permutations on~$\Bn$ of the hyperoctahedral group is the polynomial
\be\label{Ep-14} Z_n(x_1,\dots,x_\ell) = \frac{1}{|\Bn|}\sum_{h\in\Bn}\prod_{i=1}^{2\ell(n)} x_i^{t_i}. \ee
Here,~$t_i$ is the~$i$-th entry of~$t(\pi_h)$ and~$\ell$ is the {\em Landau function}, which assigns to~$n$ 
the largest order of an element from the symmetric group~$S_n$. Its values
\be 1, 1, 2, 3, 4, 6, 6, 12, 15, 20, 30, 30, 60, 60, 84, \dots \ee
can be found as sequence A000793 of the {\em Online Encyclopedia of Integer Sequences}}.
\end{Def}
Combining Table (\ref{EperminS8}) and (\ref{Ep-14}), the cycle index polynomial~$Z_3$ of~$\mathcal{B}_3$ can be found~as
\be Z_3(x_1,x_2,x_3,x_4,x_5,x_6) = \frac{1}{48} \left(x_1^8 + 6x_1^4x_2^2 + 13x_2^4 + 8 x_1^2x_3^2 + 12x_4^2+8x_2x_6\right). \ee
Further explicit expressions for the cycle index polynomials~$Z_n$ of~$\Bn$ can be found in the literature  \cite{ChGu2} 
only for~$n\leq 6$. The above rather simple algorithm implemented on a personal computer can produce the table 
corresponding to~$Z_n$ of the form (\ref{Ep-37}) for each~$n\leq 10$ within a second. In the papers \cite{Che,HaHi}, 
the cycle type~$t(\pi_h)$ of the permutation ~$\pi_h\in S_{2^n}$ induced by~$h$ is expressed in terms of the signed cycle 
type of the signed permutation corresponding to~$h\in\Bn$. Although algebraically of interest, their expressions are 
unfortunately too abstract to generate explicit numbers in a straightforward way. The above algorithm solves that problem.
%%%%%%%%%%%%%
%%%%%%%%%%%%%%
%%%%%%%%%%%%% 
\section{The 0/1-polytopes in the unit~$n$-cube}\label{ESect-2}
A {\em 0/1-polytope}\index{ 0/1-polytope} \cite{KaZi} is the convex hull of a (possibly empty) subset~$V\subset \BB^n$. Since distinct subsets 
of~$\BB^n$ give rise to distinct 0/1-polytopes, we can and prefer to define a 0/1-polytope alternatively but equivalently 
as a map~$c:\BB^n\rightarrow \{0,1\}$, using the obvious correspondence
\be c:\BB^n\rightarrow \{0,1\}:\hdrie v \mapsto 1 \Leftrightarrow v\in V. \ee
Such a map can be seen as a {\em two-coloring} of the vertices of~$I^n$ with ``colors''~$0$ and~$1$. We denote the set of all 
maps~$\BB^n\rightarrow\{0,1\}$ by~$\PP_n$, and write~$\PP_n^k\subset\PP_n$ for all~$c\in\PP_n$ with the property that precisely 
$k$ elements of~$\BB^n$ are mapped to~$1$: these correspond to the 0/1-polytopes with exactly~$k$ vertices. Observe that
\be \PP_n = \bigcup_{k=0}^{2^n}\PP_n^k \und \left|\PP_n\right|  = 
\sum_{k=0}^{2^n}\left|\PP_n^k\right| =  \sum_{k=0}^{2^n}\binom{2^n}{k} = 2^{2^n}. \ee
The {\em double-exponential} growth of~$|\PP_n|$ in~$n$ is illustrated in the below table. Already for~$n=6$ it exceeds (by one) 
the number of grains of rice that the poor merchant claimed from the king in the legend of the chess board, as displayed in Table~\ref{Etable5}.
\\
\begin{table}[h]
\begin{center}
\begin{tabular}{|c|c|c|c|c|c|c|c|c|c|c|c|c|c|c|}
\hline   
$n$ & 0 & 1 & 2 & 3 & 4 & 5 & 6\\
\hline 
$|\PP_n|$ & 2 &4 & 16 & 256 & 65536 & 4294967296 & 18446744073709551616\\
\hline
\end{tabular}
\end{center}
\caption{\small{Doubly exponential growth of the number~$|\mathcal{P}_n|$ of 0/1-polytopes in~$I^n$.}}
\label{Etable5}
\end{table}

\smallskip

We assign to each 0/1-polytope~$c$ the unique integer~$\NN(c)$ between~$0$ and~$2^{2^n}-1$ as
\be\label{Ep-5} \NN(c) = \sum_{v\in\BB^n} c(v)2^{\beta(v)}, \ee
where~$\beta$ is the binary-to-decimal numbering of $v \in \BB^n$ introduced in Figure~\ref{figureE5}.
\begin{Ex} Depicted below are the 0/1-polytopes~$c$ in the unit square. 
A circle is drawn around~$v\in\BB^2$ if and only if~$c(v)=1$. The number in the center of each square is~$\NN(c)$. 
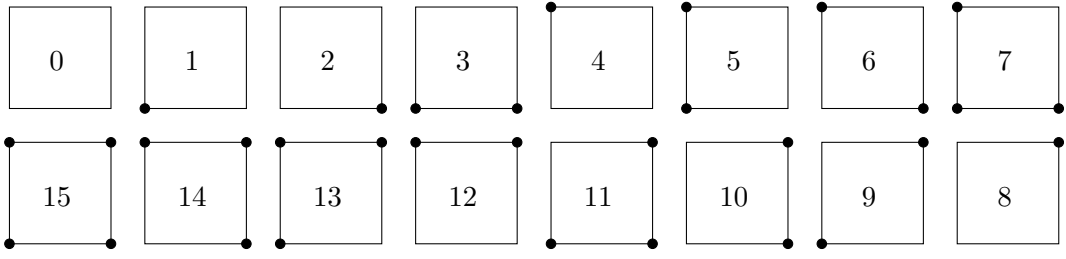
\begin{figure}[h]
\begin{center}
\begin{tikzpicture}[scale=0.9]
\begin{scope}
\draw (0,0)--(1.5,0)--(1.5,1.5)--(0,1.5)--cycle;
\node at (0.7,0.7) {$0$};
\end{scope}
\begin{scope}[shift={(2,0)}];
\draw (0,0)--(1.5,0)--(1.5,1.5)--(0,1.5)--cycle;
\node at (0.7,0.7) {$1$};
\draw[fill=black] (0,0) circle [radius=0.07];
\end{scope}
\begin{scope}[shift={(4,0)}];
\draw (0,0)--(1.5,0)--(1.5,1.5)--(0,1.5)--cycle;
\node at (0.7,0.7) {$2$};
\draw[fill=black] (1.5,0) circle [radius=0.07];
\end{scope}
\begin{scope}[shift={(6,0)}];
\draw (0,0)--(1.5,0)--(1.5,1.5)--(0,1.5)--cycle;
\node at (0.7,0.7) {$3$};
\draw[fill=black] (0,0) circle [radius=0.07];
\draw[fill=black] (1.5,0) circle [radius=0.07];
\end{scope}
\begin{scope}[shift={(8,0)}];
\draw (0,0)--(1.5,0)--(1.5,1.5)--(0,1.5)--cycle;
\node at (0.7,0.7) {$4$};
\draw[fill=black] (0,1.5) circle [radius=0.07];
\end{scope}
\begin{scope}[shift={(10,0)}];
\draw (0,0)--(1.5,0)--(1.5,1.5)--(0,1.5)--cycle;
\node at (0.7,0.7) {$5$};
\draw[fill=black] (0,0) circle [radius=0.07];
\draw[fill=black] (0,1.5) circle [radius=0.07];
\end{scope}
\begin{scope}[shift={(12,0)}];
\draw (0,0)--(1.5,0)--(1.5,1.5)--(0,1.5)--cycle;
\node at (0.7,0.7) {$6$};
\draw[fill=black] (1.5,0) circle [radius=0.07];
\draw[fill=black] (0,1.5) circle [radius=0.07];
\end{scope}
\begin{scope}[shift={(14,0)}];
\draw (0,0)--(1.5,0)--(1.5,1.5)--(0,1.5)--cycle;
\node at (0.7,0.7) {$7$};
\draw[fill=black] (0,0) circle [radius=0.07];
\draw[fill=black] (1.5,0) circle [radius=0.07];
\draw[fill=black] (0,1.5) circle [radius=0.07];
\end{scope}
\begin{scope}[shift={(0,-2)}];
\draw (0,0)--(1.5,0)--(1.5,1.5)--(0,1.5)--cycle;
\node at (0.7,0.7) {$15$};
\draw[fill=black] (1.5,0) circle [radius=0.07];
\draw[fill=black] (1.5,1.5) circle [radius=0.07];
\draw[fill=black] (0,1.5) circle [radius=0.07];
\draw[fill=black] (0,0) circle [radius=0.07];
\end{scope}
\begin{scope}[shift={(2,-2)}];
\draw (0,0)--(1.5,0)--(1.5,1.5)--(0,1.5)--cycle;
\node at (0.7,0.7) {$14$};
\draw[fill=black] (1.5,0) circle [radius=0.07];
\draw[fill=black] (1.5,1.5) circle [radius=0.07];
\draw[fill=black] (0,1.5) circle [radius=0.07];
\end{scope}
\begin{scope}[shift={(4,-2)}];
\draw (0,0)--(1.5,0)--(1.5,1.5)--(0,1.5)--cycle;
\node at (0.7,0.7) {$13$};
\draw[fill=black] (0,0) circle [radius=0.07];
\draw[fill=black] (1.5,1.5) circle [radius=0.07];
\draw[fill=black] (0,1.5) circle [radius=0.07];
\end{scope}
\begin{scope}[shift={(6,-2)}];
\draw (0,0)--(1.5,0)--(1.5,1.5)--(0,1.5)--cycle;
\node at (0.7,0.7) {$12$};
\draw[fill=black] (0,1.5) circle [radius=0.07];
\draw[fill=black] (1.5,1.5) circle [radius=0.07];
\end{scope}
\begin{scope}[shift={(8,-2)}];
\draw (0,0)--(1.5,0)--(1.5,1.5)--(0,1.5)--cycle;
\node at (0.7,0.7) {$11$};
\draw[fill=black] (0,0) circle [radius=0.07];
\draw[fill=black] (1.5,1.5) circle [radius=0.07];
\draw[fill=black] (1.5,0) circle [radius=0.07];
\end{scope}
\begin{scope}[shift={(10,-2)}];
\draw (0,0)--(1.5,0)--(1.5,1.5)--(0,1.5)--cycle;
\node at (0.7,0.7) {$10$};
\draw[fill=black] (1.5,1.5) circle [radius=0.07];
\draw[fill=black] (1.5,0) circle [radius=0.07];
\end{scope}
\begin{scope}[shift={(12,-2)}];
\draw (0,0)--(1.5,0)--(1.5,1.5)--(0,1.5)--cycle;
\node at (0.7,0.7) {$9$};
\draw[fill=black] (0,0) circle [radius=0.07];
\draw[fill=black] (1.5,1.5) circle [radius=0.07];
\end{scope}
\begin{scope}[shift={(14,-2)}];
\draw (0,0)--(1.5,0)--(1.5,1.5)--(0,1.5)--cycle;
\node at (0.7,0.7) {$8$};
\draw[fill=black] (1.5,1.5) circle [radius=0.07];
\end{scope}
\end{tikzpicture}
\end{center} 
\caption{\small{The sixteen 0/1-polytopes in~$I^2$ and their numbering given in (\ref{Ep-5}).~\hfill~$\diamondsuit$}}
\label{figureE8}
\end{figure}
\end{Ex}
Obviously,~$\NN(c)+\NN(\tilde{c})=2^{2^n}-1$ for {\em complementary} 0/1-polytopes, by 
which we mean polytopes~$c$ and~$\tilde{c}$ such that~$(c+\tilde{c})(v)=1$ for all~$v\in\BB^n$.
%%%%%%%%%%%%%
%%%%%%%%%%%%%
%%%%%%%%%%%%%
%%%%%%%%%%%%%%%%%%
\subsection{Cube symmetries acting on 0/1-polytopes: 0/1-equivalence}
Each element~$h$ of the hyperoctahedral group~$\Bn$ induces a 
permutation~$H_h$ of~$\PP_n$ by~$H_h: \PP_n \rightarrow \PP_n: c \mapsto c\circ h$. 
For each fixed~$k$ it restricts to a permutation of~$\PP_n^k\subset\PP_n$. 
Via the numbering~$\NN$ defined in (\ref{Ep-5}) it moreover induces a permutation
\be \Pi_h: \{0,\dots,2^{2^n}-1\} \rightarrow \{0,\dots,2^{2^n}-1\}: \hdrie k \mapsto \left(\NN\circ H_h \circ \NN^{-1}\right)(k). \ee
It turns out to be of interest to know the cardinalities~$|\SS|$ and~$|\SS^k|$ of the sets
\be \SS = \{ (h,c)\in\Bn\times\PP_n \sth c=c\circ h\} \und \SS^k = \{ (h,c)\in\Bn\times\PP_n^k \sth c=c\circ h\}.\ee 
Before explaining why, we present an example.
\begin{Ex} For each of the eight~$h\in\mathcal{B}_2$, the permutations~$\Pi_h$ of~$\{0,\dots,15\}$ 
are given in Table~\ref{Etable6}. The (bold) fixed points correspond to~$\SS$, and are added up per row and per column.
\\
\begin{table}[h]
\begin{center}
\small{
\begin{tabular}{|c|c|c|c|c|c|c|c|c|c|c|c|c|c|c|c|c|c|c|c|c|c|c|c|}
\hline 
$\mathcal{B}_2\times\PP_2$ & 0 & 1 & 2 & 3 & 4 & 5 & 6 & 7 & 8 & 9 & 10 & 11 & 12 & 13 & 14 & 15 & \\
 \hline
 \hline  
 $id$&{\bf 0}&{\bf 1}&{\bf 2}&{\bf 3}&{\bf 4}&{\bf 5}&{\bf 6}&{\bf 7}  &  {\bf 8}  &   {\bf 9}  &  {\bf 10}  & {\bf 11} &  {\bf 12}  &  {\bf 13}  & {\bf 14}  & {\bf 15} & 16\\
 \hline
 $c_1$ &{\bf 0}  &   2  &   1  &   {\bf 3}  &   8  &  10  &   9  &  11  &   4  &   6  &   5  &   7 &  {\bf  12}  &  14  &  13  &  {\bf 15} & 4\\
 \hline
$c_2$ & {\bf 0}& 4 & 8 & 12 & 1 & {\bf 5} & 9 & 13 & 2 & 6 & {\bf 10} & 14 & 3 & 7 & 11 & {\bf 15} & 4 \\
 \hline
 $c_1\circ c_2$          &  {\bf 0}  &   8  &   4  &  12  &   2  &  10  &  {\bf 6}  &  14  &   1  &  {\bf 9}  &   5  &  13 &  3  &  11  &   7  &  {\bf 15} & 4\\
 \hline
$ s_1$                   & {\bf 0}  &  {\bf 1}  &   4  &   5  &   2  &   3  &  {\bf 6}  &  {\bf 7}  &  {\bf 8}  &  {\bf  9}  &  12  &  13 &   10  &  11  & {\bf 14}  & {\bf 15} & 8\\
 \hline
$s_1\circ c_1$         &  {\bf 0}  &   4  &   1  &   5  &   8  &  12  &   9  &  13  &   2  &   6  &   3  &   7 &   10  &  14  &  11  & {\bf 15} & 2\\
 \hline
$s_1\circ c_2$        & {\bf 0}  &   2  &   8  &  10  &   1  &   3  &   9  &  11  &   4  &   6  &  12  &  14 &    5  &   7  &  13  & {\bf 15} & 2\\
 \hline
 $s_1\circ c_1\circ c_2$ & {\bf 0}  &   8  &  {\bf 2}  &  10  &   {\bf 4}  &  12  &  {\bf 6}  &  14  &   1  &  {\bf 9}  &   3  & {\bf 11} &    5  &  {\bf 13}  &   7  &  {\bf 15} & 8\\
 \hline
 \hline
  & 8 & 2 & 2 & 2 & 2 & 2 & 4 & 2 & 2 & 4 & 2 & 2 & 2 & 2 & 2 & 8 & 48\\
 \hline 
\end{tabular}}
\end{center}
\caption{\small{The action of~$\mathcal{B}_3$ on the~$16$ 0/1-polytopes in~$I^2$.}}
\label{Etable6}
\end{table} \\
We see directly that~$|\SS|=48$. After identifying the 0/1-polytopes with~$k$ vertices for 
given~$k\in\{0,\dots,4\}$, we moreover find that~$|\SS^0|=8$,~$|\SS^1|=8$,~$|\SS^2|=16$,~$|\SS^3|=8$, and~$|\SS^4|=8$.~\hfill~$\diamondsuit$
\end{Ex}
It may decrease complexity and uncover structure when we consider all elements in the orbit 
$\EE_n(c)$ of a 0/1-polytope~$c$ (elements in the same column of the above table) as equivalent.
\begin{Def}[0/1-equivalence] Two 0/1-polytopes~$c,\tilde{c}\in\PP_n$ for which there 
exists an~$h\in\Bn$ such that~$\tilde{c} = c\circ h$ are called {\em 0/1-equivalent}\index{0/1-equivalent}. 
\end{Def} 
It is clear that 0/1-equivalence of 0/1-polytopes implies their congruence; however, the converse 
does not hold \cite{KaZi}. Thus, 0/1-equivalence is a finer type of equivalence than congruence. 

\smallskip

We will now count the number~$\vep_n$ of 0/1-equivalence classes of 0/1-polytopes. Since 0/1-equivalent 
0/1-polytopes have the same number of vertices, we will count the number~$\vep_n^k$ of 0/1-equivalence 
classes of 0/1-polytopes with~$k$ vertices, after which~$\vep_n=\sum_k \vep_n^k$.
 \begin{Le}\label{Elem1} The number~$\vep_n^k$ of 0/1-equivalence classes of 0/1-polytopes with~$k$ vertices equals 
\be\label{Ep-6} \vep_n^k = \sum_{c\in\PP_n^k} \frac{1}{|\EE_n(c)|} = \sum_{c\in\PP_n^k}\frac{\left|\SS_c\right|}{\left|\Bn\right|},\hdrie\mbox{\rm where }\hdrie \SS_c = \{h\in\Bn \sth c = c\circ h\}. \ee
\end{Le}
{\bf Proof.} Trivially, given~$c\in\PP_n^k$, all~$\tilde{c}\in\PP_n^k$ that belong to~$\EE_n(c)$ contribute one to the sum. 
This proves the first equality in (\ref{Ep-6}). Next, if~$h\in\Bn$ is such that~$c\circ h =\tilde{c} \not=c$, then also 
$c\circ h_s\circ h=\tilde{c}$ if and only if~$h_s\in\SS_c$. Thus, for each~$\tilde{c}\in\EE_n(c)$ there are exactly 
$|\SS_c|$ elements of~$\Bn$ that map~$c$ onto~$\tilde{c}$, proving the second equality. \hfill~$\Box$
\begin{Co} \label{ECor-1} We have that
\be \vep_n^k= \frac{|\SS^k|}{|\Bn|} \und \vep_n = \frac{|\SS|}{|\Bn|} \hdrie 
\mbox{\rm due to }\hdrie  \sum_{c\in\PP_n^k} |\SS_c| = |\SS^k| \und \sum_{k=0}^{2^n} |\SS^k| = |\SS|. \ee
\end{Co}
Using this corollary we can continue to look at the example for~$n=2$.\\[2mm]
\begin{Ex}[continued] Since for~$n=2$ there are~$48$ elements in~$\SS$, by Corollary~\ref{ECor-1} 
we find~$48/8=6$ distinct 0/1-equivalence classes of 0/1-polytopes, being
\be \{0\}, \hdrie \{1,2,4,8\}, \hdrie \{3,5,10,12\}, \hdrie \{6,9\}, \hdrie \{7,11,13,14\} \und \{15\},\ee
consisting of 0/1-polytopes with zero, one, two, two, three, and four vertices, respectively. For all the 
0/1-polytopes with, for instance, two vertices, the fixed points add up to~$16$, confirming 
the existence of two distinct 0/1-equivalence classes in~$\PP_2^2$.~\hfill~$\diamondsuit$
\end{Ex}
\subsection{Counting 0/1-polytopes invariant under a given symmetry}
Counting the elements of~$\SS$ can be done by counting for each~$c\in\PP_n$ the number of elements 
of the set~$\SS_c$ from (\ref{Ep-6}). Alternatively, one can also count for each~$h\in\Bn$ the number of elements of the set
\be \SS_h = \{ c \in \PP_n \sth c = c \circ h \}. \ee
This is the Cauchy-Frobenius Lemma, also known as Burnside's Lemma \cite{Bru}. Note that 
counting~$\SS_c$ corresponds to counting the fixed points {\em per column} of Table~\ref{Etable6}.
\begin{Pro}\label{Epro-1} For given~$h\in\Bn$,~$c=c\circ h$ if and only if for each cycle of~$\pi_h$,~$c$ 
is constant on the pre-image under~$\beta$ of all numbers in that cycle. 
\end{Pro}
As a consequence, the number of 0/1-polytopes that are invariant under a given~$h\in\Bn$ 
with~$t(\pi_h)=(t_1,\dots,t_{2^n})$ equals the number of subsets of the set of the~$t_1+\dots+t_{2^n}$ cycles of~$\pi_h$.
\begin{Co}\label{Eco-5} Let~$h\in\Bn$ be given with~$t(\pi_h)=(t_1,\dots,t_{2^n})$. Then 
the cardinality~$|\SS_h|$ of the set~$\SS_h$ equals~$|\SS_h| = 2^{t_1+\dots+t_{2^n}}$.
\end{Co}
{\bf Proof. } According to Proposition~\ref{Epro-1}, the numbers within the same cycle of~$\pi_h$ must 
either all be mapped to~$0$ or all be mapped to~$1$ by~$c\circ\beta^{-1}$. \hfill~$\Box$

\smallskip

The number of 0/1-polytopes {\em with~$k$ vertices} that are invariant under a given cube 
symmetry~$h$ equals the number of subsets of the set of~$t_1+\dots+t_{2^n}$ cycles of~$\pi_h$ {\em whose lengths sum to}~$k$.
\begin{Th}\label{Eth-2} Let~$h\in\Bn$ be given with~$\tau=t(\pi_h)=(t_1,\dots,t_{2^n})$, and let~$k\leq 2^n$. 
The cardinality~$|\SS_h^k|$ of the set~$\SS_h^k=\{ c \in \PP_n^k \sth c = c \circ h \}$ equals
\be\label{Ep-c} |\SS_h^k| = \sum_{\kappa\,\vdash\, k} a(\tau,\kappa),\hdrie
\mbox{\rm where } \hdrie a(\tau,\kappa) = \prod_{j=1}^k\binom{t_j}{\kappa_j}, \ee
and where the sum ranges over all integer partitions~$\kappa$ of~$k$.
\end{Th} 
{\bf Proof. } Let~$(\kappa_1,\dots,\kappa_k)$ be a partition of~$k$. The number of ways that this partition can be 
selected from the partition~$(t_1,\dots,t_n)$ of~$n$ equals the product over all~$j\in\{1,\dots,k\}$ of the number 
of ways that~$\kappa_j$ cycles of length~$j$ can be selected from the~$t_j$ cycles of length~$j$.\hfill~$\Box$
\begin{Co}\label{Eco-7} The number of 0/1-equivalence classes of~$\PP_n^k$ equals 
\be\label{Ep-41} \vep_n^k = \frac{1}{|\Bn|}\sum_{h\in\Bn} \sum_{\kappa\,\vdash\, k} a(t(\pi_h),\kappa). \ee
\end{Co}
\begin{Ex} Consider the induced permutation~$\pi_h$ of the vertices of~$I^3$ with cycle type
\be\label{Ep-43} t(\pi_h) = (4,2,0,0,0,0,0,0), \ee
which consists of~$4+2=6$ cycles. Hence, the number of 0/1-polytopes that are mapped upon themselves by~$h$ equals~$2^6=64$, 
which illustrates Corollary~\ref{Eco-5}. To illustrate Theorem~\ref{Eth-2}, consider the five 
partitions of~$k=4$, being~$1+1+1+1 = 1+1+2 = 1+3 = 2+2 = 4$ and their corresponding cycle types,
\be (4,0,0,0), \hdrie (2,1,0,0), \hdrie (1,0,1,0), \hdrie (0,2,0,0), \hdrie (0,0,0,1). \ee
Only the first, second, and fourth partition contribute to the sum in (\ref{Ep-c}), which evaluates to
\be\label{Ep-37} \binom{4}{4}\binom{2}{0}\binom{0}{0}\binom{0}{0} +
\binom{4}{2}\binom{2}{1}\binom{0}{0}\binom{0}{0}+ \binom{4}{0}\binom{2}{2}\binom{0}{0}\binom{0}{0} = 14. \ee
Thus, each~$h\in\Bn$ with induced cycle type~$t(\pi_h)$ as in (\ref{Ep-43}), leaves 
invariant fourteen 0/1-polytopes in~$I^3$ with four vertices. See also Table~\ref{Etable7} in the next example.~\hfill~$\diamondsuit$ 
\end{Ex}
%%%%%%%%%%%%%
%%%%%%%%
%%%%%%%%%%%%%%%
\subsection{Counting the 0/1-equivalence classes of~$\PP_n^k$}
Corollary~\ref{Eco-7} in combination with the considerations in Section~\ref{ESect-1} 
give a way to compute the number~$|\PP_n^k|$ of 0/1-equivalence 
classes of 0/1-polytopes with~$k$ vertices as follows.\\[2mm]
{\bf Algorithm 2. }
Let integers~$n,k$ with~$0\leq k \leq 2^n$ be given.\\[2mm]
(1) Use Algorithm 1 from Section~\ref{ESect-1.4} to generate the cycle index~$Z_n$ of~$\Bn$ in tabulated form.\\[2mm] 
(2) Generate a second table with the~$p(k)$ partitions of~$k$, see for instance \cite{knu}.\\[2mm]
(3) For each cycle type~$\tau=(t_1,\dots,t_{2^n})\vdash 2^n$ from the first table:\\[2mm]
\hspace*{5mm} (a) sum the numbers~$a(\tau,\kappa)$ from (\ref{Ep-c}) over all~$\kappa \vdash k$;\\[2mm]
\hspace*{5mm} (b) multiply the result by the number of~$h\in\Bn$ for which~$t(\pi_h)=t$.\\[2mm]
(4) Sum over all~$\tau=(t_1,\dots,t_{2^n})\vdash 2^n$ from the first table.\\[2mm]
To illustrate this algorithm, we perform it in detail in the example below.
\begin{Ex} We consider the case~$n=3$ and~$k=4$. The part of Table~\ref{Etable7} to the left of the 
$6\times 5$ block in boldface is the table representing~$Z_3$ from Table~(\ref{EperminS8}). 
\\
\begin{table}[h]
\begin{center}
\small{
\begin{tabular}{|r||rrrr||r|r|r|r|r||r|}
\hline
 & & & & & 0 & 1 & 0 & 2 & 4 & $\kappa_1$\\
 & & & & & 0 & 0 & 2 & 1 & 0 & $\kappa_2 $\\
 & & & & & 0 & 1 & 0 & 0 & 0 & $\kappa_3$ \\
 & & & & & 1 & 0 & 0 & 0 & 0 & $\kappa_4$ \\
 \hline 
 \hline
 12 & 0 & 0 & 0 & 2 & {\bf 2} & {\bf 0} & {\bf 0} & {\bf 0} & {\bf 0} & 24\\
 \hline
 8  & 0 & 1 & 0 & 0 & {\bf 0} & {\bf 0} & {\bf 0} & {\bf 0} & {\bf 0} & 0\\
 \hline
 13 & 0 & 4 & 0 & 0 & {\bf 0} & {\bf 0} & {\bf 6} & {\bf 0} & {\bf 0} & 78\\
 \hline
  8 & 2 & 0 & 2 & 0 & {\bf 0} & {\bf 4} & {\bf 0} & {\bf 0} & {\bf 0} & 32\\
  \hline
  6 & 4 & 2 & 0 & 0 & {\bf 0} & {\bf 0} & {\bf 1} & \bf{12} & \bf{1} & 84\\
  \hline
  1 & 8 & 0 & 0 & 0 & {\bf 0} & {\bf 0} & {\bf 0} & {\bf 0} & {\bf 70} & 70\\
  \hline
  \hline
 48 & $t_1$ & $t_2$ & $t_3$ & $t_4$ & 24  & 32 & 84 & 72 & 76 & 288\\
  \hline
\end{tabular}}
\end{center}
\caption{\small{Computation of the numbers~$a(\tau,\kappa)$.}}
\label{Etable7}
\end{table}\\
The part of Table~\ref{Etable7} above the boldface part contains the five partitions of~$4$. Note that 
only the values of~$t_1,\dots,t_4$ are needed to be able to compute each of the numbers 
$a(\tau,\kappa)$. The numbers~$1,12,1$ in the fifth row in boldface are the ones computed in (\ref{Ep-37}) 
of the previous example. The sum~$288$ of the numbers in the~$6\times 5$ block divided by the order 
$48$ of~$\mathcal{B}_3$ equals~$6$, indeed the number of nonequivalent 0/1-tetrahedra in 
the cube, including two degenerate ones, as depicted in Figure~\ref{figureE9}.~\hfill~$\diamondsuit$
\begin{figure}[!]
\begin{center}
\begin{tikzpicture}[scale=0.8, every node/.style={scale=0.8}]
\begin{scope}
\draw[fill=gray!20!white] (0,0)--(2.5,1)--(0.5,3)--cycle;
\draw (0,0)--(2,0)--(2,2)--(0,2)--cycle;
\draw (0.5,1)--(2.5,1)--(2.5,3)--(0.5,3)--cycle;
\draw (0,0)--(0.5,1);
\draw (2,0)--(2.5,1);
\draw (0,2)--(0.5,3);
\draw (2,2)--(2.5,3);
\draw[fill=black] (0,0) circle [radius=0.05];
\draw[fill=black] (0.5,1) circle [radius=0.05];
\draw[fill=black] (2.5,1) circle [radius=0.05];
\draw[fill=black] (0.5,3) circle [radius=0.05];
\end{scope}
\begin{scope}[shift={(2.2,-2.3)}]
\draw[fill=gray!20!white] (0,0)--(2,0)--(0.5,3)--cycle;
\draw (0,0)--(2,0)--(2,2)--(0,2)--cycle;
\draw (0.5,1)--(2.5,1)--(2.5,3)--(0.5,3)--cycle;
\draw (0,0)--(0.5,1);
\draw (2,0)--(2.5,1);
\draw (0,2)--(0.5,3);
\draw (2,2)--(2.5,3);
\draw (2,0)--(0.5,1);
\draw[fill=black] (0,0) circle [radius=0.05];
\draw[fill=black] (0.5,1) circle [radius=0.05];
\draw[fill=black] (2,0) circle [radius=0.05];
\draw[fill=black] (0.5,3) circle [radius=0.05];
\end{scope}
\begin{scope}[shift={(5,0)}]
\draw[fill=gray!20!white] (0,0)--(0.5,1)--(2.5,3)--(2,0)--cycle;
\draw (0,0)--(2,0)--(2,2)--(0,2)--cycle;
\draw (0.5,1)--(2.5,1)--(2.5,3)--(0.5,3)--cycle;
\draw (0,0)--(0.5,1);
\draw (2,0)--(2.5,1);
\draw (0,2)--(0.5,3);
\draw (2,2)--(2.5,3);
\draw (0,0)--(2.5,3);
\draw (2,0)--(0.5,1);
\draw[fill=black] (0,0) circle [radius=0.05];
\draw[fill=black] (0.5,1) circle [radius=0.05];
\draw[fill=black] (2,0) circle [radius=0.05];
\draw[fill=black] (2.5,3) circle [radius=0.05];
\end{scope}
\begin{scope}[shift={(7.2,-2.3)}]
\draw[fill=gray!20!white] (0,0)--(2.5,1)--(2,2)--(0.5,3)--cycle;
\draw (0,0)--(2,0)--(2,2)--(0,2)--cycle;
\draw (0.5,1)--(2.5,1)--(2.5,3)--(0.5,3)--cycle;
\draw (0,0)--(0.5,1);
\draw (2,0)--(2.5,1);
\draw (0,2)--(0.5,3);
\draw (2,2)--(2.5,3);
\draw (0,0)--(2,2);
\draw (2.5,1)--(0.5,3);
\draw[fill=black] (0,0) circle [radius=0.05];
\draw[fill=black] (0.5,3) circle [radius=0.05];
\draw[fill=black] (2.5,1) circle [radius=0.05];
\draw[fill=black] (2,2) circle [radius=0.05];
\end{scope}
\begin{scope}[shift={(10,0)}]
\draw[fill=gray!20!white] (0,0)--(2,0)--(2,2)--(0,2)--cycle;
\draw (0,0)--(2,0)--(2,2)--(0,2)--cycle;
\draw (0.5,1)--(2.5,1)--(2.5,3)--(0.5,3)--cycle;
\draw (0,0)--(0.5,1);
\draw (2,0)--(2.5,1);
\draw (0,2)--(0.5,3);
\draw (2,2)--(2.5,3);
\draw (0,0)--(2,2);
\draw (2,0)--(0,2);
\draw[fill=black] (0,0) circle [radius=0.05];
\draw[fill=black] (0,2) circle [radius=0.05];
\draw[fill=black] (2,0) circle [radius=0.05];
\draw[fill=black] (2,2) circle [radius=0.05];
\end{scope}
\begin{scope}[shift={(12.2,-2.3)}]
\draw[fill=gray!20!white] (0,0)--(2,0)--(2.5,3)--(0.5,3)--cycle;
\draw (0,0)--(2,0)--(2,2)--(0,2)--cycle;
\draw (0.5,1)--(2.5,1)--(2.5,3)--(0.5,3)--cycle;
\draw (0,0)--(0.5,1);
\draw (2,0)--(2.5,1);
\draw (0,2)--(0.5,3);
\draw (2,2)--(2.5,3);
\draw (0,0)--(2.5,3);
\draw (2,0)--(0.5,3);
\draw[fill=black] (0,0) circle [radius=0.05];
\draw[fill=black] (2,0) circle [radius=0.05];
\draw[fill=black] (2.5,3) circle [radius=0.05];
\draw[fill=black] (0.5,3) circle [radius=0.05];
\end{scope}
\end{tikzpicture}
\end{center} 
\caption{\small{Representatives of each of the six 0/1-equivalence classes of all 0/1-tetrahedra.}}
\label{figureE9}
\end{figure}
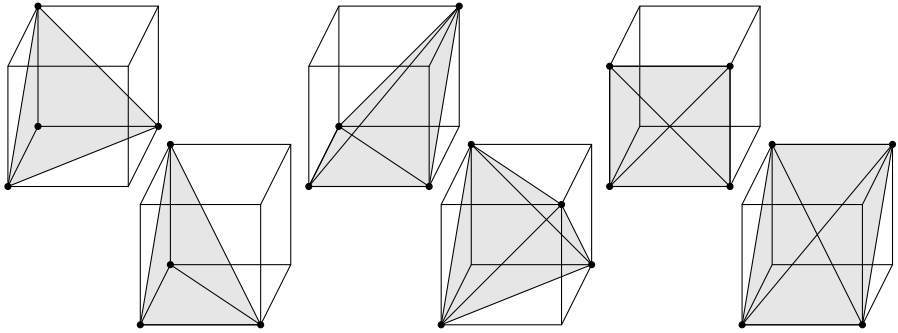
\end{Ex}
\begin{rem} The combinatorial road \cite{Bru,Polya} to arrive at the same result is to substitute 
in the cycle index polynomial~$Z_n$ of~$\Bn$ in (\ref{Ep-14}) the expressions~$x_i=b^i+w^i$. Then the coefficient 
of the monomial~$w^kb^{2^n-k}$ in the expansion equals~$|\PP_n^k|$. Although theoretically elegant 
and valuable, and widely applicable, it is not very suited for computing concrete numerical values.
\end{rem} 
The methodology described in Sections~\ref{ESect-1} and~\ref{ESect-2} leads to a way to compute the 
{\em number} of 0/1-polytopes in~$I^n$ modulo the symmetries of~$I^n$. See also Section~\ref{ESect-7} 
for some explicitly computed values. It does not yield a {\em specific element} from each 
0/1-equivalence class. In the next section we will investigate this {\em enumeration problem}.
 %%%%%%%%%%
 %%%%%%%%%%%%
 %%%%%%%%%%%%
 %%%%%%%%%%%%%%

\section{Minimal matrix representations of 0/1-polytopes}\label{ESect-3}
We will now designate in each 0/1-equivalence class~$\EE_n(c)$ of a 0/1-polytope 
$c$ some special {\em representatives}. One of them we denote as the {\em minimal representative} 
from that equivalence class. Obvious candidates for such minimal representatives are the 0/1-polytopes~$c$ for which
\be\label{Ezoniet} \NN(c) \leq \NN(c\circ h) \hdrie\mbox{\rm for all~$h\in\Bn$},\ee
where~$\NN$ is the numbering defined in (\ref{Ep-5}). However, with this definition it 
may happen that~$c\in\PP_n^k$ is a minimal representative, whereas none of its 
facets in~$\PP_n^{k-1}$ is minimal. For computational purposes, we prefer a minimal representative to have that property.
%%%%%%%
%%%%%%%
%%%%%%%
%%%%%%%
%%%%%%%
\subsection{Matrix representations of 0/1-polytopes}
A natural way to represent a 0/1-polytope~$c$, alternative to a mapping 
$c:\BB^n\rightarrow \{0,1\}$, is by means of 0/1-matrices whose colums are the vertices of~$c$.
\begin{Def}[Matrix representation]\label{Edef0}{\rm If the columns of a 
matrix~$P$ of size~$n\times k$ are precisely the~$k$ distinct vertices of a 0/1-polytope 
$c\subset\PP_n^k$, we will call~$P$ a {\em matrix representation}\index{matrix representation} of~$c$.  
With each matrix representation~$P$ we associate the integer vectors
\be \nu(P) = v_n^\top P\Pi_1 \und \mu(P) = \Pi_2 Pw_k,\ee 
where 
\be\label{Evnwk} \hdrie v_n^\top = \left(2^0,2^1,\dots,2^{n-1}\right) \und w_k^\top = \left(2^{k-1},\dots,2^1,2^0\right), \ee 
and where~$\Pi_1$ is the {\em unique}~$k\times k$ column permutation matrix sorting the~$k$ 
entries of~$v_n^\top P\Pi_1$ from left to right in {\em increasing order}, and~$\Pi_2$ {\em any} 
$n\times n$ row permutation matrix sorting the~$n$ possibly 
non-distinct entries of~$\Pi_2 Pw_k$ from top to bottom in {\em non-increasing} order}.
\end{Def}  
The permutations~$\Pi_1,\Pi_2$ depend on~$P$, but this dependence is suppressed 
from the notation. As~$P$ has distinct columns, each~$c\in\PP_n^k$ has exactly~$k!$ distinct matrix representations.
\begin{Pro}\label{Eprop2} The following statements are equivalent:\\[2mm]
$(1)$~$P_1$ and~$P_2$ are matrix representations of the same 0/1-polytope~$c\in\PP_n^k$;\\[2mm]
$(2)$ there exists a~$k\times k$ permutation matrix~$\Pi$ such that~$P_1=P_2\Pi$;\\[2mm]
$(3)$~$\nu(P_1) = \nu(P_2)$.
\end{Pro}
{\bf Proof. } This is because that no vertex of~$I^n$ is a convex combination 
of other vertices of~$I^n$, hence 0/1-polytopes are uniquely determined by their vertex set. \hfill~$\Box$

\smallskip

Due to the equivalence~$(1)\Leftrightarrow(3)$ in Proposition~\ref{Eprop2}, and with a slight abuse 
of notation, we will use~$\nu(c)$ also for a 0/1-polytope~$c\in\PP_n$, and assign to it the 
vector~$\nu(c)$ taken by any matrix representation~$P$ of~$c$. Note that~$\nu:\PP_n\rightarrow \{0,\dots,2^n\}^k$ is injective.

\smallskip

For given~$c\in\PP_n^k$, we will write~$\MM(c)$ for the set of all matrix representations of
 all~$\tilde{c}\in\EE_n(c)$. This induces an equivalence relation on the set~$\ZZ$ of all matrix 
 representations of 0/1-polytopes, that we will denote by~$P_1\sim P_2$. 
 Before studying this equivalence on~$\ZZ$, we introduce two simpler relations.
\begin{Def}[Row complementation/permutation]\label{Erow-c}{\rm A 0/1-matrix 
$P_2$ is a {\em row complementation}\index{row complementation} of~$P_1$, denoted by~$P_2 \csim P_1$, if it 
results from~$P_1$ after exchanging the zeros and ones in a subset of its rows; it 
is a {\em row permutation}\index{row permutation} of~$P_1$, denoted by~$P_2\psim P_1$ if there exists a 
permutation matrix~$\Pi$ such that~$P_2=\Pi P_1$.\hfill~$\Box$}
\end{Def}
Both~$\csim$ and~$\psim$ are equivalence relations on~$\ZZ$. The~$2^n$ row complementations applied 
to a given~$P\in\ZZ$ result in matrix representations of each of the 0/1-polytopes in an orbit under 
the action of the subgroup~$\Bn^c\subset\Bn$, whereas the~$n!$ row permutations of~$P$ are 
matrix representations of those in an orbit under the action of the subgroup~$\Bn^p\subset\Bn$. 
Thus, following Section~\ref{ESect-1}, a matrix representation of each 0/1-polytope that is in the 
same 0/1-equivalence class of a given~$c\in\PP_n^k$ can be obtained by performing each of the 
$2^nn!$ combined row complementations and permutations to a given matrix representation~$P$ of~$c$.
  %%%%%%%%%%%%%%%
  %%%%%%%%%%%%
  %%%%%%%%%%%%
  %%%%%%%%%%%%%%%
\subsection{Verification of 0/1-equivalence of matrix representations}  
For given~$P_1,P_2\in\BB^{n\times k}$ let~$r_1=P_1w_k$ and~$r_2=P_2w_k$, where~$w_k$ is the vector from (\ref{Evnwk}). 
If~$P_1\csim P_2$, then the~$j$-th entries of~$r_1$ and~$r_2$ are equal in case the~$j$-th rows of~$P_1$ and 
$P_2$ are equal, and add up to~$2^k-1$ in case these rows are complementary. Hence, verification whether 
$P_1\csim P_2$ can be done in at most~$\mathcal{O}(nk)$ operations, which is dominated by the costs of 
computing~$r_1$ and~$r_2$. Verifying whether~$P_1\psim P_2$ asks to inspect if the~$n$-vector~$r_1$ is a 
permutation of~$r_2$, requiring~$\mathcal{O}(nk+n\log n)$ operations.

\smallskip

The combination of these two observations yields the following, which can be seen as a variant 
of stating that~$r_1$ is a signed permutation of~$r_2$, see Section~\ref{ESect-1.3}.

\begin{Pro}\label{Eprop4} Let~$P_1,P_2\in\ZZ$. There exists an~$R\in\BB^{n\times k}$ such that 
\be P_1\psim R \csim P_2 \ee
if and only if the~$2n$ entries of the two~$n$-vectors~$r_1$ and~$(2^k-1)e-r_1$ are a 
permutation of the~$2n$ entries of the two~$n$-vectors~$r_2$ and~$(2^k-1)e-r_2$.

\smallskip

The verification if~$ P_1\psim R \csim P_2$ requires~$\mathcal{O}(nk+2n\log(2n))=\mathcal{O}(nk+n\log n)$ operations.
\end{Pro}
\begin{Ex} Consider matrices~$P_1$ and~$P_2$, with~$r_1=P_1w_4$ and~$r_2=P_2w_4$ computed as 
\[ P_1w_4= \left[\begin{array}{rrrr} 0 & 1 & 1 & 0\\ 0 & 0 & 1 & 0\\ 0 & 0 & 0 & 1\end{array}\right]
\left[\begin{array}{r} 8 \\ 4 \\ 2 \\ 1 \end{array}\right]  = \left[\begin{array}{r} 6 \\ 2 \\ 1 \end{array}\right]  
\und P_2w_4 =  \left[\begin{array}{rrrr} 0 & 1 & 1 & 0\\ 1 & 1 & 1 & 0\\ 1 & 1 & 0 & 1\end{array}\right]
\left[\begin{array}{r} 8 \\ 4 \\ 2 \\ 1 \end{array}\right] = \left[\begin{array}{r} 6 \\ 14 \\ 13 \end{array}\right].\]
Then~$P_1\csim R \psim P_2$ for some~$R\in\BB^{3\times 4}$, because the entries~$(6,2,1,9,13,14)$ of the 
two~$3$-vectors~$r_1$ and~$15e-r_1$ can be permuted into the entries~$(6,14,13,9,1,2)$ of~$r_2$ and~$15e-r_2$. 
Indeed,~$P_2$ is obtained by exchanging and complementing the second and third row of~$P_1$. \hfill~$\diamondsuit$
\end{Ex}
Proposition~\ref{Eprop2} showed that if~$P_1$ and~$P_2$ are matrix representations of {\em the same} 0/1-polytope, 
then~$P_1=P_2\Pi$ for a permutation matrix~$\Pi$ that can be found by inspecting if~$s_1=v_n^\top P_1$ is a 
permutation of~$s_2=v_n^\top P_2$, with~$v_n$ from (\ref{Evnwk}). To verify if~$P_1\sim P_2$, or in other words, 
if~$P_2$ is a {\em column permutation} of a row complementation and row permutation of~$P_1$, is computationally much more complex.
\begin{Pro}\label{Eprop3} Let~$P_1,P_2\in\ZZ$. Then~$P_1\sim P_2$ if and only 
if there exists an~$R\in\ZZ$ and a permutation matrix~$\Pi$ such that  
\be\label{Ep-44} P_1 \psim R \csim P_2\Pi,\ee
the verification of which can be done in~$\mathcal{O}(k!(nk+n\log n))$ operations.
\end{Pro} 
{\bf Proof. } The verification can be done by looping over all~$k!$ permutation matrices 
$\Pi$ and performing the verification in Proposition~\ref{Eprop4} for each of them. \hfill~$\Box$ 

\smallskip

\begin{rem}{\rm  Relation (\ref{Ep-44}) holds if there exist permutation matrices~$\Pi_1$ and~$\Pi_2$ such that 
\be \label{Ep-45} \Pi_1P_1 \csim P_2\Pi_2. \ee 
This more symmetric formulation suggests that in order to verify if~$P_1\sim P_2$, one can 
establish the existence of permutation matrices~$\Pi_1$ and~$\Pi_2$ such that~$\Pi_1 P_1\csim P_2\Pi_2$ in two ways:\\[2mm]
$\bullet$ for each~$\Pi_2$, verify if there exists~$\Pi_1$ such that~$\Pi_1 P_1\csim P_2\Pi_2$;\\[2mm]
$\bullet$ for each~$\Pi_1$, verify if there exists~$\Pi_2$ such that~$\Pi_1 P_1\csim P_2\Pi_2$.\\[2mm]
The second strategy would require an efficient way to verify the existence of a {\em column} permutation of 
$P_2$ such that it equals a row complementation of the given matrix~$\Pi_1 P_1$. This verification is far less 
trivial than the one in Proposition~\ref{Eprop4}. Nevertheless, if~$k>n+1$, there are ways to repair this and 
make the second strategy more economic than the first. Because our main interest is 0/1-{\em simplices} 
for which~$k\leq n+1$, we will not go into detail.}
\end{rem}
%%%%%%%%%%
%%%%%%%%%%%%
%%%%%%%%%%%%%%
%%%%%%%%%%%%%%%
\subsection{Minimal matrix representations and their properties}
The lexicographical order~$\prec$ on the integer vectors~$\nu(P)$ associated 
with the matrix representations~$P$ of 0/1-polytopes~$c\in\PP_n$ induces 
a total order on~$\PP_n$ as well as on~$\ZZ$.
\begin{Def}[Minimal representative]\label{Edef3}{\rm The {\em minimal representative}\index{minimal representative}  
of a $c\in\PP_n^k$ in~$\EE_n(c)$ is the unique 0/1-polytope $c^\ast\in\EE_n(c)$ for 
which~$\nu(c^\ast)\prec \nu(d)$ for all~$d\in\EE_n(c^\ast), d\not=c^\ast$}.
\end{Def}
The minimal representative~$c^\ast$ of~$c\in\PP_n^k$ has~$k!$ distinct matrix representations, 
of which we designate one as the minimal matrix representation. 
\begin{Def}[Minimal matrix representation]\label{Edef4}{\rm The 
{\em minimal matrix representation}\index{minimal matrix representation} of $c\in\PP_n$ in~$\MM(c)$ is 
the unique matrix representation~$P^\ast$ of~$c^\ast$ for which 
$v_n^\top P^\ast=\nu(P^\ast)$, in other words, whose column 
numbers~$v_n^\top P^\ast$ are strictly increasing}.
\end{Def}
We will now study further properties of minimal matrix representations of 
0/1-equivalence classes of 0/1-polytopes. The following result proves 
a desirable property, mentioned already at the beginning of this section.
\begin{Le}\label{Elem6} Let~$P^\ast$ be a minimal matrix representation of a 
0/1-polytope~$c\in\PP_n^k$. Then for each~$j\in\{1,\dots,k-1\}$, the 
submatrix~$P_j^\ast$ of~$P^\ast$ consisting of its~$j$ leftmost columns is 
a minimal matrix representation of a 0/1-polytope~$c_j\in\PP_n^j$.
\end{Le}
{\bf Proof. } Let~$P^\ast\in\BB^{n\times k}$ be a minimal matrix representation. 
Then by Definition~\ref{Edef4}, ~$v_n^\top P^\ast$ is increasing, and hence, so is 
$v_n^\top P^\ast_{k-1}$. To arrive at a contradiction, assume that~$P^\ast_{k-1}$ 
is not a minimal matrix representation. Then there exists a row permutation~$\Pi$ such 
that~$\Pi P^\ast_{k-1}\csim P_{k-1}$ and~$\nu(P_{k-1})\prec\nu(P^\ast_{k-1})$. But this 
means that~$\Pi P^\ast \csim P$, where~$P$ is a matrix whose~$k-1$ leftmost columns equal 
$P_{k-1}$. Irrespective of the rightmost column of~$P$, this implies that~$\nu(P)\prec\nu(P^\ast)$, 
contradicting the minimality of~$P^\ast$. This proves that~$P^\ast_{k-1}$ is a minimal matrix 
representation, and hence inductively, the minimality of all~$P^\ast_j$.\hfill~$\Box$
\begin{Co} Any minimal representative of a 0/1-polytope with~$k$ vertices 
contains a minimal representative of a 0/1-polytope with~$k-1$ vertices. 
\end{Co}
\begin{Co}\label{Eco5} The first column of a minimal matrix 
representation~$P^\ast$ of a 0/1-polytope~$c\in\PP_n^k$ equals~$0\in\RR^n$.
\end{Co} 
{\bf Proof. } According to Lemma~\ref{Elem6}, the first column of~$P^\ast$ is a 
minimal matrix representation of a 0/1-polytope with one vertex. 
Clearly, this is the zero vector.\hfill~$\Box$

\smallskip

By Definition~\ref{Edef3},~$\nu(P^\ast)=v_n^\top P^\ast$, which means that the 
integer vector~$v_n^\top P^\ast$ is increasing. The next lemma proves that 
additionally,~$P^\ast w_p$ is non-increasing from top to bottom.
\begin{Le}\label{Elem4} Let~$P^\ast$ be a minimal matrix representation of a 0/1-polytope~$c\in\PP_n^k$. Then 
\be P^\ast w_k=\mu(P^\ast) \ee
or equivalently,~$P^\ast w_k$ is non-increasing from top to bottom.
\end{Le}  
{\bf Proof. } Write~$p_i^\ast$ for row~$i$ of~$P^\ast$ and~$p_j^\ast$ for 
row~$j$ of~$P^\ast$. Assume that~$1\leq i<j\leq n$ and~$p_i^\ast w_k < p_j^\ast w_k$, contradicting the statement to prove. Then 
\be p = \min_{\ell\in\{1,\dots,k\}}\left\{ p_i^\ast e_\ell \not= p_j^\ast e_\ell\right\} \ee
exists and equals the index of the leftmost column in which~$p_i^\ast$ and 
$p_j^\ast$ differ. The assumption~$p_i^\ast w_k < p_j^\ast w_k$ implies 
that~$p_i^\ast e_k=0$ and~$p_j^\ast e_k=1$. Write~$P$ for the matrix that 
results after the transposition of rows~$i$ and~$j$ of~$P^\ast$. Then 
the first~$p-1$ columns of~$P^\ast$ and~$P$ coincide. However, in column 
$p$, the one in row~$j$ is exchanged with the zero in row~$i$ above it. As 
a result,~$v_n^\top P\prec v_n^\top P^\ast$, contradicting that~$P^\ast$ 
is a minimal matrix representation. \hfill~$\Box$ 

\smallskip

In Figure \ref{figureE10} we display four elements~$c_1,\dots,c_4$ from~$\EE_3(c_1)$ 
of a tetrahedron~$c_1$.  The stabilizer~$\SS_{c_1}$ of~$c_1$ in 
$\mathcal{B}_3$ contains two elements, hence~$|\EE_3(c_1)|=24$. We also display 
matrix representations for~$c_1,\dots,c_4$, all with increasing column numbers.
\begin{figure}[h]
\begin{center}
\begin{tikzpicture}[scale=0.9, every node/.style={scale=0.9}]
\begin{scope}[shift={(4,0)}]
\draw[fill=gray!20!white] (0,0)--(2,0)--(2.5,1)--(0,2)--cycle;
\draw (0,0)--(2,0)--(2,2)--(0,2)--cycle;
\draw (0.5,1)--(2.5,1)--(2.5,3)--(0.5,3)--cycle;
\draw (0,0)--(0.5,1);
\draw (2,0)--(2.5,1);
\draw (0,2)--(0.5,3);
\draw (2,2)--(2.5,3);
\draw (2,0)--(0,2);
\draw (0,0)--(2.5,1);
\draw[fill=black] (0,0) circle [radius=0.05];
\draw[fill=black] (2.5,1) circle [radius=0.05];
\draw[fill=black] (2,0) circle [radius=0.05];
\draw[fill=black] (0,2) circle [radius=0.05];
\node[scale=0.8] at (1,-1) {$\left[\begin{array}{rrrr}0& 1 & 1 & 0 \\0&  0 & 1 & 0\\ 0& 0 & 0 & 1\end{array}\right]$};
\node[scale=0.8] at (2.5,-1) {$\left[\begin{array}{rrr} 6\\2\\1\end{array}\right]$};
\node[scale=0.8] at (1,-1.9) {$\left[\begin{array}{rrrr}0 & 1 & 3 & 4\end{array}\right]$};
\node at (1.5,2.5) {$c_2$};
\end{scope}
\begin{scope}[shift={(0,0)}]
\draw[fill=gray!20!white] (0,0)--(2,0)--(0.5,3)--cycle;
\draw (0,0)--(2,0)--(2,2)--(0,2)--cycle;
\draw (0.5,1)--(2.5,1)--(2.5,3)--(0.5,3)--cycle;
\draw (0,0)--(0.5,1);
\draw (2,0)--(2.5,1);
\draw (0,2)--(0.5,3);
\draw (2,2)--(2.5,3);
\draw (2,0)--(0.5,1);
\draw[fill=black] (0,0) circle [radius=0.05];
\draw[fill=black] (0.5,1) circle [radius=0.05];
\draw[fill=black] (2,0) circle [radius=0.05];
\draw[fill=black] (0.5,3) circle [radius=0.05];
\node[scale=0.8] at (1,-1) {$\left[\begin{array}{rrrr}0& 1 & 0 & 0 \\0&  0 & 1 & 1\\ 0& 0 & 0 & 1\end{array}\right]$};
\node[scale=0.8] at (2.5,-1) {$\left[\begin{array}{rrr} 4\\3\\1\end{array}\right]$};
\node[scale=0.8] at (1,-1.9) {$\left[\begin{array}{rrrr}0 & 1 & 2 & 6\end{array}\right]$};
\node at (1.5,2.5) {$c_1$};
\end{scope}
\begin{scope}[shift={(8,0)}]
\draw[fill=gray!20!white] (0,0)--(2,0)--(2,2)--(0.5,1)--cycle;
\draw (0,0)--(2,0)--(2,2)--(0,2)--cycle;
\draw (0.5,1)--(2.5,1)--(2.5,3)--(0.5,3)--cycle;
\draw (0,0)--(0.5,1);
\draw (2,0)--(2.5,1);
\draw (0,2)--(0.5,3);
\draw (2,2)--(2.5,3);
\draw (2,0)--(0.5,1);
\draw (0,0)--(2,2);
\draw[fill=black] (0,0) circle [radius=0.05];
\draw[fill=black] (2,2) circle [radius=0.05];
\draw[fill=black] (2,0) circle [radius=0.05];
\draw[fill=black] (0.5,1) circle [radius=0.05];
\node[scale=0.8] at (1,-1) {$\left[\begin{array}{rrrr}0& 1 & 0 & 1 \\0&  0 & 1 & 0\\ 0& 0 & 0 & 1\end{array}\right]$};
\node[scale=0.8] at (2.5,-1) {$\left[\begin{array}{rrr} 5\\2\\1\end{array}\right]$};
\node[scale=0.8] at (1,-1.9) {$\left[\begin{array}{rrrr}0 & 1 & 2 & 5\end{array}\right]$};
\node at (1.5,2.5) {$c_3$};
\end{scope}
\begin{scope}[shift={(12,0)}]
\draw[fill=gray!20!white] (0,0)--(0,2)--(2.5,3)--(2,2)--cycle;
\draw (0,0)--(2,0)--(2,2)--(0,2)--cycle;
\draw (0.5,1)--(2.5,1)--(2.5,3)--(0.5,3)--cycle;
\draw (0,0)--(0.5,1);
\draw (2,0)--(2.5,1);
\draw (0,2)--(0.5,3);
\draw (2,2)--(2.5,3);
\draw (0,0)--(2,2);
\draw (0,0)--(2.5,3);
\draw[fill=black] (0,0) circle [radius=0.05];
\draw[fill=black] (0,2) circle [radius=0.05];
\draw[fill=black] (2.5,3) circle [radius=0.05];
\draw[fill=black] (2,2) circle [radius=0.05];
\node[scale=0.8] at (1,-1) {$\left[\begin{array}{rrrr}0& 0 & 1 & 1 \\0&  0 & 0 & 1\\ 0& 1 & 1 & 1\end{array}\right]$};
\node[scale=0.8] at (2.5,-1) {$\left[\begin{array}{rrr} 3\\1\\7\end{array}\right]$};
\node[scale=0.8] at (1,-1.9) {$\left[\begin{array}{rrrr}0 & 4 & 5 & 7\end{array}\right]$};
\node at (1.5,0.5) {$c_4$};
\end{scope}
\end{tikzpicture}
\end{center} 
\caption{\small{Four 0/1-equivalent 0/1-tetrahedra with corresponding matrix representations.}}
\label{figureE10}
\end{figure}
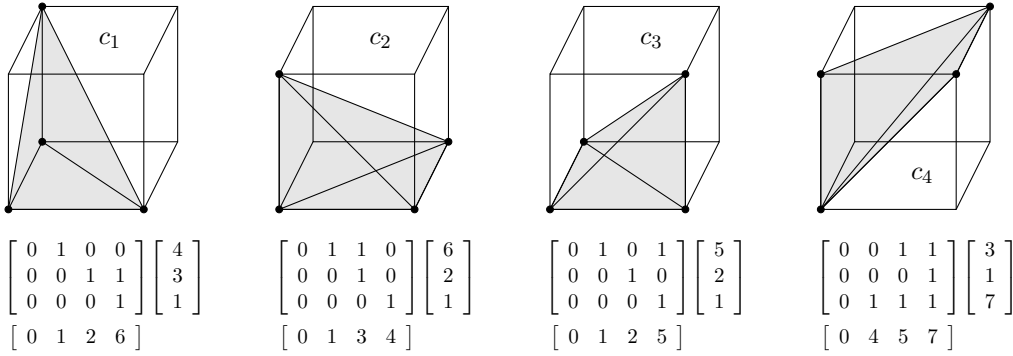\\[2mm]
Without proof, we mention the following facts for illustration.\\[2mm]
$\bullet$~$\NN(c_2)\leq\NN(c)$ for all~$c\in\EE_3(c_1)$, but none of its triangular facets is~$\NN$-minimal;\\[2mm]
$\bullet$~$c_3$ is the unique minimal representative of~$\EE_3(c_1)$;\\[2mm]
$\bullet$~$\mu(c_4)$ is not nonincreasing, hence~$c_4$ is not the minimal representative of~$\EE_3(c_1)$;\\[2mm]
$\bullet$ no 0/1-polytope formed by two or three vertices of~$c_4$ is a minimal representative;\\[2mm]
$\bullet$~$c_1$ has all properties proved above of the minimal representative, but is not it.\\[2mm]
Now, let~$P$ be any matrix representation of a 0/1-simplex~$c$ with~$k\leq n+1$ vertices. 
For each~$j\in\{1,\dots,k\}$, write~$P_j$ for the matrix obtained from~$P$ by 
first complementing those rows in~$P$ that have an entry~$1$ in column~$j$, 
and then swapping columns~$1$ and~$j$ of the result. Each matrix~$P_j$ 
corresponds to~$c$ with one of its vertices placed at the origin. Hence, due to Corollary~\ref{Eco5}, there exists a~$j\in\{1,\dots,k\}$ such that 
\[ P^\ast = \Pi_1P_j\Pi_2 \] 
for some permutation matrices~$\Pi_1$ and~$\Pi_2$, and where~$\Pi_2$ leaves 
the first column of~$P_j$ invariant. Instead of applying all the~$n!(k-1)!$ 
permutations and verifying which of them result in~$P^\ast$, we can use 
that by Lemma~\ref{Elem4}, the row numbers~$P^\ast w_k$ are non-increasing. 
Therefore, for each of the~$(k-1)!$ choices for~$\Pi_2$, it suffices to simply sort 
the rows of~$P_j\Pi_2$. Then~$P^\ast$ must be among the resulting~$k(k-1)! = k!$ 
matrices,~$(k-1)!$ for each value of~$j\in\{1,\dots,k\}$.\\[2mm]
Summarized in algorithmic form, this reads as follows.\\[2mm]
{\bf Algorithm 3}: Computing the minimal matrix representation~$P^\ast$ of a 0/1-simplex~$c$.\\[2mm]
Let~$P$ be any matrix representation of~$c$. Define~$P_1,\dots,P_k$ as described above.\\[2mm]
For~$j=1,\dots,k$;\\[2mm]
$(1)$ apply all~$(k-1)!$ column permutations~$\Pi_2$ to~$P_j$ that leave its first column invariant;\\[2mm]
$(2)$ for each of those, apply any row permutation~$\Pi_1$ for which~$\Pi_1 P_j\Pi_2w_p$ is non-increasing;\\[2mm]
$(3)$ store the matrix~$P_j^\ast = \Pi_1 P_j \Pi_2$ for which~$\Pi_1 P_j\Pi_2w_p$ is lexicographically minimal.\\[2mm]
Each~$P_j^\ast$ can be seen as a {\em local minimizer} over all matrices 
that can be obtained from~$P_j$ by permuting its rows and columns. 
The minimal among all~$k$ local minima is then~$P^*$.
\begin{Co} A 0/1-polytope with~$k$ vertices has at most~$k!$ 
distinct matrix representations~$P$ with~$Pe_1=0$ and with~$Pw_k$ nonincreasing.
\end{Co} 
We will now use Algorithm 3 to enumerate the minimal representatives and 
their minimal matrix representations of all 0/1-{\em triangles} in~$I^n$, 
and in particular of the subset of all {\em acute} 0/1-triangles. These 
minimal matrix representations of 0/1-triangles will be extended by a 
computer code to minimal matrix representations of nonobtuse and acute 0/1-{\em simplices}.
%%%%%%%%%%%%%%
%%%%%%%%%%%%
%%%%%%%%%%%%%
%%%%%%%%%%%%%%
\subsection{All minimal matrix representations of 0/1-triangles}\label{Esect-4.4}
Let~$\TT\in\PP_n^3$ be a 0/1-triangle. We will characterize its minimal 
matrix representation~$T^\ast~$. By Definition~\ref{Edef4}, we know 
that~$v_n^\top T^\ast~$ is increasing; by Corollary~\ref{Eco5}, the 
first column of~$T^\ast~$ equals zero; and by Lemma~\ref{Elem4}, 
$T^\ast w_3$ is non-increasing. Therefore, we know that 
\be\label{Eadd-1}\small  T^\ast  = \left[ \begin{array}{c|ccc} a & 0 & 1 & 1 \\ 
b & 0 & 1 & 0 \\ c & 0 & 0 & 1 \\ d & 0 & 0 & 0\end{array}\right] \hdrie\mbox{\rm for certain}\hdrie a+b+c+d = n,\ee
and where the right-hand side stands for the~$n\times 3$ matrix 
whose top~$a$ rows equal $[0\ 1\ 1]$, whose next~$b$ 
rows equal~$[\,0\,\, 1\,\, 0\,]$ and so on. Of course,~$T^\ast$ is 
not minimal for {\em all} values of~$a,b,c,d$. For instance, if~$b>c$ 
it is not. In that case, swapping the second and third column of~$T^\ast$ 
and sorting the rows, leads to a matrix with a smaller second column number:~$2^{a+c}<2^{a+b}$.

\smallskip

To further specify~$a,b,c,d$, we compute the~$k!=6$ matrices that are generated 
by Algorithm 3, with start matrix~$T_1^\ast = T^\ast$ for some choice of~$a,b,c,d$. 
Each of those six matrices is of the same form as in (\ref{Eadd-1}), but with the 
numbers~$a,b,c,d$ of repeated rows permuted. Instead of writing down the 
matrices, we only present in Table~\ref{Etable8} their corresponding permutations of 
$a,b,c,d$, in three sets of two, each pair belonging to one of the 
matrices~$T^\ast_1,T^\ast_2,T^\ast_3$, where~$T_1^\ast=T^\ast$,
\begin{table}[h]
\begin{center}
\small{
\begin{tabular}{|c|c||c|c||c|c|}
\hline
$a$ & $a$ & $b$ & $b$ & $c$ & $c$\\
$b$ & $c$ & $a$ & $c$ & $a$ & $b$\\
$c$ & $b$ & $c$ & $a$ & $b$ & $a$\\
$d$ & $d$ & $d$ & $d$ & $d$ & $d$\\ 
\hline 
\end{tabular}}
\end{center}
\caption{\small{Induced (block-)permutations.}}
\label{Etable8}
\end{table}\\
Observe that Table \ref{Etable8} basically consists, in fact, of all six permutations of~$a,b,c$. 
\begin{Th}\label{Eth1} The matrix~$T^\ast$ is the minimal matrix representation 
of a 0/1-triangle~$\TT$ in~$I^n$ if and only if it is of the form (\ref{Eadd-1}) with
\be\label{Ezes} 1 \leq a+b \und a+b+c \leq n \und a \leq b \leq c. \ee
\end{Th} 
{\bf Proof. } Consider the column numbers of~$T^\ast$,
\be v_n^\top T_{11} = \left( 0, \hdrie 2^{a+b}-1, \hdrie 2^{a+b+c}-2^{a+b}+2^a-1\right). \ee
Necessary and sufficient conditions for this vector to be lexicographically 
minimal over all permutations of~$a,b,c$ are as follows. The second entry 
is minimal if and only if~$a+b\leq a+c$ and~$a+b\leq b+c$, hence if and 
only if~$b\leq c$ and~$a\leq c$. If this is the case, additionally the third 
entry is minimal if and only if~$a\leq b$. The fact that~$a+b+c$ must be 
bounded above by~$n$ is trivial. The additional bound~$1\leq a+b$ is a 
necessary and sufficient condition for the three vertices of~$\TT$ to be distinct. \hfill~$\Box$
 \begin{Co}\label{Eco1} Let~$T^\ast$ be a minimal matrix representation of a 0/1-triangle~$\TT$ in~$I^n$. Then:\\[2mm]
$\bullet$ if~$a=0$ then~$\TT$ is a right triangle;\\[2mm]
$\bullet$ if~$a>0$ then~$\TT$ has acute angles only.
\end{Co}   
{\bf Proof. } Suppose that~$a=0$. Then due to~$1 \leq a+b$ and~$a\leq b \leq c$ in (\ref{Ezes}), we have that~$1\leq b~$ and~$1\leq c$ and thus, 
$T^\ast$ in (\ref{Eadd-1}) obviously represents a nondegenerate right triangle. If~$0<a$, then again due to (\ref{Ezes}), also~$0<b\leq c$. This 
shows that the difference between the second and third column is not orthogonal to either one of them, and thus is~$T^\ast$ not right. Finally, 
since no triangle in~$I^n$ can have obtuse angles, also the second bullet is proved. \hfill~$\Box$

\smallskip

Theorem~\ref{Eth1} establishes a bijection between the minimal matrix representations of all 0/1-triangles in~$I^n$ and the set all points (except 
the origin due to~$1\leq a+b$) with integer coordinates in the polyhedron~$K$ in the nonnegative octant of~$\RR^3$ defined by the inequalities 
\be 0\leq a,\hdrie 0\leq b,  \hdrie 0\leq c  \und a+b+c \leq n \und a \leq b \leq c. \ee
A closer inspection shows that~$K$ is a tetrahedron, the intersection of the so-called {\em path-tetrahedron}~$P$ defined 
by the inequalities~$0\leq a\leq b\leq c\leq n$, and the {\em cube-corner}~$C$ defined by 
$0\leq a, 0\leq b, 0\leq c$ and~$a+b+c \leq n$. This is depicted in Figure~\ref{figureE11}. 
\begin{rem}{\rm The {\em right} triangles correspond to the integer points in the bottom facet of~$K$.}
\end{rem}
Because the cube~$[0,n]^3$ can be subdivided into six congruent path tetrahedra all sharing the same long diagonal,~$K$ is one of six congruent parts of 
the cube corner. In fact, each of those six parts corresponds to exactly one of the matrices in Table~\ref{Etable8}.
\begin{figure}[h]
\begin{center}
\begin{tikzpicture}
\begin{scope}
\draw[fill=gray!20!white] (1,1)--(3,3)--(4,1)--(3,0)--cycle;
\draw[thick] (1,1)--(3,0)--(4,1)--(3,3)--cycle;
\draw[thick] (1,1)--(4,1);
\draw[thick] (3,0)--(3,3);
\draw[gray] (0,0)--(3,0)--(3,3)--(0,3)--cycle;
\draw[gray] (1,1)--(4,1)--(4,4)--(1,4)--cycle;
\draw[gray] (0,0)--(1,1);
\draw[gray] (3,0)--(4,1);
\draw[gray] (3,3)--(4,4);
\draw[gray] (0,3)--(1,4);
\draw[fill=white] (1,1) circle [radius=0.05];
\draw[fill=black] (3,0) circle [radius=0.05];
\draw[fill=black] (4,1) circle [radius=0.05];
\draw[fill=black] (3,3) circle [radius=0.05];
\node at (0.7,1) {$o$};
\node at (0.6,2) {$a\uparrow$};
\node at (2,3.8) {$\rightarrow b$};
\node[rotate=225] at (0.4,3.7) {$\rightarrow$};
\node at (0.1,3.4) {$c$};
\node at (2.5,1.7) {$P$};
\end{scope}
\begin{scope}[shift={(5,0)}]
\draw[fill=gray!20!white] (1,4)--(4,1)--(0,0)--cycle;
\draw[thick] (1,4)--(4,1)--(0,0)--cycle;
\draw[thick] (1,1)--(1,4);
\draw[thick] (1,1)--(4,1);
\draw[thick] (1,1)--(0,0);
\draw[gray] (0,0)--(3,0)--(3,3)--(0,3)--cycle;
\draw[gray] (1,1)--(4,1)--(4,4)--(1,4)--cycle;
\draw[gray] (0,0)--(1,1);
\draw[gray] (3,0)--(4,1);
\draw[gray] (3,3)--(4,4);
\draw[gray] (0,3)--(1,4);
\draw[fill=white] (1,1) circle [radius=0.05];
\draw[fill=black] (4,1) circle [radius=0.05];
\draw[fill=black] (0,0) circle [radius=0.05];
\draw[fill=black] (1,4) circle [radius=0.05];
\node at (0.7,1) {$o$};
\node at (1.4,2) {$\uparrow a$};
\node at (2,3.8) {$\rightarrow b$};
\node[rotate=225] at (0.4,3.7) {$\rightarrow$};
\node at (0.1,3.4) {$c$};
\node at (-2.5,-0.5) {$P:\,\,0\leq a\leq b\leq c \leq n$ \hspace{1mm} and};
\node at (2.2,1.7) {$C$};
\end{scope}
\begin{scope}[shift={(10,0)}]
\draw[fill=gray!20!white] (1,1)--(1.66,1.66)--(4,1)--(2,0.5)--cycle;
\draw[thick] (1,1)--(1.66,1.66)--(4,1)--(2,0.5)--cycle;
\draw[thick] (1,1)--(4,1);
\draw[thick] (1.66,1.66)--(2,0.5);
\draw[gray] (0,0)--(3,0)--(3,3)--(0,3)--cycle;
\draw[gray] (1,1)--(4,1)--(4,4)--(1,4)--cycle;
\draw[gray] (0,0)--(1,1);
\draw[gray] (3,0)--(4,1);
\draw[gray] (3,3)--(4,4);
\draw[gray] (0,3)--(1,4);
\draw[fill=white] (1,1) circle [radius=0.05];
\draw[fill=black] (2,0.5) circle [radius=0.05];
\draw[fill=black] (1.66,1.66) circle [radius=0.05];
\draw[fill=black] (4,1) circle [radius=0.05];
\node at (0.7,1) {$o$};
\node at (0.6,2) {$a\uparrow$};
\node at (2,3.8) {$\rightarrow b$};
\node[rotate=225] at (0.4,3.7) {$\rightarrow$};
\node at (0.1,3.4) {$c$};
\node at (-3.5,-0.5) {$C:\,\,a+b+c \leq n$};
\node at (2.5,2) {$K$};
\node at (2,-0.5) {$K=P\cap C$};
\end{scope}
\end{tikzpicture}
\end{center} 
\caption{\small{The tetrahedron~$K$ as intersection a path tetrahedron~$P$ and a cube corner~$C$.}}
\label{figureE11}
\end{figure}
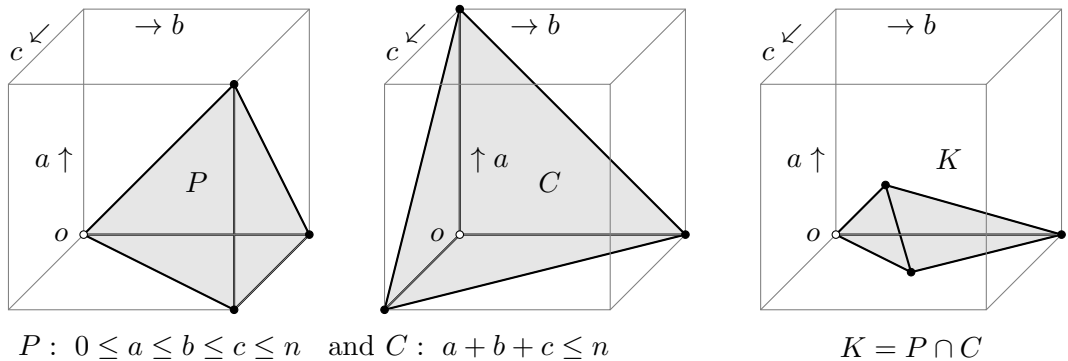\\[2mm]
Obviously, with the above characterizations, the enumeration of all minimal 
matrix representations of all 0/1-triangles, or of those of all acute 0/1-triangles is now a trivial matter.
 %%%%%%%%%%%%%%%%%%%%%%%
 %%%%%%%%%%%%%%%%
 %%%%%%%%%%%%%%%
 \subsection{Simple explicit expressions for~$\vep_n^2$}
The number of equivalence classes~$\vep_n^2$ of 0/1-triangles in~$I^n$ was in principle counted in Section~\ref{ESect-2}, using P\'olya's theory of 
counting. However, it did not provide a simple expression for this number, nor did it count the number of {\em acute} triangles. To do this here, we use 
the following change of variables,
\be\label{Echange} p = a+b, \hdrie q = a+c, \und r = a. \ee
The conditions on~$a,b$ and~$c$ in Theorem~\ref{Eth1} are equivalent to the conditions
\be\label{Ef7} 1\leq p\leq q \leq p+q-2r\leq n-r\leq n. \ee
We will count the triples~$(p,q,r)$ satisfying (\ref{Ef7}) by fixing a value for~$r$ and counting the tupels~$(p,q)$ that satisfy the resulting equation, and 
summing over the feasible values of~$r$.
\begin{Le}\label{Elem2} Let~$m\geq 1$ be an integer. The number of integer tupels~$(x,y)$ satisfying
\be\label{Ef12} 1\leq x \leq y \leq m-x \ee
equals
\be \left\lfloor \frac{m}{2} \right\rfloor \left\lceil \frac{m}{2} \right\rceil,\ee
where~$\lfloor \cdot \rfloor$ is the floor-operator and~$\lceil\cdot\rceil$ the ceil-operator. 
\end{Le} 
{\bf Proof.} Only for values of~$x$ with~$1\leq x \leq \lfloor m/2\rfloor$, we have that~$x\leq m-x$. The number of integers between such an~$x$ and 
$m-x$ equals~$m+1-2x$. This leads to a total of
\be \sum_{x=1}^{\lfloor m/2\rfloor} (m+1-2x) = \left\lfloor \frac{m}{2} \right\rfloor (m+1) - 2\cdot\half \left\lfloor \frac{m}{2} 
\right\rfloor \left( \left\lfloor \frac{m}{2} \right\rfloor +1 \right) \ee
tupels~$(x,y)$ that satisfy (\ref{Ef12}). Splitting~$m =\lceil\frac{m}{2}\rceil +\lfloor\frac{m}{2}\rfloor$ leads to the statement. \hfill~$\Box$
\begin{Co}\label{Eco2} The number of 0/1-equivalence classes of right triangles in~$I^n$ equals
\be\label{Ef13}  \left\lfloor \frac{n}{2} \right\rfloor \left\lceil \frac{n}{2} \right\rceil. \ee
\end{Co} 
{\bf Proof. } Combining Corollary~\ref{Eco1} with the change of variables in (\ref{Echange}) shows that we may set~$r=0$ in (\ref{Ef7}) and continue to count 
to number of tupels~$(p,q)$ satisfying
\be 1 \leq p \leq q \leq p+q \leq n. \ee
Since the inequality~$q\leq p+q$ is always valid, it can be removed. Thus, we only need to count the number 
of tupels~$(p,q)$ such that~$1\leq p\leq q \leq n-p$, which was done in Lemma~\ref{Elem2}. \hfill~$\Box$

\smallskip

In the next lemma we will count equivalence classes of triangles for fixed values of~$r\geq 1$. It will turn 
out that if~$3r>n$, there are no solutions. Moreover, substituting~$r=0$ in (\ref{Ef10}) below does 
not yield the result of Corollary~\ref{Eco2}. After its proof it is explained why not.
\begin{Le}\label{Ele3} For given~$r\geq 1$ with~$3r \leq n$, the number of tupels~$(p,q)$ satisfying (\ref{Ef7}) equals
\be \label{Ef10}\left\lfloor \frac{n-3r+2}{2} \right\rfloor \left\lceil \frac{n-3r+2}{2} \right\rceil.\ee  
\end{Le}
{\bf Proof. } Let~$r\geq 1$ be fixed. If~$p<2r$, there are no integers~$q$ that satisfy the 
third inequality~$q\leq p+q-2r$ in (\ref{Ef7}). If~$p\geq 2r$, this inequality holds for all~$q$ 
and can thus be removed. Thus, we only need to count the tupels~$(p,q)$ for which
\be\label{Ef11} 2r \leq p\leq q \leq n+r-p. \ee
For such tupels to exist, we need that~$2r\leq n+r-p$, but since~$p\geq 2r$ this 
translates into~$2r\leq n+r-2r$. This explains the requirement~$3r\leq n$ 
in the statement of this lemma. To count the tupels, subtract~$2r-1$ from each 
term in (\ref{Ef11}), and define~$x=p-(2r-1), y=q-(2r-1)$, and~$z=n-3r+2$, then
\be 1\leq x \leq y \leq n+r-x-2(2r-1) = n-3r+2-x = z-x. \ee
Applying Lemma~\ref{Elem2} gives the number of tupels~$(x,y)$ satisfying 
these inequalities in terms of~$z$, and substituting back~$z=n-3r+2$ proves the statement. \hfill~$\Box$
\begin{rem}{\rm Choosing~$r=0$ in (\ref{Ef10}) does not give (\ref{Ef13}). This is because 
setting~$r=0$ in (\ref{Ef11}) does not imply~$1\leq p$, as is required, whereas for~$r\geq 1$, it does.}
\end{rem}
We will now count the number of equivalence classes of acute triangles. First another lemma.
\begin{Le}\label{Elem3} For nonnegative integers~$k$ we have that~$(k\mod 2)^2 = k\mod 2$, and hence
\be \left\lfloor\frac{k}{2}\right\rfloor\left\lceil\frac{k}{2}\right\rceil  = 
\left(\frac{k-k\mod 2}{2}\right)\left(\frac{k+k\mod 2}{2}\right) = \frac{1}{4}(k^2-k\mod 2). \ee
Moreover,
\be \sum_{k=1}^n k\mod 2 = \left\lfloor\frac{n+1}{2} \right\rfloor, \und 
\left\lfloor \frac{n-\lfloor \frac{n}{3}\rfloor}{2}\right\rfloor = \left\lfloor \frac{n+1}{3}\right\rfloor.\ee
\end{Le}
{\bf Proof. } Elementary, and thus left to the reader. \hfill~$\Box$   
\begin{Th} \label{Eth4} The number of 0/1-equivalence classes of acute triangles in~$I^n$ equals
\be\label{Ef14}  \left\lfloor\frac{2n^3 +3n^2-6n+9}{72}\right\rfloor. \ee
\end{Th}
{\bf Proof. } We need to sum the expression in (\ref{Ef10}) over all~$r\geq 1$ satisfying~$3r\leq n$. 
Now, since~$(n-3r+2)\mod 2 = (n-r)\mod 2$, we find using Lemma~\ref{Elem3} that
\be \label{Ef17}\sum_{r=1}^{\lfloor \frac{n}{3} \rfloor} \left\lfloor \frac{n-3r+2}{2} \right\rfloor \left\lceil \frac{n-3r+2}{2} \right\rceil = \frac{1}{4} 
\sum_{r=1}^{\lfloor \frac{n}{3} \rfloor}(n-3r+2)^2 - \frac{1}{4} \sum_{r=1}^{\lfloor \frac{n}{3} \rfloor}(n-r)\mod 2. \ee 
The first sum in the right-hand side of (\ref{Ef17}) can be evaluated using standard expressions for sums of squares as
\be\label{Ef15} \sum_{r=1}^{\lfloor \frac{n}{3}\rfloor}(n-3r+2)^2 = \left\lfloor \frac{n}{3} \right\rfloor(n+2)\left(n-1-3\left\lfloor \frac{n}{3} 
\right\rfloor\right)+ \frac{3}{2}\left\lfloor \frac{n}{3} \right\rfloor\left(\left\lfloor \frac{n}{3} \right\rfloor+1\right)\left(2\left\lfloor\frac{n}{3} 
\right\rfloor+1\right).\ee
Using Lemma~\ref{Elem3} again, the second sum in the right-hand side of (\ref{Ef17}) evaluates to
\be \label{Ef16} \sum_{r=1}^{\lfloor \frac{n}{3}\rfloor}(n-r)\mod 2 = \sum_{r=1}^{n-1}r\mod 2 - 
\sum_{r=1}^{n-\lfloor \frac{n}{3}\rfloor-1}r\mod 2 = \left\lfloor \frac{n}{2}\right\rfloor - \left\lfloor \frac{n+1}{3}\right\rfloor. \ee
Combining (\ref{Ef17}), (\ref{Ef15}) and (\ref{Ef16}), the number of equivalence classes of acute 0/1-triangles equals
\be \label{Ef20} \frac{\left( \left\lfloor \frac{n}{3} \right\rfloor(n+2)\left(n-1-3\left\lfloor \frac{n}{3} \right\rfloor\right)+ 
\frac{3}{2}\left\lfloor \frac{n}{3} \right\rfloor\left(\left\lfloor \frac{n}{3} \right\rfloor+1\right)
\left(2\left\lfloor\frac{n}{3} \right\rfloor+1\right)-\left\lfloor \frac{n}{2}\right\rfloor +\left\lfloor \frac{n+1}{3}\right\rfloor\right)}{4}.\ee
To verify that this expression equals (\ref{Ef14}) is a tedious task, but can be done 
as follows. First, we substitute~$n=6k+\ell$ with~$\ell\in\{0,\dots,5\}$ into (\ref{Ef14}), which after simplifications results in 
\be \label{Ef19} 6k^3+ \frac{3}{2}\left(2\ell+1\right)k^2+\half\left(\ell^2+\ell-1\right)k+
\left\lfloor\frac{1}{36}\ell^3+\frac{1}{24}\ell^2-\frac{1}{12}\ell+\frac{1}{8}\right\rfloor,\ee  
where we have used that~$2\ell+1$ and~$\ell^2+\ell-1=\ell(\ell+1)-1$ are both 
odd, which implies that the sum of the first three terms in (\ref{Ef19}) is indeed an integer for all~$k$ and~$\ell$. 

\smallskip

Next, substitute~$n=6k+\ell$ with~$\ell\in\{0,1,2\}$ in (\ref{Ef20}), and note that it simplifies to
\be\label{Ef18} 6k^3+ \frac{3}{2}\left(2\ell+1\right)k^2+\half\left(\ell^2+\ell-1\right)k, \ee 
which equals the expression in (\ref{Ef19}) because for~$\ell\in\{0,1,2\}$ the 
floor results in zero. Finally, set~$n=6k+\ell$ with~$\ell\in\{3,4,5\}$ in (\ref{Ef20}). After simplification there remains 
\be\label{Ef21} 6k^3+ \frac{3}{2}\left(2\ell+1\right)k^2+\half\left(\ell^2+\ell-1\right)k+
\frac{1}{4}\left(\ell^2-2\ell+1-\left\lfloor\frac{\ell}{2}\right\rfloor+\left\lfloor\frac{\ell+1}{3}\right\rfloor\right). \ee 
Comparing (\ref{Ef19}) with (\ref{Ef21}), it can be easily verified that for~$\ell\in\{3,4,5\}$,
\be \frac{1}{4}\left(\ell^2-2\ell+1-\left\lfloor\frac{\ell}{2}\right\rfloor+
\left\lfloor\frac{\ell+1}{3}\right\rfloor\right) =\left\lfloor\frac{1}{36}\ell^3+\frac{1}{24}\ell^2-\frac{1}{12}\ell+\frac{1}{8}\right\rfloor. \ee
And this proves the theorem. \hfill~$\Box$

\smallskip

In Table \ref{Etable9} are listed the numbers~$r_n$ and~$a_n$ of 0/1-equivalence 
classes of right and acute 0/1-triangles and their sum~$d_n$ for small values of~$n$.
\begin{table}[h]
\begin{center}
\small{
\begin{tabular}{|r||r|r|r|r|r|r|r|r|r|r|r|r|r|r|r|r|r|r|r|r|}
\hline
$n$ &  2 & 3 & 4 & 5 & 6 & 7 & 8 & 9 & 10 & 11 & 12 & 13 & 14 & 15 & 16 & 17  \\
\hline
\hline
$r_n$ & 1 & 2 & 4 & 6 & 9 & 12 & 16 & 20 & 25& 30 & 36 & 42 & 49 & 56 & 64 & 72  \\
\hline
$a_n$ & 0 & 1 & 2 & 4 & 7 & 11 & 16 & 23 & 31 & 41 & 53 & 67 & 83 & 102 & 123 & 147 \\
\hline
$d_n$ & 1 & 3 & 6 & 10 & 16 & 23 & 32 & 43 & 56 & 71 & 89 & 109 & 132 & 158 & 187 & 219  \\
\hline
\end{tabular}}
\end{center}
\caption{\small{Right, acute, and all 0/1-triangles in~$I^n$ modulo the action of~$\Bn$.}}
\label{Etable9}
\end{table}\\
In the OEIS, the sequence~$r_n$ has label A002620, sequence~$a_n$ has label 
A181120, and~$d_n$ has label A034198. Only the latter has as 
description ``number of distinct triangles on vertices of~$n$-dimensional cube'', 
the other two are not associated with counting triangles in~$I^n$.
%%%%%%%%%%%%%%%
%%%%%%%%%%%%%%%%
%%%%%%%%%%%%%%%
\section{Minimal representatives of acute 0/1-simplices}\label{ESect-5}
We will now describe how to generate by means of a computer program 
minimal matrix representations of each 0/1-equivalence class of so-called 
{\em acute 0/1-simplices}, which are 0/1-simplices having only acute 
dihedral angles. They form the higher dimensional generalizations of the 
{\em acute 0/1-triangles} from the previous section.
\begin{Def}[Acute~$k$-simplex] \label{Edef-acute}{\rm Let~$c\in\PP_n^{k}$ 
with~$1\leq k\leq n$ be a nondegenerate~$k$-simplex in~$I^n$. Let 
$R^\ast\in\BB^{n\times(k+1)}$ be the minimal matrix representation 
of~$c$. Let~$P$ be the~$n\times k$ matrix with the nonzero columns 
of~$R^\ast$. If the~$k\times k$ Gramian~$G = P^\top P$ satisfies:\\[2mm]
$(1)$ each off-diagonal entry of~$G^{-1}$ is negative ($G^{-1}$ is {\em strictly Stieltjes}),\\[2mm]
$(2)$ each row sum of~$G^{-1}$ is positive ($G^{-1}$ is {\em diagonally dominant}),\\[2mm]
then~$c$ is called an {\em acute} 0/1-$k$-simplex.}
\end{Def}
The properties (1) and (2) are purely {\em geometric}, and concern the 
{\em dihedral angles} of the simplex, for which we refer to \cite{BrKoKr,BrKoKrSo} 
for details. These angles are invariant under the action of~$\Bn$. This guarantees 
that the concept of acute 0/1-simplex is well-defined using only the minimal 
matrix representation. Note that~$G$ is invertible as~$c$ is assumed nondegenerate.

\smallskip

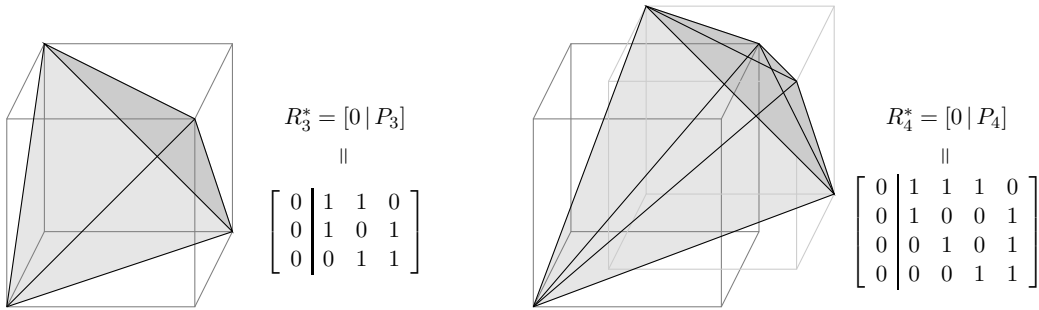
\begin{figure}[h]
\begin{center}
\begin{tikzpicture}
\begin{scope}
\begin{scope}[rotate=180]
\draw[gray!20!white, fill=gray!20!white] (2.5,0)--(0.5,1)--(0,2.5)--(3,3.5)--cycle;
\draw[gray!40!white, fill=gray!40!white] (2.5,0)--(0.5,1)--(0,2.5)--cycle;
\draw[gray] (0,0)--(2.5,0)--(2.5,2.5)--(0,2.5)--cycle;
\draw[gray] (0.5,1)--(3,1)--(3,3.5)--(0.5,3.5)--cycle;
\draw[gray] (0,0)--(0.5,1);
\draw[gray] (2.5,0)--(3,1);
\draw[gray] (2.5,2.5)--(3,3.5);
\draw[gray] (0,2.5)--(0.5,3.5);
\draw (2.5,0)--(0.5,1)--(0,2.5)--cycle;
\draw (3,3.5)--(2.5,0);
\draw (3,3.5)--(0.5,1);
\draw (3,3.5)--(0,2.5);
\end{scope}
\begin{scope}
\node[scale=0.8] at (1.5,-1) {$R_3^\ast=[0\,|\, P_3]$};
\node[scale=0.8,rotate=90] at (1.5,-1.5) {$=$};
\node[scale=0.8] at (1.5,-2.5) {$\left[\begin{array}{r|rrr} 0 & 1 & 1 & 0\\ 0 & 1 & 0 & 1\\ 0 & 0 & 1 & 1\end{array}\right]$};
\end{scope}
\begin{scope}
\node[scale=0.8] at (9.5,-1) {$R_4^\ast=[0\,|\, P_4]$};
\node[scale=0.8,rotate=90] at (9.5,-1.5) {$=$};
\node[scale=0.8] at (9.5,-2.5) {$\left[\begin{array}{r|rrrr} 0 & 1 & 1 & 1 & 0\\ 0 & 1 & 0 & 0 & 1\\ 0 & 0 & 1 & 0 & 1\\ 0 & 0 & 0 & 1 & 1 \end{array}\right]$};
\end{scope}
\end{scope}
\begin{scope}[shift={(8,0.5)}]
\begin{scope}[rotate=180]
\draw[gray!20!white, fill=gray!20!white] (2.5,0)--(1,0.5)--(0.5,1)--(0,2.5)--(4,4)--cycle;
\draw[gray!40!white, fill=gray!40!white] (2.5,0)--(1,0.5)--(0.5,1)--(0,2.5)--cycle;
\begin{scope}[shift={(1,0.5)}]
\draw[gray] (0,0)--(2.5,0)--(2.5,2.5)--(0,2.5)--cycle;
\draw[gray] (0.5,1)--(3,1)--(3,3.5)--(0.5,3.5)--cycle;
\draw[gray] (0,0)--(0.5,1);
\draw[gray] (2.5,0)--(3,1);
\draw[gray] (2.5,2.5)--(3,3.5);
\draw[gray] (0,2.5)--(0.5,3.5);
\end{scope}
\draw[gray!40!white] (0,0)--(2.5,0)--(2.5,2.5)--(0,2.5)--cycle;
\draw[gray!40!white] (0.5,1)--(3,1)--(3,3.5)--(0.5,3.5)--cycle;
\draw[gray!40!white] (0,0)--(0.5,1);
\draw[gray!40!white] (2.5,0)--(3,1);
\draw[gray!40!white] (2.5,2.5)--(3,3.5);
\draw[gray!40!white] (0,2.5)--(0.5,3.5);
\draw (2.5,0)--(1,0.5)--(0.5,1)--(0,2.5)--cycle;
\draw (2.5,0)--(0.5,1);
\draw (0,2.5)--(1,0.5);
\draw (4,4)--(2.5,0);
\draw (4,4)--(0.5,1);
\draw (4,4)--(0,2.5);
\draw (4,4)--(1,0.5);
\end{scope}
\end{scope}
\end{tikzpicture}
\end{center} 
\caption{\small{An acute antipodal~$n$-simplex for~$n=3$ (left) and~$n=4$ 
(right) together with their minimal matrix representations~$R_3^\ast$ and~$R_4^\ast$ 
and their nonsingular parts~$P_3$ and~$P_4$ that satisfy the conditions in Definition~\ref{Edef-acute}.}}
\label{figureE12}
\end{figure}

As examples of acute simplices, in Figure~\ref{figureE12} we display on the left the {\em only} 
acute  0/1-{\em tetrahedron} in~$I^3$, and on the right the {\em only} acute 
0/1-$4$-{\em simplex} in~$I^4$. Both are members of the family of so-called 
{\em antipodal}  simplices in~$I^n$. An antipodal~$n$-simplex in~$I^n$ is 
0/1-equivalent with the simplex whose vertices are the origin and all 
$v\in\BB^n$ with exactly one entry equal to zero. For this family, the 
matrices~$P$ and~$G$ from Definition~\ref{Edef-acute} are, indexed by~$n$,
\be\label{Efive} P_n = \left[\begin{array}{c|l} e^\top & 0 \\ \hline I_{n-1} & e \end{array}\right], 
\und G_n = I_n+ee^\top+(n-3)e_ne_n^\top.\ee
As before,~$e$ is the all-ones vector of appropriate length and~$I_\ell$ is the~$\ell\times\ell$ 
identity matrix. It is easy to verify that~$G^{-1}$ satisfies the criteria (1) and (2) in 
Definition~\ref{Edef-acute} for the family of antipodal simplices to be acute.
\\
 \begin{rem}\label{Erem-8} {\rm The tetrahedral facet~$T$ of the antipodal 
~$4$-simplex represented by the first four columns~$K$ of the matrix~$R_4^\ast$ 
 is {\em congruent} to the regular tetrahedron in the left picture. They are, 
 however, not 0/1-{\em equivalent}. Indeed,~$T$ does not lie in a 
 three-dimensional cubic facet of~$I^4$, and this property is invariant 
 under the action of~$\mathcal{B}_4$. A congruence~$Q$ mapping one onto the other is, for instance,
\be\label{Etwo} \small QK = \frac{1}{2}\left[\begin{array}{rrrr}1 & 1 & 1 & -1 \\ 1 & 1 & -1 & 1 \\ 1 & -1 & 1 & 1 \\ -1 & 1 & 1 & 1\end{array}\right] 
\left[\begin{array}{r|rrr} 0 & 1 & 1 & 1 \\ \hline 0 & 1 & 0 & 0  \\ 0 &0 & 1 & 0  \\0 & 0 & 0 & 1  \end{array}\right] =
\left[\begin{array}{r|rrr}0 & 1 & 1 & 0  \\ 0 &1 & 0 & 1  \\0 & 0 & 1 & 1 \\ \hline 0 & 0 & 0 & 0\end{array}\right], \ee
but this congruence~$Q$ is not a member of~$\mathcal{B}_4$.} 
\end{rem} 
%%%%%%%%%%%%
%%%%%%%%%%%%%%
%%%%%%%%%%%%%%%
%%%%%%%%%%%%%%%%%
\subsection{Acute 0/1-simplices and their candidate acute extensions}
Here we will list a number of properties of acute simplices that are relevant in 
the context of their computational enumeration. Some of them are new, others 
are simply valid for acute simplices in general \cite{BrKoKr,BrKoKrSo,Fie}. 
\begin{Pro}[\cite{Fie}]\label{Eprop-6} Each~$\ell$-facet of an acute 0/1-$k$-simplex is an acute 0/1-$\ell$-simplex.
\end{Pro}
This corresponds to the well-known linear algebraic statement that the inverse of 
each principal~$k\times k$ submatrix of~$G$ is also a diagonally 
dominant strictly Stieltjes matrix. Together with Lemma~\ref{Elem6}, this proves the following.
\begin{Co} \label{Ecor-6} Let~$P^\ast$ be the minimal matrix representation of an acute simplex~$c\in\PP_n^k$. Then
for each~$j\in\{1,\dots,k-1\}$, the submatrix~$P_j^\ast$ of~$P^\ast$ consisting of its~$j$ leftmost columns is a
minimal matrix representation of an acute simplex~$c_j\in\PP_n^j$.
\end{Co}
Corollary~\ref{Ecor-6} shows in particular that the first three columns of any 
minimal matrix representation of an acute 0/1-simplex form a minimal matrix representation of an acute 0/1-triangle.
\begin{Def}[Acute extensions of~$S$]{\rm Let~$S\subset I^n$ be an \index{acute extensions}
acute 0/1-simplex with~$k\leq n$ vertices. The set~$\AA^n(S)$ of 
{\em acute extensions} of~$S$ consists of all~$v\in\BB^n$ such that~$\conv(S,v)$ is an acute 0/1-simplex with~$k+1$ vertices.}
\end{Def}
The following classical result formulates a necessary condition 
for a vertex~$v\in\BB^n$ of~$I^n$ to be an element of the set~$\AA^n(S)$ just defined.
\begin{Pro}[\cite{Fie}] \label{Epro-9} Let~$S$ be an acute~$n$-simplex. 
Then each vertex of~$S$ projects orthogonally into the interior of its opposite~$(n\!-\!1)$-dimensional facet.
\end{Pro} 
\begin{Def}[Candidate acute extensions of~$S$]\label{Edef-cand}{\rm Let~$S\subset I^n$ 
be an acute 0/1-simplex with~$k\leq n$ vertices. The set~$\CC^n(S)$ of 
{\em candidate acute extensions}\index{candidate acute extensions} of~$S$ consists of all~$v\in\BB^n$ such that~$v$ projects orthogonally into the interior of~$S$.} 
\end{Def}

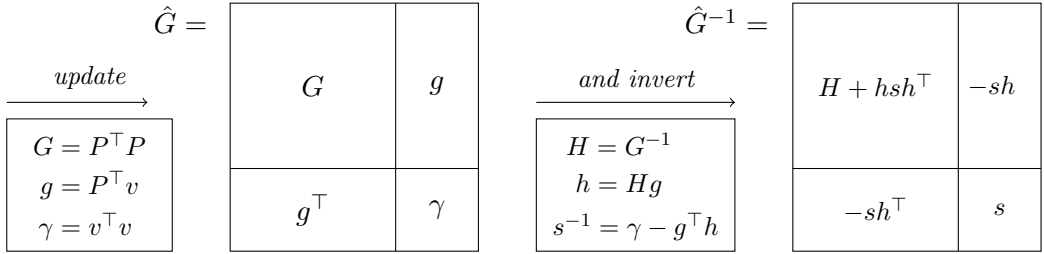
\begin{figure}[h]
\begin{center}
\begin{tikzpicture}[scale=1.1]
\begin{scope}[shift={(0,-0.7)}]
\draw (0,1.7)--(2,1.7)--(2,3.3)--(0,3.3)--cycle;
\node[scale=0.9] at (1,3) {$G=P^\top P$};
\node[scale=0.9] at (1,2.5) {$g=P^\top v$};
\node[scale=0.9] at (0.95,2) {$\gamma=v^\top v$};
\draw[->] (0,3.5)--(1.7,3.5);
\node[scale=0.9] at (1,3.8) {\em update};
\node at (2.1,4.5) {$\hat{G}=$};
\end{scope}
\begin{scope}[shift={(2.7,0)}]
\draw (0,1)--(3,1)--(3,4)--(0,4)--cycle;
\draw (2,1)--(2,4);
\draw (0,2)--(3,2);
\node at (1,3) {$G$};
\node at (1,1.5) {$g^\top$};
\node at (2.5,3) {$g$};
\node at (2.5,1.5) {$\gamma$};
\end{scope}
\begin{scope}[shift={(6.4,-0.7)}]
\draw (0,1.7)--(2.4,1.7)--(2.4,3.3)--(0,3.3)--cycle;
\node[scale=0.9] at (1,3) {$H=G^{-1}$};
\node[scale=0.9] at (1,2.5) {$h=Hg$};
\node[scale=0.9] at (1.2,2) {$s^{-1}=\gamma-g^\top h$};
\draw[->] (0,3.5)--(2.4,3.5);
\node[scale=0.9] at (1.2,3.8) {\em and invert};
\node at (2.3,4.5) {$\hat{G}^{-1}=$};
\end{scope}
\begin{scope}[shift={(9.5,0)}]
\draw (0,1)--(3,1)--(3,4)--(0,4)--cycle;
\draw (2,1)--(2,4);
\draw (0,2)--(3,2);
\node[scale=0.9] at (1,3) {$H+hsh^\top$};
\node[scale=0.9] at (1,1.5) {$-sh^\top$};
\node[scale=0.9] at (2.4,3) {$-sh$};
\node[scale=0.9] at (2.5,1.5) {$s$};
\end{scope}
\end{tikzpicture}
\end{center} 
\caption{\small{Updating the inverse of an updated Gramian.}}
\label{figureE13}
\end{figure}

\begin{rem}{\rm Due to Proposition~\ref{Epro-9}, we have that~$\AA^n(S)\subset \CC^n(S)$. 
The sets are in general not equal. This can be seen in Figure~\ref{figureE12}. For each of the acute 
triangular facets~$\TT$ of the antipodal 0/1-tetrahedron in~$I^3$, there exist {\em two} 
vertices of~$I^3$ that project in the interior of~$\TT$, but only the convex hull of~$\TT$ 
with {\em one} of them yields an acute tetrahedron.}
\end{rem}    
We will now investigate on a linear algebraic level when~$v\in\CC^n(S)$ or 
even~$v\in \AA^n(S)$. For this, let~$P\in\BB^{n\times k}$ with~$1\leq k<n$  
be such that~$G=P^\top P$ satisfies the conditions (1) and (2) in 
Definition~\ref{Edef-acute}, and let~$v\in\BB^n$. Consider the matrix 
$[P|v]$. Its Gramian~$\hat{G}$ is a simple update of~$G$. Also its inverse 
$\hat{G}^{-1}$ is an update of the inverse~$H$ of~$G$, as depicted in Figure~\ref{figureE13}.
\\[2mm]
Note that~$\hat{G}$ is positive semi-definite. It is invertible if and only 
if~$s>0$. This condition will turn out to be automatically fulfilled if~$v$ projects in the interior of its opposite facet.
\begin{Le}\label{Elem-9} The vertex~$v$ is an element of~$\in\CC^n(S)$ if and only if~$h>0$ and~$e^\top h < 1$.
\end{Le}
{\bf Proof}.  Observe that the orthogonal projection of~$v$ on the column span of~$P$ equals~$Ph$, because
\be Ph = P(P^\top P)^{-1}P^\top v. \ee
To lie in the {\em interior} of the corresponding facet of~$S$,~$Ph$ must be 
a convex combination in which all vertices of that facet, including the origin, 
participate nontrivially. Thus, the entries of~$h$ must be positive and up to less than one.\hfill~$\Box$
\begin{Co} If~$h>0$ and~$e^\top h<1$ then~$s>0$. 
\end{Co} 
{\bf Proof.} If~$v$ projects in the interior of its opposite facet, then in particular, 
$v$ is not equal to a vertex of that facet. Also, no vertex of~$I^n$ is a convex 
combination of any of the others. Thus, the convex hull of the facet and~$v$ has nonzero volume.\hfill~$\Box$

\smallskip

The diagram in Figure~\ref{figureE13} and Lemma~\ref{Elem-9} again show that even if~$v\in\CC^n(S)$, 
it does not need to be in~$\AA^n(S)$. Indeed, because~$hsh^\top >0$, the updated 
matrix~$H+hsh^\top$ may have nonnegative off-diagonal entries and violate condition 
(1) in Definition~\ref{Edef-acute}. Moreover, condition (2) may also be violated, as the row sums of~$\hat{G}^{-1}$ equal
\be\label{Erow-sums}  \left[\begin{array}{r}r\\\rho\end{array}\right]  = 
\hat{G}^{-1}\left[\begin{array}{r}e\\1\end{array}\right] = 
\left[\begin{array}{c} He + sh(h^\top e-1) \\ s(1-h^\top e)\end{array}\right], \ee
and if~$v\in\CC^n(S)$ then according to Lemma~\ref{Elem-9},~$h^\top e-1<0$. 
Although this implies that~$\rho>0$, some of the remaining entries of~$r$ may be negative, in spite of~$He>0$.

\smallskip

Suppose now that~$[P|v]$ is indeed such, that~$v\in\CC^n(S)$ but that~$v\not\in\AA^n(S)$. For some~$w\in\BB^m$, consider the matrix
\[ \left[\begin{array}{r|r} P & v \\ \hline 0 & w\end{array}\right]\hdrie\mbox{\rm with Gramian}\hdrie 
\tilde{G} = \left[\begin{array}{r|r} P^\top P & P^\top v\\ \hline 
{v^\top P}^{\phantom{\!\!\!\!\!X}}& v^\top v+w^\top w\end{array}\right]. \]
In comparison with the Gramian~$\hat{G}$ in Figure~\ref{figureE13}, only the bottom right 
entry has changed. Obviously, it~$w^\top w$ is large enough, the corresponding 
value of~$s$ will decrease so much, that the off-diagonal entries of~$H+hsh^\top$ 
are negative, and the row sums in (\ref{Erow-sums}) positive.\\[2mm]
In other words, if a vertex~$v$ projects in the interior of an acute facet~$F$, then 
by moving~$v$ orthogonally away from~$F$, the simplex~$\conv(F,v)$ will ultimately always become acute.\\[2mm]
The above discussion proves the following theorem. 
\begin{Th}\label{Eth-7} Let~$S$ be an acute 0/1-simplex in 
$I^n$. Consider~$I^n$ as a facet of~$I^{n+m}$. Let the first~$n$ 
entries of~$v\in\BB^{n+m}$ correspond to the vertices of~$I^n$. Then:
\be \label{Eeq-72-1}  \CC^{n+m}(S)  = \left\{ \left[\begin{array}{c} v \\ 
w\end{array}\right] \sth v \in \CC^{n}(S),\hdrie w\in\BB^m\right\}. \ee
and
\be \label{Eeq-72-2}   \AA^{n+m}(S)  \supset \left\{ \left[\begin{array}{c} v \\ 
w\end{array}\right] \sth v \in \AA^{n}(S),\hdrie w\in\BB^m\right\}. \ee
Moreover, for each~$v\in\CC^n(S)$ there exists an~$\ell$ such that
\be \label{Eeq-72-3}  w^\top w \geq \ell \hdrie\Leftrightarrow \hdrie 
\left[\begin{array}{c} v \\ w\end{array}\right]\in\AA^{n+m}(S), \ee
provided that~$m$ is large enough.
\end{Th}   
{\bf Proof. } Statement (\ref{Eeq-72-1}) follows because the right-hand side consists 
precisely of those vertices of~$I^{n+m}$ whose orthogonal projection on~$I^n$ land 
in~$\CC^n(S)$. The claims in (\ref{Eeq-72-2}) and (\ref{Eeq-72-3}) follow from the above discussion.\hfill~$\Box$

\smallskip

Note that the optimal value of~$\ell$ in (\ref{Eeq-72-3}) can, in principle, be 
computed as soon as the data~$H,h,g$ and~$\gamma$ are available, as is visible from Figure~\ref{figureE13}.

\smallskip

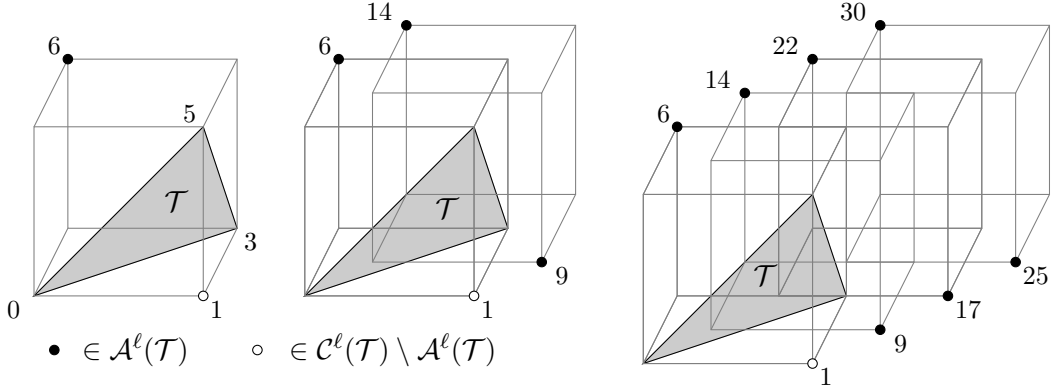
\begin{figure}[h]
\begin{center}
\begin{tikzpicture}[scale=0.9]
\begin{scope}
\draw[gray!40!white, fill=gray!40!white] (0,0)--(3,1)--(2.5,2.5)--cycle;
\node at (2.1,1.4) {$\mathcal{T}$};
\node[scale=0.9] at (-0.3,-0.2) {$0$};
\node[scale=0.9] at (2.7,-0.2) {$1$};
\node[scale=0.9] at (3.2,0.8) {$3$};
\node[scale=0.9] at (2.3,2.7) {$5$};
\node[scale=0.9] at (0.3,3.7) {$6$};
\draw[fill=black] (0.3,-0.8) circle [radius=0.07];
\draw[fill=white] (3.3,-0.8) circle [radius=0.07];
\node at (1.5,-0.8) {$\in\mathcal{A}^\ell(\mathcal{T})$};
\node at (5.3,-0.8) {$\in\mathcal{C}^\ell(\mathcal{T})\setminus \mathcal{A}^\ell(\mathcal{T})$};
\draw (0,0)--(3,1)--(2.5,2.5)--cycle;
\draw[gray] (0,0)--(2.5,0)--(2.5,2.5)--(0,2.5)--cycle;
\draw[gray] (0.5,1)--(3,1)--(3,3.5)--(0.5,3.5)--cycle;
\draw[gray] (0,0)--(0.5,1);
\draw[gray] (2.5,0)--(3,1);
\draw[gray] (2.5,2.5)--(3,3.5);
\draw[gray] (0,2.5)--(0.5,3.5);
\draw[fill=white] (2.5,0) circle [radius=0.07];
\draw[fill=black] (0.5,3.5) circle [radius=0.07];
\end{scope}
\begin{scope}[shift={(4,0)}]
\draw[gray!40!white, fill=gray!40!white] (0,0)--(3,1)--(2.5,2.5)--cycle;
\node at (2.1,1.3) {$\mathcal{T}$};
\node[scale=0.9] at (2.7,-0.2) {$1$};
\node[scale=0.9] at (3.8,0.3) {$9$};
\node[scale=0.9] at (0.3,3.7) {$6$};
\node[scale=0.9] at (1.1,4.2) {$14$};
\draw (0,0)--(3,1)--(2.5,2.5)--cycle;
\draw[gray] (0,0)--(2.5,0)--(2.5,2.5)--(0,2.5)--cycle;
\draw[gray] (0.5,1)--(3,1)--(3,3.5)--(0.5,3.5)--cycle;
\draw[gray] (0,0)--(0.5,1);
\draw[gray] (2.5,0)--(3,1);
\draw[gray] (2.5,2.5)--(3,3.5);
\draw[gray] (0,2.5)--(0.5,3.5);
\begin{scope}[shift={(1,0.5)}]
\draw[gray] (0,0)--(2.5,0)--(2.5,2.5)--(0,2.5)--cycle;
\draw[gray] (0.5,1)--(3,1)--(3,3.5)--(0.5,3.5)--cycle;
\draw[gray] (0,0)--(0.5,1);
\draw[gray] (2.5,0)--(3,1);
\draw[gray] (2.5,2.5)--(3,3.5);
\draw[gray] (0,2.5)--(0.5,3.5);
\end{scope}
\draw[gray] (0,0)--(2.5,0)--(2.5,2.5)--(0,2.5)--cycle;
\draw[gray] (0.5,1)--(3,1)--(3,3.5)--(0.5,3.5)--cycle;
\draw[gray] (0,0)--(0.5,1);
\draw[gray] (2.5,0)--(3,1);
\draw[gray] (2.5,2.5)--(3,3.5);
\draw[gray] (0,2.5)--(0.5,3.5);
\draw[fill=white] (2.5,0) circle [radius=0.07];
\draw[fill=black] (0.5,3.5) circle [radius=0.07];
\draw[fill=black] (3.5,0.5) circle [radius=0.07];
\draw[fill=black] (1.5,4) circle [radius=0.07];
\end{scope}
\begin{scope}[shift={(9,-1)}]
\draw[gray!40!white, fill=gray!40!white] (0,0)--(3,1)--(2.5,2.5)--cycle;
\node at (1.8,1.3) {$\mathcal{T}$};
\node[scale=0.9] at (2.7,-0.2) {$1$};
\node[scale=0.9] at (3.8,0.3) {$9$};
\node[scale=0.9] at (4.8,0.8) {$17$};
\node[scale=0.9] at (5.8,1.3) {$25$};
\node[scale=0.9] at (0.3,3.7) {$6$};
\node[scale=0.9] at (1.1,4.2) {$14$};
\node[scale=0.9] at (2.1,4.7) {$22$};
\node[scale=0.9] at (3.1,5.2) {$30$};
\draw (0,0)--(3,1)--(2.5,2.5)--cycle;
\draw[gray] (0,0)--(2.5,0)--(2.5,2.5)--(0,2.5)--cycle;
\draw[gray] (0.5,1)--(3,1)--(3,3.5)--(0.5,3.5)--cycle;
\draw[gray] (0,0)--(0.5,1);
\draw[gray] (2.5,0)--(3,1);
\draw[gray] (2.5,2.5)--(3,3.5);
\draw[gray] (0,2.5)--(0.5,3.5);
\begin{scope}[shift={(1,0.5)}]
\draw[gray] (0,0)--(2.5,0)--(2.5,2.5)--(0,2.5)--cycle;
\draw[gray] (0.5,1)--(3,1)--(3,3.5)--(0.5,3.5)--cycle;
\draw[gray] (0,0)--(0.5,1);
\draw[gray] (2.5,0)--(3,1);
\draw[gray] (2.5,2.5)--(3,3.5);
\draw[gray] (0,2.5)--(0.5,3.5);
\end{scope}
\draw[gray] (0,0)--(2.5,0)--(2.5,2.5)--(0,2.5)--cycle;
\draw[gray] (0.5,1)--(3,1)--(3,3.5)--(0.5,3.5)--cycle;
\draw[gray] (0,0)--(0.5,1);
\draw[gray] (2.5,0)--(3,1);
\draw[gray] (2.5,2.5)--(3,3.5);
\draw[gray] (0,2.5)--(0.5,3.5);
\draw[fill=white] (2.5,0) circle [radius=0.07];
\draw[fill=black] (0.5,3.5) circle [radius=0.07];
\draw[fill=black] (3.5,0.5) circle [radius=0.07];
\draw[fill=black] (1.5,4) circle [radius=0.07];
\end{scope}
\begin{scope}[shift={(11,0)}]
\draw[gray] (0,0)--(2.5,0)--(2.5,2.5)--(0,2.5)--cycle;
\draw[gray] (0.5,1)--(3,1)--(3,3.5)--(0.5,3.5)--cycle;
\draw[gray] (0,0)--(0.5,1);
\draw[gray] (2.5,0)--(3,1);
\draw[gray] (2.5,2.5)--(3,3.5);
\draw[gray] (0,2.5)--(0.5,3.5);
\begin{scope}[shift={(1,0.5)}]
\draw[gray] (0,0)--(2.5,0)--(2.5,2.5)--(0,2.5)--cycle;
\draw[gray] (0.5,1)--(3,1)--(3,3.5)--(0.5,3.5)--cycle;
\draw[gray] (0,0)--(0.5,1);
\draw[gray] (2.5,0)--(3,1);
\draw[gray] (2.5,2.5)--(3,3.5);
\draw[gray] (0,2.5)--(0.5,3.5);
\end{scope}
\draw[gray] (0,0)--(2.5,0)--(2.5,2.5)--(0,2.5)--cycle;
\draw[gray] (0.5,1)--(3,1)--(3,3.5)--(0.5,3.5)--cycle;
\draw[gray] (0,0)--(0.5,1);
\draw[gray] (2.5,0)--(3,1);
\draw[gray] (2.5,2.5)--(3,3.5);
\draw[gray] (0,2.5)--(0.5,3.5);
\draw[fill=black] (2.5,0) circle [radius=0.07];
\draw[fill=black] (0.5,3.5) circle [radius=0.07];
\draw[fill=black] (3.5,0.5) circle [radius=0.07];
\draw[fill=black] (1.5,4) circle [radius=0.07];
\end{scope}
\end{tikzpicture}
\end{center} 
\caption{\small{Impression of the structure of the sets~$\AA^\ell(\TT)$ and 
$\CC^\ell(\TT)$  for increasing values of~$\ell$. The white vertex is not in 
$\AA^\ell(\TT)$ but the ones ``orthogonally above'' it, are in~$\AA^\ell(\TT)$. }}
\label{figureE14}
\end{figure}

Figure \ref{figureE14} serves to illustrate the claims of Theorem~\ref{Eth-7}. Consider the acute 
0/1-triangle~$\TT$ with vertices~$0,3,5$, in the numbering of Figure~\ref{figureE14}. The set 
$\CC^3(\TT)$ consists of vertices~$1$ and~$6$, as both project in the interior of 
$\TT$. Only vertex~$6$ is an element of~$\AA^3(\TT)$. Indeed, the tetrahedron 
formed by~$1$ and~$\TT$ is not acute: it has {\em right} dihedral angles. 
However, each of the vertices~$9,17,25$ {\em orthogonally above} vertex~$1$ 
forms an acute 0/1-tetrahedron with~$\TT$ and thus belong to~$\AA^5(\TT)$, 
as do the ones~$14,22,30$ orthogonally above vertex~$6$.
\\[2mm]
The value of Theorem~\ref{Eth-7} is that in order to determine the set~$\AA^n(S)$ 
of a given minimal representative of an acute 0/1-simplex~$S$, the computational work can be reduced to:\\[2mm]
$\bullet$ find the smallest~$k\leq n$ for which~$S\in I^k$ and determine~$\CC^k(S)$;\\[2mm]
$\bullet$ determine which~$v\in\CC^k(S)$ are in~$\AA^k(S)$;\\[2mm]
$\bullet$ for each~$v\in\CC^k(S)\setminus\AA^k(S)$, determine the value of~$\ell$ in (\ref{Eeq-72-3}).\\[2mm]
After doing so, all remaining vertices of~$v\in\BB^n$ that are in~$\AA^n(S)$ 
can now be easily listed without having to verify acuteness of the simplex~$\conv(S,v)$.

\smallskip

The next theorem is not difficult, but will play a central role in the enumeration problem.
\begin{Th}\label{Eth-6} Let~$S\subset I^n$ be an acute 0/1-simplex with 
$k$ vertices. If~$\hat{S}$ is an acute 0/1-simplex in~$I^n$ having~$S$ 
as a facet, then each vertex of~$\hat{S}$ belongs to~$S$ or to~$\AA^n(S)$.
\end{Th}
{\bf Proof.} Let~$v$ be a vertex of~$\hat{S}$ that does not belong to~$S$. 
Then~$\conv(S,v)$ is a facet of~$\hat{S}$. Since~$\hat{S}$ is acute, 
Proposition~\ref{Eprop-6} shows that~$\conv(S,v)$ is acute, and thus~$v\in\AA^n(S)$.\hfill~$\Box$

\smallskip

In the language of Figure \ref{figureE14}, Theorem~\ref{Eth-6} expresses that 
each acute 0/1-simplex having~$\TT$ as a triangular facet, has 
all its vertices amongst the black bullets. Note that for each of 
these black bullets, its projection on the triangle~$\TT$ is the orthocenter of~$\TT$. 
\begin{Co}\label{Ecor-7} Let~$\TT$ be a facet of an acute 
0/1-simplex~$S$. Then~$\AA^n(S)\subset\AA^n(\TT)$.
\end{Co} 
Theorem~\ref{Eth-6} shows the importance of administrating the set of 
acute extensions in the process of building acute 0/1-simplices from 
the starting point of a minimal representative an acute 0/1-triangle 
$\TT$. Adding vertices to~$\TT$, the set of acute extensions of the 
resulting simplices becomes smaller and smaller as the dimension of the 
simplex become larger, hence reducing the amount of work to be done 
to build all minimal representatives of 0/1-simplices having~$\TT$ as minimal triangular facet.

\smallskip

It remains necessary to work with the concept of minimal matrix representations, 
to reduce the amount of data to be computed. Not only {\em after} the construction process but also {\em during}.
%%%%%%%%%%%%%%%%%%%%%%%%
%%%%%%%%%%%%%
\subsection{Minimal acute extensions of acute 0/1-simplices}
Let~$\TT^\ast$ be a minimal representative of an acute 0/1-triangle 
in~$I^n$ with minimal matrix representation~$T^\ast$.  Thus, the vertices 
of~$\TT^\ast$ are the column vectors of~$T^\ast$. Now, consider the~$p$ 
matrices of size~$n\times 4$ defined by
\be\label{Eminimals} [T^\ast| t_1], \hdrie \dots\hdrie , \hdrie [T^\ast| t_p], \ee 
where~$\AA^n(\TT^\ast)=\{t_1,\dots,t_p\}$ is the set of acute extensions of~$\TT^\ast$.
\begin{Pro} Each minimal matrix representation of each acute 0/1-tetrahedron having 
$\TT^\ast$ as minimal triangular facet, is among the matrices in (\ref{Eminimals}).
\end{Pro}
{\bf Proof. } This follows immediately from Corollary~\ref{Ecor-6} and the fact that 
there are no other~$t\in I^n$ such that~$[T^\ast| t]$ is an acute 0/1-tetrahedron. \hfill~$\Box$

\smallskip

Consequently, we can subdivide the set of acute extensions of a minimal representative of an acute 
0/1-simplex into a minimal and a non-minimal part.
\begin{Def}[Minimal acute extensions]\label{Emae}{\rm The set \index{Minimal acute extension}
$\AA^n(\TT^\ast)$ of acute extensions of a minimal representative 
$\TT^\ast$ of an acute 0/1-simplex with~$k\leq n$ vertices, with 
minimal matrix representation~$T^\ast$, is subdivided as
\[ \AA_\ast^n(\TT^\ast) = \{ t \in \AA^n(\TT^\ast) \sth [T^\ast|t] \hdrie
\mbox{\rm is a minimal matrix representation}\}, \]
and its complement~$\AA^n_\circ(\TT^\ast)$ in~$\AA^n(\TT^\ast)$.}
\end{Def} 
The results of Section~\ref{ESect-3} immediately show that the following 
matrices in (\ref{Eminimals}) are in~$\AA^n_\circ(\TT^\ast)$:\\[2mm]
$(1)$ the ones for which the {\em column number}~$v_n^\top t_j$ of 
$t_j$ is {\em smaller} than~$v_n^\top T^\ast e_3$;\\[2mm]
$(2)$ the ones for which the {\em row numbers}~$[T^\ast| t_j]w_4$ 
are {\em not} non-increasing.\\[2mm]
To make the subdivision of~$\AA^n(\TT^\ast)$ in (\ref{Eminimals}) into 
$\AA_\ast^n(\TT^\ast)$ and~$\AA_\circ^n(\TT^\ast)$ complete, one 
may use an adapted version of Algorithm 3; adapted in the sense that 
it should be aborted as soon as a matrix representation is encountered 
that proves that~$[T^\ast| t_j]$ is not minimal.
\begin{rem}{\rm Note that if
\[   [T^\ast| t_j] = \left[\begin{array}{r|r} T^\ast & t_j^1 \\ \hline 
0 & t_j^2 \end{array}\right] \hdrie\mbox{\rm with} \hdrie t_j^2\in\BB^m,\]
then~$[T^\ast| t_j]w_4$ is not non-increasing if~$t_j^2$ itself is not 
non-increasing. Thus, for a number of~$t\in\AA^n(\TT^\ast)$ it may 
be directly indicated that they do not belong to~$t\in\AA_\ast^n(\TT^\ast)$.}
\end{rem}
Suppose that it has been established that the matrix~$[T^\ast| t_\ell]$ is a minimal 
matrix representation of a 0/1-tetrahedron with minimal representative 
$\hat{\TT^\ast}$. In order to continue the construction process of acute 
0/1-simplices efficiently, either in a {\em depth-first} or a {\em breadth-first} 
fashion, the data structures of acute candidates and acute extensions of~$\hat{\TT^\ast}$ need to be updated. 
\begin{rem}{\em It may happen that while~$t$ is not a minimal acute 
extension of some acute 0/1-simplex~$\TT^\ast$, it {\em is} indeed a 
minimal acute extension of an acute simplex having~$\TT^\ast$ as minimal facet. Indeed, let
\be \small  T^\ast = \left[\begin{array}{rrr} 0 & 1 & 1 \\0 & 1 & 0 \\ 
0 & 0 & 1 \\ 0 & 0 & 0\\ 0 & 0 & 0\end{array}\right], \hdrie t_1=\left[\begin{array}{r} 1 \\ 
0 \\ 0 \\ 1 \\ 0\end{array}\right] \und t_2=\left[\begin{array}{r} 1 \\ 0 \\ 0 \\ 0 \\ 1\end{array}\right]. \ee
Then both~$t_1$ and~$t_2$ are acute extensions of~$T^\ast$. Only~$t_1$ is a minimal 
acute extension of~$T^\ast$, whereas~$t_2$ is not due to criterion (2) above. But~$t_2$ is 
a minimal acute extension of the minimal matrix representation~$[T^\ast|t_1]$.}
\end{rem}
Due to Corollary~\ref{Ecor-7} we have that~$\AA^n(\hat{\TT}^\ast) \subset \AA^n(\TT^\ast)$. 
To determine~$\AA^n(\hat{\TT}^\ast)$ {\em exactly}, it may not be necessary to verify for 
{\em each}~$t\in\AA^n(\TT^\ast)$ whether the convex hull~$\conv(\hat{\TT}^\ast,t)$ of~$\hat{\TT}^\ast$ 
and~$t$ is an acute 0/1-simplex. Indeed, if~$\hat{\TT}^\ast\subset I^k$ for some~$k<n$, it suffices to 
find out which~$t\in\AA^k(\TT^\ast)$ are in~$\CC^k(\hat{\TT}^\ast)$ and which of these are in 
$\AA^k(\hat{\TT}^\ast)$ and then use Theorem~\ref{Eth-7}. 
 \begin{Ex} The matrix~$T^\ast\in\BB^{6\times 3}$ in (\ref{Eexample}) with vertex numbers 
$0,3$ and~$13$, is a minimal matrix representation of a minimal representative~$\TT^\ast$ of an 
acute 0/1-triangle. Clearly~$\TT^\ast\subset I^4$. The vertices from~$I^4$ that are in~$\CC^4(\TT^\ast)$ 
are listed by their vertex numbers~$1,5,6,9,10$ and~$14$ as well as their 0/1-vectors. The numbers in 
the row indicated by~$\AA^4(\TT^\ast)$ correspond to the smallest value~$m$ of entries equal to one 
need to be appended to the vector above it such that it becomes an element of~$\AA^{n+m}(\TT^\ast)$. 
\be\label{Eexample}\small T^\ast =\begin{array}{|r|r|r|}\hline 0 & 1 & 1 \\ 0 & 1 & 0 \\ 0 & 0 & 1 \\ 
0 & 0 & 1\\ 0 & 0 & 0 \\ 0 & 0 & 0\\ \hline\hline 0 & 3 & 13 \\ \hline \end{array} \hdrie \hdrie
\begin{array}{|r||r|r|r|r|r|r|} 
\hline
\CC^4(\TT^\ast) & 1 & 5 & 6 & 9 & 10 & 14\\
\hline
\hline 
& 1  &   1  &   0  &   1  &   0  &   0 \\
& 0  &   0  &   1  &   0  &   1  &   1 \\
& 0  &   1  &   1  &   0  &   0  &   1 \\
& 0  &   0  &   0  &   1  &   1  &   1 \\
\hline
\hline
\AA^4(\TT^\ast) & 1  &   1  &   0  &   1  &   0  &   0 \\
\hline
\hline
\AA^5(\TT^\ast) & 17 & 21 & 22 & 25 & 26 & 30\\
\hline
\hline
\AA^6(\TT^\ast) & 33 & 37 & 38 & 41 & 42 & 46 \\
                & 49 & 53 & 54 & 57 & 58 & 62\\
\hline
\end{array}\ee 
Thus, the vertices~$6,10$ and~$14$ are in~$\AA^4(\TT^\ast)$, and the vertices~$1,5$ and~$9$ are 
candidates that need only one additional~$1$ to become acute extensions. This is visible in the next 
row, where the vertex numbers of {\em additional} vertices in~$\AA^5(\TT^\ast)$ are displayed, 
which are the ones from~$\CC^4(\TT^\ast)$ plus~$2^4=16$. Finally, the elements from 
$\AA^6(\TT^\ast)$ {\em additional} to the ones from~$\AA^4(\TT^\ast)$ and~$\AA^5(\TT^\ast)$ 
are precisely those with~$2^5=32$ added. It can easily be verified that only~$14, 17, 21, 22, 30, 49, 53, 54, 62$ 
remain after removing the ones that fall under the items (1) or (2) below Definition~\ref{Emae}. It is also 
not hard to see that the matrix~$\hat{T}^\ast$ with column numbers~$0,3,13,21$ is indeed a 
minimal matrix representation of a 0/1-tetrahedron~$\hat{\TT}^\ast\subset I^5$. To determine 
$\AA^6(\hat{\TT}^\ast)$, we re-investigate the each of the vertices~$1, 5, 6, 9, 10, 14, 17, 21, 22, 25, 26, 30$ 
and indicate whether it belongs to~$\CC^5(\hat{\TT}^\ast)$ or~$\AA^5(\hat{\TT}^\ast)$ or neither. This 
determines which of the vertices in~$33, 37, 38, 41, 42, 46, 49, 53, 54, 57$,~$58, 62$ are {\em additionally} 
in~$\AA^6(\hat{\TT}^\ast)$.~\hfill~$\diamondsuit$
\end{Ex}
%%%%%%%%%
%%%%%%%%%%%
%%%%%%%%%%%%%%
\section{A special class of acute 0/1-simplices}\label{ESect-6}
Here we analyze the computational results of the codes presented in 
Section~\ref{ESect-7}. Looking at the structure of the 0/1-matrices 
presented there, we observe some patterns. Although not all patters can 
be mathematically accounted for, there is one pattern that can be fully explained.
\begin{rem}{\rm Since each first column of a minimal matrix representation 
is zero according to Corollary~\ref{Eco5}, we will omit this redundant column 
from the notation. What remains is a square matrix, whose Gramian has an 
inverse that is a diagonally dominant strictly Stieltjes matrix, and which 
we will also call a minimal matrix representation. See Definition~\ref{Edef-acute}.}
\end{rem}
\subsection{Acute simplices (upper Hessenberg matrix representations)}\label{Esect-7.1}
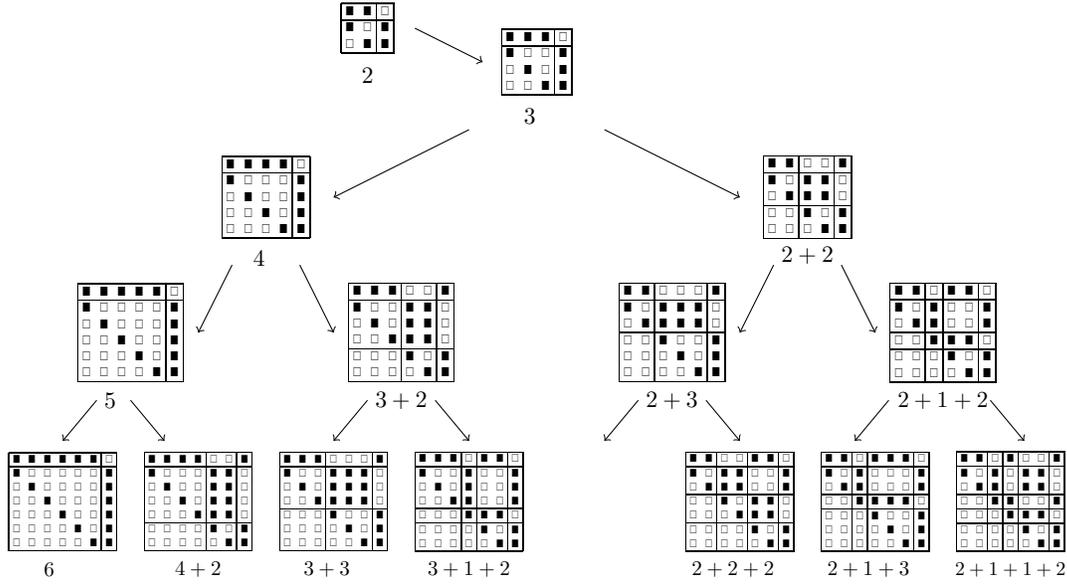
\begin{figure}[h]
\begin{center}
\begin{tikzpicture}[scale=0.9, every node/.style={scale=0.9}]
\node[scale=0.9] at (4.5,-2.2) {$2$};
\node[xscale=0.4, yscale=0.5] at (4.5,-1.5) {$\begin{array}{|cc|c|}
\hline
\blacksquare&\blacksquare&\square\\\hline
\blacksquare&\square&\blacksquare\\
\square&\blacksquare&\blacksquare\\\hline
\end{array}$};
\draw[->] (5.2,-1.5)--(6.2,-2);
\node[scale=0.9] at (6.9,-2.8) {$3$};
\node[xscale=0.4, yscale=0.5] at (7,-2) {$\begin{array}{|ccc|c|}
\hline
\blacksquare&\blacksquare&\blacksquare&\square\\\hline
\blacksquare&\square&\square&\blacksquare\\
\square&\blacksquare&\square&\blacksquare\\
\square&\square&\blacksquare&\blacksquare\\\hline
\end{array}$};
\draw[->] (6,-3)--(4,-4);
\draw[->] (8,-3)--(10,-4);
\node[scale=0.9] at (2.9,-4.9) {$4$};
\node[xscale=0.4, yscale=0.5] at (3,-4) {$\begin{array}{|cccc|c|}
\hline
\blacksquare&\blacksquare&\blacksquare&\blacksquare&\square\\\hline
\blacksquare&\square&\square&\square&\blacksquare\\
\square&\blacksquare&\square&\square&\blacksquare\\
\square&\square&\blacksquare&\square&\blacksquare\\
\square&\square&\square&\blacksquare&\blacksquare\\\hline
\end{array}$};
\node[scale=0.9] at (11,-4.85) {$2+2$};
\node[xscale=0.4, yscale=0.5] at (11,-4) {$\begin{array}{|cc|cc|c|}
\hline
\blacksquare&\blacksquare&\square&\square&\blacksquare\\\hline
\blacksquare&\square&\blacksquare&\blacksquare&\square\\
\square&\blacksquare&\blacksquare&\blacksquare&\square\\\hline
\square&\square&\blacksquare&\square&\blacksquare\\
\square&\square&\square&\blacksquare&\blacksquare\\\hline
\end{array}$};
\draw[->] (2.5,-5)--(2,-6);
\draw[->] (3.5,-5)--(4,-6);
\draw[->] (10.5,-5)--(10,-6);
\draw[->] (11.5,-5)--(12,-6);
\node[scale=0.9] at (0.7,-7) {$5$};
\node[xscale=0.4, yscale=0.5] at (1,-6) {$\begin{array}{|ccccc|c|}
\hline
\blacksquare&\blacksquare&\blacksquare&\blacksquare&\blacksquare&\square\\\hline
\blacksquare&\square&\square&\square&\square&\blacksquare\\
\square&\blacksquare&\square&\square&\square&\blacksquare\\
\square&\square&\blacksquare&\square&\square&\blacksquare\\
\square&\square&\square&\blacksquare&\square&\blacksquare\\
\square&\square&\square&\square&\blacksquare&\blacksquare\\
\hline
\end{array}$};
\node[scale=0.9] at (5,-7) {$3+2$};
\node[xscale=0.4, yscale=0.5] at (5,-6) {$\begin{array}{|ccc|cc|c|}
\hline
\blacksquare&\blacksquare&\blacksquare&\square&\square&\blacksquare\\\hline
\blacksquare&\square&\square&\blacksquare&\blacksquare&\square\\
\square&\blacksquare&\square&\blacksquare&\blacksquare&\square\\
\square&\square&\blacksquare&\blacksquare&\blacksquare&\square\\\hline
\square&\square&\square&\blacksquare&\square&\blacksquare\\
\square&\square&\square&\square&\blacksquare&\blacksquare\\
\hline
\end{array}$};
\node[scale=0.9] at (9,-7) {$2+3$};
\node[xscale=0.4, yscale=0.5] at (9,-6) {$\begin{array}{|cc|ccc|c|}
\hline
\blacksquare&\blacksquare&\square&\square&\square&\blacksquare\\\hline
\blacksquare&\square&\blacksquare&\blacksquare&\blacksquare&\square\\
\square&\blacksquare&\blacksquare&\blacksquare&\blacksquare&\square\\\hline
\square&\square&\blacksquare&\square&\square&\blacksquare\\
\square&\square&\square&\blacksquare&\square&\blacksquare\\
\square&\square&\square&\square&\blacksquare&\blacksquare\\
\hline
\end{array}$};
\node[scale=0.88] at (13,-7) {$2+1+2$};
\node[xscale=0.4, yscale=0.5] at (13,-6) {$\begin{array}{|cc|c|cc|c|}
\hline
\blacksquare&\blacksquare&\square&\blacksquare&\blacksquare&\square\\\hline
\blacksquare&\square&\blacksquare&\square&\square&\blacksquare\\
\square&\blacksquare&\blacksquare&\square&\square&\blacksquare\\\hline
\square&\square&\blacksquare&\blacksquare&\blacksquare&\square\\\hline
\square&\square&\square&\blacksquare&\square&\blacksquare\\
\square&\square&\square&\square&\blacksquare&\blacksquare\\
\hline
\end{array}$};
\draw[->] (0.5,-7)--(0,-7.6);
\draw[->] (1,-7)--(1.5,-7.6);
\draw[->] (4.5,-7)--(4,-7.6);
\draw[->] (5.5,-7)--(6,-7.6);
\draw[->] (8.5,-7)--(8,-7.6);
\draw[->] (9.5,-7)--(10,-7.6);
\draw[->] (12.2,-7)--(11.7,-7.6);
\draw[->] (13.7,-7)--(14.2,-7.6);
\node[scale=0.8] at (-0.2,-9.5) {$6$};
\node[xscale=0.35, yscale=0.43] at (0,-8.5) {$\begin{array}{|cccccc|c|}
\hline
\blacksquare&\blacksquare&\blacksquare&\blacksquare&\blacksquare&\blacksquare&\square\\\hline
\blacksquare&\square&\square&\square&\square&\square&\blacksquare\\
\square&\blacksquare&\square&\square&\square&\square&\blacksquare\\
\square&\square&\blacksquare&\square&\square&\square&\blacksquare\\
\square&\square&\square&\blacksquare&\square&\square&\blacksquare\\
\square&\square&\square&\square&\blacksquare&\square&\blacksquare\\
\square&\square&\square&\square&\square&\blacksquare&\blacksquare\\
\hline \end{array}$};
\node[scale=0.8] at (2,-9.5) {$4+2$};
\node[xscale=0.35, yscale=0.43] at (2,-8.5) {$\begin{array}{|cccc|cc|c|}
\hline
\blacksquare&\blacksquare&\blacksquare&\blacksquare&\square&\square&\blacksquare\\\hline
\blacksquare&\square&\square&\square&\blacksquare&\blacksquare&\square\\
\square&\blacksquare&\square&\square&\blacksquare&\blacksquare&\square\\
\square&\square&\blacksquare&\square&\blacksquare&\blacksquare&\square\\
\square&\square&\square&\blacksquare&\blacksquare&\blacksquare&\square\\\hline
\square&\square&\square&\square&\blacksquare&\square&\blacksquare\\
\square&\square&\square&\square&\square&\blacksquare&\blacksquare\\
\hline \end{array}$};
\node[scale=0.8] at (3.9,-9.5) {$3+3$};
\node[xscale=0.35, yscale=0.43] at (4,-8.5) {$\begin{array}{|ccc|ccc|c|}
\hline
\blacksquare&\blacksquare&\blacksquare&\square&\square&\square&\blacksquare\\\hline
\blacksquare&\square&\square&\blacksquare&\blacksquare&\blacksquare&\square\\
\square&\blacksquare&\square&\blacksquare&\blacksquare&\blacksquare&\square\\
\square&\square&\blacksquare&\blacksquare&\blacksquare&\blacksquare&\square\\\hline
\square&\square&\square&\blacksquare&\square&\square&\blacksquare\\
\square&\square&\square&\square&\blacksquare&\square&\blacksquare\\
\square&\square&\square&\square&\square&\blacksquare&\blacksquare\\
\hline \end{array}$};
\node[scale=0.8] at (6,-9.5) {$3+1+2$};
\node[xscale=0.35, yscale=0.43] at (6,-8.5) {$\begin{array}{|ccc|c|cc|c|}
\hline
\blacksquare&\blacksquare&\blacksquare&\square&\blacksquare&\blacksquare&\square\\\hline
\blacksquare&\square&\square&\blacksquare&\square&\square&\blacksquare\\
\square&\blacksquare&\square&\blacksquare&\square&\square&\blacksquare\\
\square&\square&\blacksquare&\blacksquare&\square&\square&\blacksquare\\\hline
\square&\square&\square&\blacksquare&\blacksquare&\blacksquare&\square\\\hline
\square&\square&\square&\square&\blacksquare&\square&\blacksquare\\
\square&\square&\square&\square&\square&\blacksquare&\blacksquare\\
\hline \end{array}$};
\node[scale=0.8] at (9.9,-9.5) {$2+2+2$};
\node[xscale=0.35, yscale=0.43] at (10,-8.5) {$\begin{array}{|cc|cc|cc|c|}
\hline
\blacksquare&\blacksquare&\square&\square&\blacksquare&\blacksquare&\square\\\hline
\blacksquare&\square&\blacksquare&\blacksquare&\square&\square&\blacksquare\\
\square&\blacksquare&\blacksquare&\blacksquare&\square&\square&\blacksquare\\\hline
\square&\square&\blacksquare&\square&\blacksquare&\blacksquare&\square\\
\square&\square&\square&\blacksquare&\blacksquare&\blacksquare&\square\\\hline
\square&\square&\square&\square&\blacksquare&\square&\blacksquare\\
\square&\square&\square&\square&\square&\blacksquare&\blacksquare\\
\hline \end{array}$};
\node[scale=0.8] at (11.9,-9.5) {$2+1+3$};
\node[xscale=0.35, yscale=0.43] at (12,-8.5) {$\begin{array}{|cc|c|ccc|c|}
\hline
\blacksquare&\blacksquare&\square&\blacksquare&\blacksquare&\blacksquare&\square\\\hline
\blacksquare&\square&\blacksquare&\square&\square&\square&\blacksquare\\
\square&\blacksquare&\blacksquare&\square&\square&\square&\blacksquare\\\hline
\square&\square&\blacksquare&\blacksquare&\blacksquare&\blacksquare&\square\\\hline
\square&\square&\square&\blacksquare&\square&\square&\blacksquare\\
\square&\square&\square&\square&\blacksquare&\square&\blacksquare\\
\square&\square&\square&\square&\square&\blacksquare&\blacksquare\\
\hline \end{array}$};
\node[scale=0.75] at (14,-9.5) {$2+1+1+2$};
\node[xscale=0.35, yscale=0.43] at (14,-8.5) {$\begin{array}{|cc|c|c|cc|c|}
\hline
\blacksquare&\blacksquare&\square&\blacksquare&\square&\square&\blacksquare\\\hline
\blacksquare&\square&\blacksquare&\square&\blacksquare&\blacksquare&\square\\
\square&\blacksquare&\blacksquare&\square&\blacksquare&\blacksquare&\square\\\hline
\square&\square&\blacksquare&\blacksquare&\square&\square&\blacksquare\\\hline
\square&\square&\square&\blacksquare&\blacksquare&\blacksquare&\square\\\hline
\square&\square&\square&\square&\blacksquare&\square&\blacksquare\\
\square&\square&\square&\square&\square&\blacksquare&\blacksquare\\
\hline \end{array}$};
\end{tikzpicture}
\end{center} 
\caption{\small{Upper Hessenberg matrix representations and integer compositions.}}
\label{figureE15}
\end{figure}

The computational results in Section~\ref{ESect-7} show that all acute 0/1-simplices in 
$I^3,I^4$ and~$I^5$ have a minimal matrix representation that is an {\em unreduced 
upper Hessenberg matrix}. For~$n\geq 6$, many, but not all of them are unreduced upper 
Hessenberg. A closer inspection of these matrices shows that each of them corresponds 
to a unique {\em composition} of the integer~$\nmo$.
\begin{Def}[Integer composition] {\rm A {\em composition}\index{integer composition} of an integer~$n$ in 
$k$ parts is an ordered~$k$-tupel~$\lambda=\langle \lambda_1,\dots,\lambda_k\rangle$ 
with~$\lambda_j\in\mathbb{N}$ with the property that~$n = \lambda_1 + \lambda_2 + \dots + \lambda_k$.}
\end{Def}
In Figure \ref{figureE15} we depict the observed correspondence, restricted to~$n\leq 7$, as a 
binary tree. The first~$\nmo$ entries in the first row of each matrix  form an integer 
composition~$\langle \lambda_1,\dots,\lambda_k\rangle$ of~$\nmo$, by considering 
consecutive entries with the {\em same value} as belonging to the {\em same part}. 
The last~$\nmo$ entries of the last column show the same composition.

\smallskip

The horizontal and vertical lines separating the parts of both compositions, 
subdivide the matrix in blocks. There are~$k$ {\em identity matrices} of 
consecutive sizes~$\lambda_1\times\lambda_1, \dots,\lambda_k\times\lambda_k$ 
containing the first sub-diagonal. The blocks above those identity matrices alternatingly 
contain only zeros or only ones, starting with ones directly above the identity matrices.
\begin{rem}{\rm In Figure \ref{figureE15}, the matrix corresponding to the composition~$2+4$ of~$6$ is 
missing. This is because it is not a {\em minimal} matrix representation. The matrix that 
{\em is} a minimal matrix representation of the corresponding simplex is {\em not} upper 
Hessenberg. In Figure~\ref{figureE16}, it is depicted to its right. The same is done for the compositions~$2+5$ and~$2+1+4$ of~$7$.}
\end{rem}
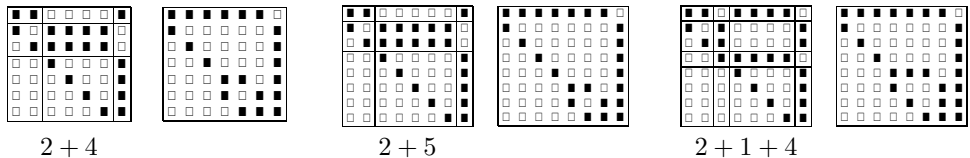
\begin{figure}[h]
\begin{center}
\begin{tikzpicture}[scale=0.9, every node/.style={scale=0.9}]
\node[xscale=0.4, yscale=0.5] at (2.3,0.03) {$\begin{array}{|ccccccc|}
\hline
\blacksquare&\blacksquare&\blacksquare&\blacksquare&\blacksquare&\blacksquare&\square\\
\blacksquare&\square&\square&\square&\square&\square&\blacksquare\\
\square&\blacksquare&\square&\square&\square&\square&\blacksquare\\
\square&\square&\blacksquare&\square&\square&\square&\blacksquare\\
\square&\square&\square&\blacksquare&\blacksquare&\square&\blacksquare\\
\square&\square&\square&\blacksquare&\square&\blacksquare&\blacksquare\\
\square&\square&\square&\square&\blacksquare&\blacksquare&\blacksquare\\
\hline \end{array}$};
\node at (0,-1.2) {$2+4$};
\node at (5,-1.2) {$2+5$};
\node at (10,-1.2) {$2+1+4$};
\node[xscale=0.4, yscale=0.5] at (0,0.03) {$\begin{array}{|cc|cccc|c|}
\hline
\blacksquare&\blacksquare&\square&\square&\square&\square&\blacksquare\\\hline
\blacksquare&\square&\blacksquare&\blacksquare&\blacksquare&\blacksquare&\square\\
\square&\blacksquare&\blacksquare&\blacksquare&\blacksquare&\blacksquare&\square\\\hline
\square&\square&\blacksquare&\square&\square&\square&\blacksquare\\
\square&\square&\square&\blacksquare&\square&\square&\blacksquare\\
\square&\square&\square&\square&\blacksquare&\square&\blacksquare\\
\square&\square&\square&\square&\square&\blacksquare&\blacksquare\\
\hline \end{array}$};
\node[xscale=0.37, yscale=0.46] at (7.3,0) {$\begin{array}{|cccccccc|}
\hline
\blacksquare&\blacksquare&\blacksquare&\blacksquare&\blacksquare&\blacksquare&\blacksquare&\square\\
\blacksquare&\square&\square&\square&\square&\square&\square&\blacksquare\\
\square&\blacksquare&\square&\square&\square& \square&\square&\blacksquare\\
\square&\square&\blacksquare&\square&\square&\square&\square&\blacksquare\\
\square&\square&\square&\blacksquare&\square&\square&\square&\blacksquare\\
\square&\square&\square&\square&\blacksquare&\blacksquare&\square&\blacksquare\\
\square&\square&\square&\square&\blacksquare&\square&\blacksquare&\blacksquare\\
\square&\square&\square&\square&\square&\blacksquare&\blacksquare&\blacksquare\\
\hline \end{array}$};
\node[xscale=0.37, yscale=0.46] at (5,0) {$\begin{array}{|cc|ccccc|c|}
\hline
\blacksquare&\blacksquare&\square&\square&\square&\square&\square&\blacksquare\\\hline
\blacksquare&\square&\blacksquare&\blacksquare&\blacksquare&\blacksquare&\blacksquare&\square\\
\square&\blacksquare&\blacksquare&\blacksquare&\blacksquare&\blacksquare&\blacksquare&\square\\\hline
\square&\square&\blacksquare&\square&\square&\square&\square&\blacksquare\\
\square&\square&\square&\blacksquare&\square&\square&\square&\blacksquare\\
\square&\square&\square&\square&\blacksquare&\square&\square&\blacksquare\\
\square&\square&\square&\square&\square&\blacksquare&\square&\blacksquare\\
\square&\square&\square&\square&\square&\square&\blacksquare&\blacksquare\\
\hline \end{array}$};
\node[xscale=0.37, yscale=0.46] at (12.3,0) {$\begin{array}{|cccccccc|}
\hline
\blacksquare&\blacksquare&\blacksquare&\blacksquare&\blacksquare&\blacksquare&\blacksquare&\square\\
\blacksquare&\square&\square&\square&\square&\square&\square&\blacksquare\\
\square&\blacksquare&\square&\square&\square&\square&\square&\blacksquare\\
\square&\square&\blacksquare&\square&\square&\square&\square&\blacksquare\\
\square&\square&\square&\blacksquare&\blacksquare&\blacksquare&\square&\blacksquare\\
\square&\square&\square&\blacksquare&\square&\square&\blacksquare&\blacksquare\\
\square&\square&\square&\square&\blacksquare&\square&\blacksquare&\blacksquare\\
\square&\square&\square&\square&\square&\blacksquare&\blacksquare&\blacksquare\\
\hline\end{array}$};
\node[xscale=0.37, yscale=0.46] at (10,0) {$\begin{array}{|cc|c|cccc|c|}
\hline
\blacksquare&\blacksquare&\square&\blacksquare&\blacksquare&\blacksquare&\blacksquare&\square\\\hline
\blacksquare&\square&\blacksquare&\square&\square&\square&\square&\blacksquare\\
\square&\blacksquare&\blacksquare&\square&\square&\square&\square&\blacksquare\\\hline
\square&\square&\blacksquare&\blacksquare&\blacksquare&\blacksquare&\blacksquare&\square\\\hline
\square&\square&\square&\blacksquare&\square&\square&\square&\blacksquare\\
\square&\square&\square&\square&\blacksquare&\square&\square&\blacksquare\\
\square&\square&\square&\square&\square&\blacksquare&\square&\blacksquare\\
\square&\square&\square&\square&\square&\square&\blacksquare&\blacksquare\\
\hline \end{array}$};
\end{tikzpicture}
\end{center} 
\caption{\small{Upper Hessenberg matrices corresponding to integer 
compositions~$2+4, 2+5$, and~$2+1+4$. Next to them are depicted the 
corresponding minimal matrix representations.}}
\label{figureE16}
\end{figure}
\begin{Pro} The matrix~$H_\lambda$ corresponding to an integer partition 
$\langle \lambda_1,\dots,\lambda_k\rangle$ is not a minimal matrix representation 
if~$\lambda_j>\lambda_1+1$ for some~$j\in\{2,\dots,k\}$.
\end{Pro}
{\bf Proof. } The first~$\lambda_1$ columns of~$H_\lambda$ should, together with 
the origin, form the largest possible subset of vertices with mutual distances equal to 
two. If~$\lambda_j>\lambda_1+1$, this is not the case. \hfill~$\Box$

\smallskip

Now, let~$n\geq 4$ be given, and let~$\lambda=\langle \lambda_1,\dots,\lambda_k\rangle$ 
be an integer composition of~$\nmo$ with the property that~$\lambda_1\not=1\not=\lambda_k$. 
It can easily be verified that there are~$2^{n-4}$ such compositions.

\smallskip

We will write~$H_\lambda$ for the upper Hessenberg matrix that corresponds to~$\lambda$ according 
to the above description and examples. See also the introduction to this paper. The rule defining the 
tree in Figure~\ref{figureE15} is now depicted in Figure~\ref{figureE17}.
\begin{figure}[h]
\begin{center}
\begin{tikzpicture}
\node at (7,0) {$\lambda = \langle \lambda_1,\dots,\lambda_{k-1},\lambda_k\rangle$};
\node at (4,-1.5) {$\lambda = \langle \lambda_1,\dots,\lambda_{k-1},\lambda_k+1\rangle$};
\node at (10,-1.5) {$\lambda = \langle \lambda_1,\dots,\lambda_{k-1},\lambda_k-1,2\rangle$};
\draw[->] (7,-0.5)--(4.1,-0.9);
\draw[->] (7,-0.5)--(9.9,-0.9);
\draw[fill=black] (7,-0.5) circle [radius=0.05];
\draw[fill=black] (4,-1) circle [radius=0.05];
\draw[fill=black] (10,-1) circle [radius=0.05];
\end{tikzpicture}
\end{center} 
\caption{\small{Splitting rule that defines the binary tree in Figure \ref{figureE15}.}}
\label{figureE17}
\end{figure}
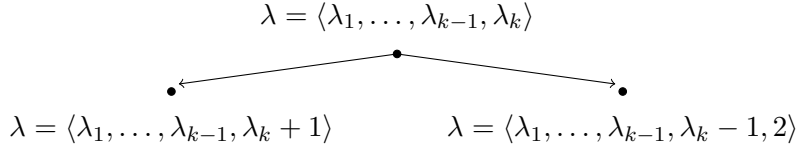\\[2mm]
In the next section we will prove that for given~$n$, each of the~$2^{n-4}$ matrices 
$H_\lambda$ corresponding to a composition~$\lambda$ of~$n-1$ with first and last 
part at least~$2$, represents an acute 0/1-simplex. Conversely, we will show that if~$H$ 
is an~$n\times n$ unreduced upper Hessenberg matrix that represents an acute 0/1-simplex, 
then~$H\sim H_\lambda$ for some composition~$\lambda$ of~$\nmo$.
%%%%%%%%%
%%%%%%%%%
%%%%%%%%
%%%%%%%%%
\subsection{An application of the one neighbor theorem}\label{Esect-7.2}
We first recall the following theorem, which limits the number of candidate acute 
extensions from Definition~\ref{Edef-cand} of an acute simplex in~$I^n$ with~$n$ vertices to two.
\begin{Th}[\cite{BrDiHaKr}]\label{Eth-8} Let~$F$ be an acute 0/1-simplex in~$I^n$ with~$n$ vertices. 
Then $\CC^n(F)\subset\BB^n$ consists of at most two points. If it consists of two points, they add up to~$e=(1,\dots,1)^\top$.
\end{Th}
In the context of triangulations this result is called the {\em one neighbor theorem}, 
since it proves that an acute simplex in~$I^n$ has at most one face-to-face neighbor 
in~$I^n$. See \cite{BrDiHaKr,Ci} for applications of this result in nonobtuse triangulations 
of~$I^n$, and of 0/1-polytopes in~$I^4$, respectively. Here, we will apply the result to 
prove the observed structures in Section~\ref{Esect-7.1}.
\begin{Le}\label{Elem-a1} Let~$H$ be an~$n\times(\nmo)$ unreduced upper 
Hessenberg 0/1-matrix, whose columns together with the origin form an 
acute~$(\nmo)$-simplex in~$I^n$. Then there exist at most two~$(\npo)\times n$ 
unreduced upper Hessenberg matrices whose columns together with the origin form 
an acute~$n$-simplex in~$I^{\npo}$, that have~$H$ as top left~$n\times(\nmo)$ part.
\end{Le}  
{\bf Proof. } Let~$H$ be an~$n\times(\nmo)$ unreduced upper Hessenberg 0/1-matrix, 
whose columns together with the origin form an acute~$(\nmo)$-simplex. Then due to 
Theorem~\ref{Eth-8}, there exists {\em at most} two vertices~$g,h\in\BB^n$ with~$g+h=e$ 
such that the~$n\times n$ matrices~$[H|g]$ and~$[H|h]$ represent acute~$n$-simplices in~$I^n$. As a result, only
\be\label{Eeq-82} \left[\begin{array}{r|r} H & g\\\hline 0 & 1\end{array}\right] \und 
\left[\begin{array}{r|r} H & h\\\hline 0 & 1\end{array}\right] \ee
may be~$(\npo)\times n$ unreduced upper Hessenberg matrices whose columns 
together with the origin form acute~$n$-simplices in~$I^{\npo}$. \hfill~$\Box$
\begin{Co}\label{Ecor-a2} There exist at most~$2^{n-2}$ unreduced upper 
Hessenberg matrices of size~$n\times(\nmo)$ whose columns together with 
the origin are the vertices of an acute~$(\nmo)$-simplex.
\end{Co} 
{\bf Proof. } One can easily verify that in~$I^3$, the matrices
\[ H_1 = \left[\begin{array}{rr} 1 & 1 \\ 1 & 0 \\ 0 & 1 \end{array}\right] \und 
H_2 = \left[\begin{array}{rr} 1 & 0 \\ 1 & 1 \\ 0 & 1 \end{array}\right]\]
are the only two~$3\times 2$ upper Hessenberg matrices whose columns together 
with the origin are acute triangles in~$I^3$, the statement is now proved by 
induction based on Lemma~\ref{Elem-a1}. \hfill~$\Box$

\smallskip

Since~$H_2$ is obtained from~$H_1$ by swapping its first two rows, we see 
that any unreduced upper Hessenberg matrix that is a {\em minimal} matrix 
representation, has~$H_1$ as its top~$3\times 2$ block.
\begin{Co}\label{Ecor-a1} The only two~$(\npo)\times(\npo)$ 
unreduced upper Hessenberg matrices with~$n\times(\nmo)$ 
top left part equal to~$H$ that may represent acute~$(\npo)$-simplices are 
\be \label{Eeq-83} \left[\begin{array}{r|r|r} H & g & h\\\hline 0 & 1 & 1\end{array}\right] 
\und \left[\begin{array}{r|r|r} H & h & g\\\hline 0 & 1 & 1\end{array}\right].  \ee 
In case they do, these two matrices obviously represent the same 0/1-simplex, and 
thus at most one of them can be a minimal matrix representation.
\end{Co}  
{\bf Proof. } Suppose that both~$(\npo)\times n$ matrices in (\ref{Eeq-82}) indeed 
represent acute~$n$-simplices in~$I^{\npo}$. Suppose moreover that adding 
$v\in\BB^{n+1}$ as a~$(\npo)$-st column results in a matrix representing an 
acute~$(\npo)$-simplex. Then due to Theorem~\ref{Eth-8}, the top~$n$ entries 
of~$v$ should consist of~$g$ or~$h$. For the left matrix in (\ref{Eeq-82}) this leads to four options,
\be\label{Eeq-a1}  \left[\begin{array}{r|r|r} H & g & g\\\hline 0 & 1 & 0\end{array}\right], \hdrie  
\left[\begin{array}{r|r|r} H & g & g\\\hline 0 & 1 & 1\end{array}\right], \hdrie  
\left[\begin{array}{r|r|r} H & g & h\\\hline 0 & 1 & 0\end{array}\right] \und  
\left[\begin{array}{r|r|r} H & g & h\\\hline 0 & 1 & 1\end{array}\right].\ee  
We claim that only the rightmost matrix in (\ref{Eeq-a1}) may represent an acute~$(\npo)$-simplex. 
Indeed, the difference between the last two columns of the first matrix is orthogonal to the last. 
Thus, it has a right triangular facet, and thus due to Proposition~\ref{Eprop-6} it cannot represent 
an acute simplex. The second matrix is obviously singular. The last two columns of the third matrix 
are orthogonal because~$g+h=e$ and thus this simplex too has a right triangular facet. Thus, the 
fourth matrix remains. For the right matrix in (\ref{Eeq-82}) a similar analysis can be made. Finally, 
note that the matrices in (\ref{Eeq-83}) differ only by swapping the last columns.\hfill~$\Box$

\smallskip

To be able to fully explain the tree in Figure~\ref{figureE15}, we will need to go one step further, and even 
describe which~$(\npt)\times(\npt)$ unreduced upper Hessenberg matrices with~$n\times(\nmo)$ 
part equal to~$H$ have the potential to be a minimal matrix representation of an acute~$(\npt)$-simplex.
\begin{Co}\label{Ecor-a3} Assume that the right matrix in (\ref{Eeq-83}) is not a minimal matrix representation. Then 
\be\label{Eeq-a2} \left[\begin{array}{c|c|c|c} H & g & g & h \\ \hline 0 & 1 & 0 & 1 \\ \hline 0 & 0 & 1 & 1\end{array}\right]
 \und \left[\begin{array}{c|c|c|c} H & h & h & g \\ \hline 0 & 1 & 0 & 1 \\ \hline 0 & 0 & 1 & 1\end{array}\right]\ee
are the only two~$(\npt)\times(\npt)$ unreduced upper Hessenberg matrices with 
the~$n\times(\nmo)$ unreduced upper Hessenberg matrix~$H$ as top left part, that 
may be minimal matrix representations of an acute~$(\npt)$-simplex. 
\end{Co} 
{\bf Proof. } We follow the lines of the proofs of Lemma~\ref{Elem-a1} and 
Corollary~\ref{Ecor-a1} but then with~$H$ consecutively replaced by each of the 
two matrices in (\ref{Eeq-82}). Consider first the left matrix in (\ref{Eeq-82}). It 
gives rise to the following four candidates,
\be \left[\begin{array}{c|c|c} H & g & v \\ \hline 0 & 1 & a \\ \hline 
0 & 0 & 1\end{array}\right] \hdrie \mbox{\rm with }\hdrie v\in\{g,h\} \und a\in\{0,1\}. \ee
The option~$(v,a)=(g,1)$ and~$(v,a)=(h,0)$ both lead to a right triangular facet 
and are thus infeasible. The remaining options are~$(v,a)=(g,0)$ and~$(v,a)=(h,1)$, 
which, in line with the proof of Lemma~\ref{Elem-a1}, add up to the all-ones vector. 
Moreover, in line with Corollary~\ref{Ecor-a1}, they account for the left matrix in 
 (\ref{Eeq-a2}). For the right matrix in (\ref{Eeq-82}), conversely, only the options 
$(v,a)=(h,0)$ and~$(v,a)=(g,1)$ are feasible, and lead to the right matrix in (\ref{Eeq-a2}). \hfill~$\Box$ 

\smallskip

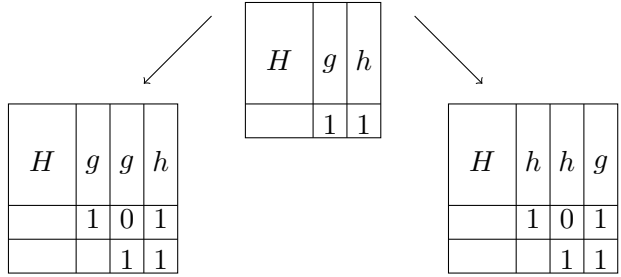
\begin{figure}[h]
\begin{center}
\begin{tikzpicture}[scale=0.9]
\begin{scope}
\draw (0,0)--(2,0)--(2,2)--(0,2)--cycle;
\draw (1,0)--(1,2);
\draw (1.5,0)--(1.5,2);
\draw (0,0.5)--(2,0.5);
\node at (0.5,1.15) {$H$};
\node at (1.25,1.1) {$g$};
\node at (1.75,1.15) {$h$};
\node at (1.25,0.2) {$1$};
\node at (1.75,0.2) {$1$};
\draw[->] (-0.5,1.8)--(-1.5,0.8);
\draw[->] (2.5,1.8)--(3.5,0.8);
\end{scope}
\begin{scope}[shift={(-3.5,-2)}]
\draw (0,0)--(2.5,0)--(2.5,2.5)--(0,2.5)--cycle;
\draw (1,0)--(1,2.5);
\draw (1.5,0)--(1.5,2.5);
\draw (2.0,0)--(2.0,2.5);
\draw (0,0.5)--(2.5,0.5);
\draw (0,1)--(2.5,1);
\node at (0.5,1.65) {$H$};
\node at (1.25,1.6) {$g$};
\node at (1.75,1.6) {$g$};
\node at (2.25,1.65) {$h$};
\node at (1.25,0.8) {$1$};
\node at (1.75,0.8) {$0$};
\node at (2.25,0.8) {$1$};
\node at (1.75,0.2) {$1$};
\node at (2.25,0.2) {$1$};
\end{scope}
\begin{scope}[shift={(3,-2)}]
\draw (0,0)--(2.5,0)--(2.5,2.5)--(0,2.5)--cycle;
\draw (1,0)--(1,2.5);
\draw (1.5,0)--(1.5,2.5);
\draw (2.0,0)--(2.0,2.5);
\draw (0,0.5)--(2.5,0.5);
\draw (0,1)--(2.5,1);
\node at (0.5,1.65) {$H$};
\node at (1.25,1.65) {$h$};
\node at (1.75,1.65) {$h$};
\node at (2.25,1.6) {$g$};
\node at (1.25,0.8) {$1$};
\node at (1.75,0.8) {$0$};
\node at (2.25,0.8) {$1$};
\node at (1.75,0.2) {$1$};
\node at (2.25,0.2) {$1$};
\end{scope}
\end{tikzpicture}
\end{center} 
\caption{\small{Splitting rule that defines the binary tree in Figure \ref{figureE15}, in matrix form.}}
\label{figureE18}
\end{figure}

Corollary~\ref{Ecor-a3} explicitly proves that an unreduced upper Hessenberg matrix 
$H_\lambda$ of size~$(\npo)\times(\npo)$ that represents an acute simplex, can 
have at most two unreduced upper Hessenberg descendants of size 
$(\npt)\times(\npt)$ that represent an acute simplex, and who share their 
$n\times(n-1)$ top left parts. This is depicted in Figure~\ref{figureE18}.
\\[2mm]
This also proves that for each unreduced upper Hessenberg matrix that represents 
an acute 0/1-simplex, there exists a matrix~$H_\lambda$ in the tree in Figure~\ref{figureE15}  
with which it is 0/1-equivalent.

\smallskip

We will now proceed to prove that each of the~$2^{n-4}$ unreduced upper Hessenberg 
matrices in the tree indeed represents an acute simplex. For this, we will use the concept 
of strictly ultrametric matrix, as defined in Section~\ref{sect-E1.1}.
\begin{Th} Let~$H_\lambda$ be the unreduced upper Hessenberg 
matrix corresponding to the integer decomposition 
$\lambda=\langle\lambda_1,\dots,\lambda_k\rangle$ of 
$n-1$. Then~$G_\lambda=H_\lambda^\top H_\lambda$ is strictly ultrametric.
\end{Th}
{\bf Proof.} We use the splitting rule proved in Corollaries~\ref{Ecor-a1} and 
\ref{Ecor-a3} and depicted in Figure~\ref{figureE18} as starting point for an inductive proof. 
Consider the~$n+1$ columns of the parent matrix in Figure~\ref{figureE18}, and write them as
\be\label{Eultr} \left[\begin{array}{c}h_1 \\ \hline 0 \end{array}\right], \dots, 
\left[\begin{array}{c}h_{n-1} \\ \hline 0 \end{array}\right], \hdrie 
\left[\begin{array}{c} g \\ \hline 1 \end{array}\right] \und 
\left[\begin{array}{c}h \\ \hline 1 \end{array}\right]. \ee
By definition of strict ultrametricity, there is no triple~$u,v,w$ of distinct columns 
taken from (\ref{Eultr}) such that one of the three numbers~$u^\top v, v^\top w, w^\top u$ 
is smaller than the other two. We will prove this property also for the~$n+2$ columns 
of the left descendant in Figure~\ref{figureE18}, which are
\be \label{Eultr-2}\left[\begin{array}{c}h_1 \\ \hline 0 \\ 0 \end{array}\right], \dots, 
\left[\begin{array}{c}h_{n-1} \\ \hline 0 \\ 0 \end{array}\right], \hdrie \left[\begin{array}{c} g \\ \hline 
1  \\ 0\end{array}\right],\hdrie \left[\begin{array}{c} h \\ \hline 1  \\ 1\end{array}\right] \hdrie
\mbox{\rm and the new column }\hdrie  \left[\begin{array}{c}g \\ \hline 0 \\ 1\end{array}\right]. \ee
Obviously, each triple taken from (\ref{Eultr-2}) that does not contain the new 
column has the same mutual inner products as a triple from (\ref{Eultr}). The same clearly holds for each of the triples 
\[ \left\{\left[\begin{array}{c}h_i \\ \hline 0 \\ 0 \end{array}\right],\hdrie 
\left[\begin{array}{c}h_j \\ \hline 0 \\ 0 \end{array}\right],\hdrie 
\left[\begin{array}{c}g \\ \hline 0 \\ 1\end{array}\right]  \right\} \und  
\left\{\left[\begin{array}{c}h_i \\ \hline 0 \\ 0 \end{array}\right],\hdrie
\left[\begin{array}{c} h \\ \hline 1  \\ 1\end{array}\right],\hdrie  
\left[\begin{array}{c}g \\ \hline 0 \\ 1\end{array}\right]\right\}. \]
Two possible triples remain to be discussed, being      
\[ \left\{\left[\begin{array}{c}h_i \\ \hline 0 \\ 0 \end{array}\right],\hdrie 
\left[\begin{array}{c}g \\ \hline 1 \\ 0 \end{array}\right],\hdrie 
\left[\begin{array}{c}g \\ \hline 0 \\ 1\end{array}\right]  \right\} \und  
\left\{\left[\begin{array}{c}g \\ \hline 1 \\ 0 \end{array}\right],\hdrie
\left[\begin{array}{c} h \\ \hline 1  \\ 1\end{array}\right],\hdrie  
\left[\begin{array}{c}g \\ \hline 0 \\ 1\end{array}\right]\right\}. \]       
For the left triple, we use the generally valid fact that~$x^\top y\leq x^\top x$ for 
all~$x,y\in\BB^n$ to conclude that~$h_j^\top g \leq g^\top g$, which proves the 
required property. For the right triple it suffices to note that~$g^\top h=0$ and 
$g^\top g\geq 1$. Next, we consider the right descendant in Figure~\ref{figureE18}, with columns
\be \label{Eultr-3}\left[\begin{array}{c}h_1 \\ \hline 0 \\ 0 \end{array}\right], \dots, 
\left[\begin{array}{c}h_{n-1} \\ \hline 0 \\ 0 \end{array}\right], \hdrie \left[\begin{array}{c} h \\ \hline 
1  \\ 0\end{array}\right],\hdrie \left[\begin{array}{c} g \\ \hline 1  \\ 1\end{array}\right] \hdrie
\mbox{\rm and the new column }\hdrie  \left[\begin{array}{c}h \\ \hline 0 \\ 1\end{array}\right]. \ee
Compared to (\ref{Eultr-2}), only the roles of~$g$ and~$h$ have been exchanged. This does not affect the 
validity of the above arguments. Thus, both descendents are strictly ultrametric. Since the~$3\times 3$ 
matrix~$H_\lambda$ with~$\lambda=\langle 2\rangle$ has a strictly ultrametric Gramian, this proves the 
statement for all members~$H_\lambda$ of the tree in Figure~\ref{figureE15} by induction.~\hfill~$\Box$
%%%%%%%%%%%%%%
%%%%%%%%%%%%%\
%%%%%%%%%
\subsection{Determinant of~$H_\lambda$ as continued fraction numerator}
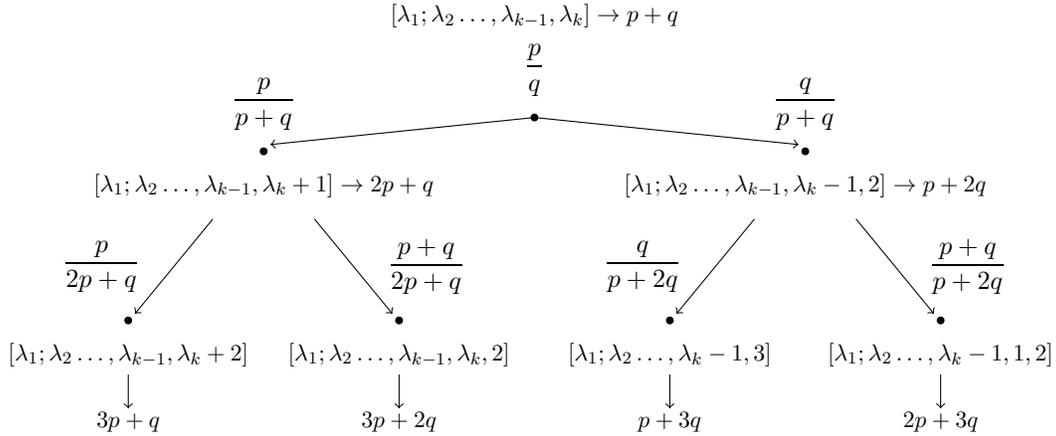
\begin{figure}[h]
\begin{center}
\begin{tikzpicture}[scale=0.9, every node/.style={scale=0.9}]
\node[scale=0.9] at (7,1) {$[\lambda_1;\lambda_2\dots,\lambda_{k-1},\lambda_k]\rightarrow  p+q$};
\node at (7,0.2) {$\dfrac{p}{q}$};
\node[scale=0.9] at (3,-1.5) {$[\lambda_1;\lambda_2\dots,\lambda_{k-1},\lambda_k+1]\rightarrow  2p+q$};
\node[scale=0.9] at (11,-1.5) {$[\lambda_1;\lambda_2\dots,\lambda_{k-1},\lambda_k-1,2]\rightarrow  p+2q$};
\draw[->] (7,-0.5)--(3.1,-0.9);
\draw[->] (7,-0.5)--(10.9,-0.9);
\draw[fill=black] (7,-0.5) circle [radius=0.05];
\draw[fill=black] (3,-1) circle [radius=0.05];
\draw[fill=black] (11,-1) circle [radius=0.05];
\node at (3,-0.3) {$\dfrac{p}{p+q}$};
\node at (11,-0.3) {$\dfrac{q}{p+q}$};
\draw[fill=black] (1,-3.5) circle [radius=0.05];
\draw[fill=black] (5,-3.5) circle [radius=0.05];
\node at (0.6,-2.7) {$\dfrac{p}{2p+q}$};
\node at (5.4,-2.7) {$\dfrac{p+q}{2p+q}$};
\draw[->] (2.25,-2)--(1.1,-3.4);
\draw[->] (3.75,-2)--(4.9,-3.4);
\draw[->] (10.25,-2)--(9.1,-3.4);
\draw[->] (11.75,-2)--(12.9,-3.4);
\draw[fill=black] (9,-3.5) circle [radius=0.05];
\draw[fill=black] (13,-3.5) circle [radius=0.05];
\node at (8.6,-2.7) {$\dfrac{q}{p+2q}$};
\node at (13.4,-2.7) {$\dfrac{p+q}{p+2q}$};
\node[scale=0.9] at (1,-4) {$[\lambda_1;\lambda_2\dots,\lambda_{k-1},\lambda_k+2]$};
\node[scale=0.9] at (5,-4) {$[\lambda_1;\lambda_2\dots,\lambda_{k-1},\lambda_k,2]$};
\node[scale=0.9] at (9,-4) {$[\lambda_1;\lambda_2\dots,\lambda_{k}-1,3]$};
\node[scale=0.9] at (13,-4) {$[\lambda_1;\lambda_2\dots,\lambda_{k}-1,1,2]$};
\node[scale=0.9] at (1,-5) {$3p+q$};
\node[scale=0.9] at (5,-5) {$3p+2q$};
\node[scale=0.9] at (9,-5) {$p+3q$};
\node[scale=0.9] at (13,-5) {$2p+3q$};
\draw[->] (1,-4.3)--(1,-4.8);
\draw[->] (5,-4.3)--(5,-4.8);
\draw[->] (9,-4.3)--(9,-4.8);
\draw[->] (13,-4.3)--(13,-4.8);
\end{tikzpicture}
\end{center} 
\caption{\small{If parent and grandparent satisfy the statement, so do the four descendants.}}
\label{figureE19}
\end{figure}

We are now able to derive an explicit expression for the determinant of 
$H_\lambda$ for any given integer composition~$\lambda$. For this, we 
associate with the parent matrix in Figure~\ref{figureE18} two integers:
\be H_\lambda = \left[\begin{array}{c|c|c} H & g & h \\ \hline 0 & 1 & 1\end{array}\right] \rightarrow  
\left[\begin{array}{r} p \\ q\end{array}\right] = \left[\begin{array}{r} \det(H|g) \\ \det(H|h)\end{array}\right]. \ee
By developing the last row of~$H_\lambda$, we have that~$\det(H_\lambda) = p-q$, 
whereas for the descendants~$H_\lambda^\ell$ and~$H_\lambda^r$ of~$H_\lambda$, 
also by development of their last rows, we find that
\be \left[\begin{array}{c|c|c|c} H & g & g & h \\ \hline 0 & 1 & 0 & 1 \\ \hline 
0 & 0 & 1 & 1\end{array}\right] \rightarrow \left[\begin{array}{c} -p  \\ p-q \end{array}\right] \und 
\left[\begin{array}{c|c|c|c} H & h & h & g \\ \hline 0 & 1 & 0 & 1 \\ \hline 
0 & 0 & 1 & 1\end{array}\right] \rightarrow \left[\begin{array}{c} -q  \\ q-p\end{array}\right]. \ee
Since~$p=1$ and~$q=-2$ for the matrix~$H_\lambda$ with~$\lambda=\langle 3\rangle$ at 
the root of the tree in Figure~\ref{figureE15}, we see that this explains the correspondence between 
the absolute determinants of the matrices~$H_\lambda$ and Kepler's Tree of Fractions, 
as claimed in Section~\ref{sect-E1.2}. To additionally  prove the statement in 
Theorem~\ref{Eth-main} that~$\det(H_\lambda)$ equals the numerator~$f_k$ of the continued fraction
\be [\lambda_1;\lambda_2,\dots,\lambda_k] = \lambda_1 + \cfrac{1}{\lambda_2 + \cfrac{1}{\ddots + \cfrac{1}{\lambda_k} } } 
=\frac{f_k}{g_k} \hdrie \mbox{\rm with~$f_k,g_k$ coprime,} \ee
we use the well-known result from continued fraction theory that~$f_k$ can be computed from the two-term recursion 
\be\label{Ecf-rec} f_j = \lambda_j f_{j-1} + f_{j-2}, \hdrie \mbox{\rm with} \hdrie f_{0}=1 \und f_{-1}=0. \ee
We inductively assume that the statement holds for both the parent 
and the grandparent of a vertex in the tree, and prove the statement for the descendants, as depicted in Figure~\ref{figureE19}.
\\[2mm]
First observe that all seven continued fractions in Figure~\ref{figureE19} start with 
the same~$k-1$ numbers~$\lambda_1,\dots,\lambda_{k-1}$. Denote the 
numerator of~$[\lambda_1;\lambda_2,\dots,\lambda_{k-2}]$ by~$f_{k-2}$ 
and the numerator of~$[\lambda_1;\lambda_2,\dots,\lambda_{k-1}]$ by 
$f_{k-1}$. Then the induction hypothesis on the grand parent in Figure~\ref{figureE19} together with (\ref{Ecf-rec}) imply that
\[ p+q = \lambda_k f_{k-1} + f_{k-2}, \]
whereas the induction hypothesis on its left descendant translates to
\[ 2p+q = (\lambda_k+1) f_{k-1} + f_{k-2}. \]
From thes two relations we can solve~$f_{k-1}$ and~$f_{k-2}$ as 
\[ f_{k-1} = p \und f_{k-2} = (1-\lambda_k)p+q. \]
It remains to verify whether these values for~$f_{k-1}$ and~$f_{k-2}$ are 
consistent with the remaining five continued fractions in Figure~\ref{figureE19} and the 
expressions of their numerators in terms of~$p$ and~$q$. First, for the continued 
fraction~$[\lambda_1;\lambda_2,\dots,\lambda_{k}-1,2]$ we find, taking two steps of (\ref{Ecf-rec}) that its numerator indeed equals
\[ 2\cdot\left[(\lambda_k-1) f_{k-1} + f_{k-2}\right]+f_{k-1} = 2\cdot \left[(\lambda_k-1)p + (1-\lambda_k)p+q\right]+p = p+2q.\]
Thus, if a vertex and its left descendant satisfy the statement, then so does 
its right descendant. Consequently, the only two continued fractions to verify 
are the two left descendants at the lowest level in Figure~\ref{figureE19}. For the continued 
fraction~$[\lambda_1;\lambda_2,\dots,\lambda_k+2]$ we find
\[ (\lambda_k+2)f_{k-1}+f_{k-2} = (\lambda_k+2)p +(1-\lambda_k)p+q = 3p+q, \]
and for~$[\lambda_1;\lambda_2,\dots,\lambda_k-1,3]$
\[ 3\cdot[(\lambda_k-1)f_{k-1}+f_{k-2}]+f_{k-1} = 3\cdot[(\lambda_k-1)p+(1-\lambda_k)p+q]+p=p+3q.\]
Since both are consistent, this finishes the induction step. As the induction basis for the 
first two levels of the tree in Figure~\ref{figureE15} is easily verified, this completes the induction proof.
%%%%%%%%%%%%%%

\section{Computational results}\label{ESect-7}
In this final section of this paper we present a selection of the computational 
data obtained by implementations of the algorithms presented. For simplicity, 
we chose Matlab as programming environment, as Matlab contains useful built-in 
functionalities in the area of linear algebra. Faster imlementations can of course 
be obtained using a lower level programming language.
 
\subsection{The cycle index~$Z_n$ of~$\Bn$ for  the values~$n\in\{3,\dots,9\}$}
In Algorithm 1 in Section~\ref{ESect-1}, we described how to compute the cycle index 
$Z_n$ for the induced permutations of~$\BB^n$ by the hyperoctahedral group~$\Bn$. 
The implementation of this algorithm yields each cycle index~$Z_n$ as a table, 
see Figure~\ref{figureE20},~\ref{figureE21},~\ref{figureE22}. These tables 
are a condensed form Table~\ref{EperminS8}, in the sense that zero columns have 
been removed, and zero entries disregarded, see Table~\ref{figureE20}, \ref{figureE21}, \ref{figureE22}.

\smallskip

To generate partitions needed in Algorithm 1, we used Algorithm P 
in \cite{knu}. We also used the most efficient way to determine the cycle type 
of a given permutation, which is~$\mathcal{O}(p)$ for a permutation of~$p$ 
objects. If~$n\geq 10$ then all cycle type computations take more than ninety 
percent of the total computational time in computing~$Z_n$, and this percentage 
increases for increasing~$n$. Thus, no additional improvements of the algorithm can be expected.

\subsection{The number of 0/1-polytopes with~$k$ vertices} 
Using Algorithm 2 from Section~\ref{ESect-2}, we computed the number 
of 0/1-polytopes with~$k$ vertices for~$0\leq k\leq 2^n$, for the values 
$n\leq 5$. They are displayed in Table~\ref{figureE23}. For~$n=5$, only half of the 
results are displayed, as the results for~$k=\ell$ and~$k=2^n-\ell$ are the same. 

In Table~\ref{Etable14}, we zoom in on the 0/1-simplices in~$I^n$ with~$k\leq n+1$ 
vertices. In Table~\ref{Etable15} we present the number of ways to choose~$k$ points 
from~$\BB^n$. Comparing these numbers with the corresponding numbers 
in Table~\ref{Etable14} shows the large gain of working modulo the action of~$\Bn$.
\begin{rem} {\rm Note that for~$k\geq 4$ these numbers include 
{\em degenerate}~$k$-simplices, which lie in a hyperplane of dimension 
less than~$k$. For~$k\in\{2,3\}$ such degenerate cases do not exist: three distinct point in~$I^n$ are never colinear}.
\end{rem}
Finally, in Table~\ref{Etable16} we compare the number~$a(n)$ of {\em acute} 
0/1-simplices in~$I^n$ with their {\em total} number~$s(n)$, both modulo the action of~$\Bn$. 
%%%%%%%%%
%%%%%%%%%
%%%%%%%%%
\subsection{Minimal matrix representations of acute 0/1-simplices}\label{Esect-6.2}
Here we present the computed minimal matrix representations 
(without their zero first column) of all acute 0/1-simplices with~$n+1$ 
vertices in~$I^n$ for~$3\leq n \leq 9$ together with the absolute values of their determinants. 
%%%%%%%%%%%%%%%
%%%%%%%%%%%%%%
%%%%%%%%%%%%%%%
%%%%%%%%%%%%%%%
%\subsubsection{Minimal matrix representations of acute 0/1-$3,4,5,6,7$-simplices}
There are~$1,1,2$ acute 0/1-simplices in~$I^3,I^4,I^5$ respectively modulo 
the action of the hyperoctahedral group. The absolute values of the determinants 
of their minimal matrix representatives given below are in the set
\be \det_3 = \{2\}, \hdrie \det_4 = \{3\}, \hdrie \det_5 = \{4,5\}.\ee
There are~$6$ acute 0/1-simplices in~$I^6$ modulo the action of 
$\mathcal{B}_6$.  The absolute determinants of their minimal matrix representatives given below are in the set
\be \det_6 = \{ 5,7,8,9\}.\ee
In~$I^7$ there are~$13$ acute 0/1-simplices modulo the action of 
$\mathcal{B}_7$.  The determinants of their minimal matrix representatives given below are in the set
\be \det_7 = \{ 6, 9, 10, 11, 12, 13, 14, 24, 32\}.\ee

%\subsubsection{Minimal matrix representations of acute 0/1-$8$-simplices}\label{ESect-7.3.2}
There are~$29$ acute 0/1-$8$-simplices in~$I^8$ modulo the action 
of~$\mathcal{B}_8$.  The determinants of their minimal matrix representatives given below are in the set
\be \det_8 = \{7, 11, 13, 14, 15, 16, 17, 18, 19, 20, 21, 22, 23, 40, 44, 56\}.\ee

\begin{figure}[!]
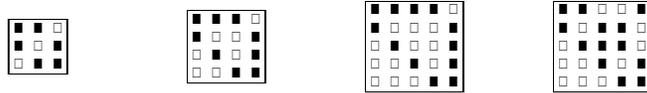

\begin{center}
%
\end{center} 
\caption{\small{Minimal matrix representations of all acute 0/1-simplices in~$I^3, I^4$ and~$I^5$.}}
\label{figureE224}
\end{figure}

\begin{figure}[!]
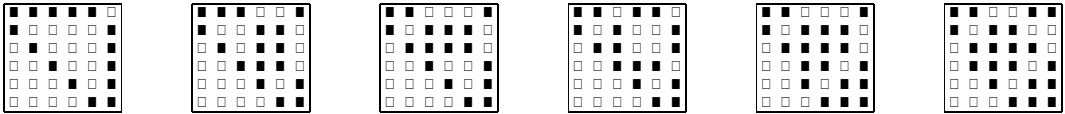

\begin{center}
%
\end{center} 
\caption{\small{Minimal matrix representations of all six acute 0/1-simplices in~$I^6$.}}
\label{figureE225}
\end{figure}

\begin{figure}[!]
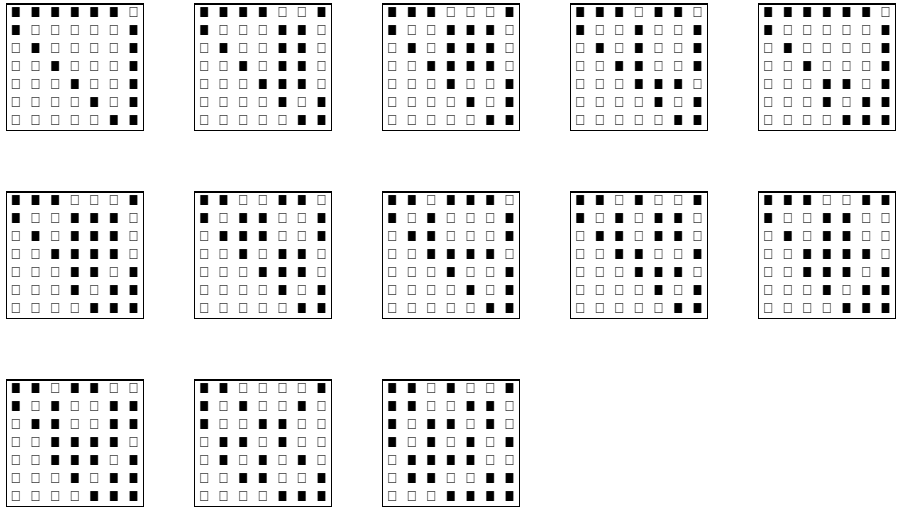

\begin{center}
%
\end{center} 
\caption{\small{Minimal matrix representations of all thirteen acute 0/1-simplices in~$I^7$.}}
\label{figureE226}
\end{figure}

\begin{figure}[!]
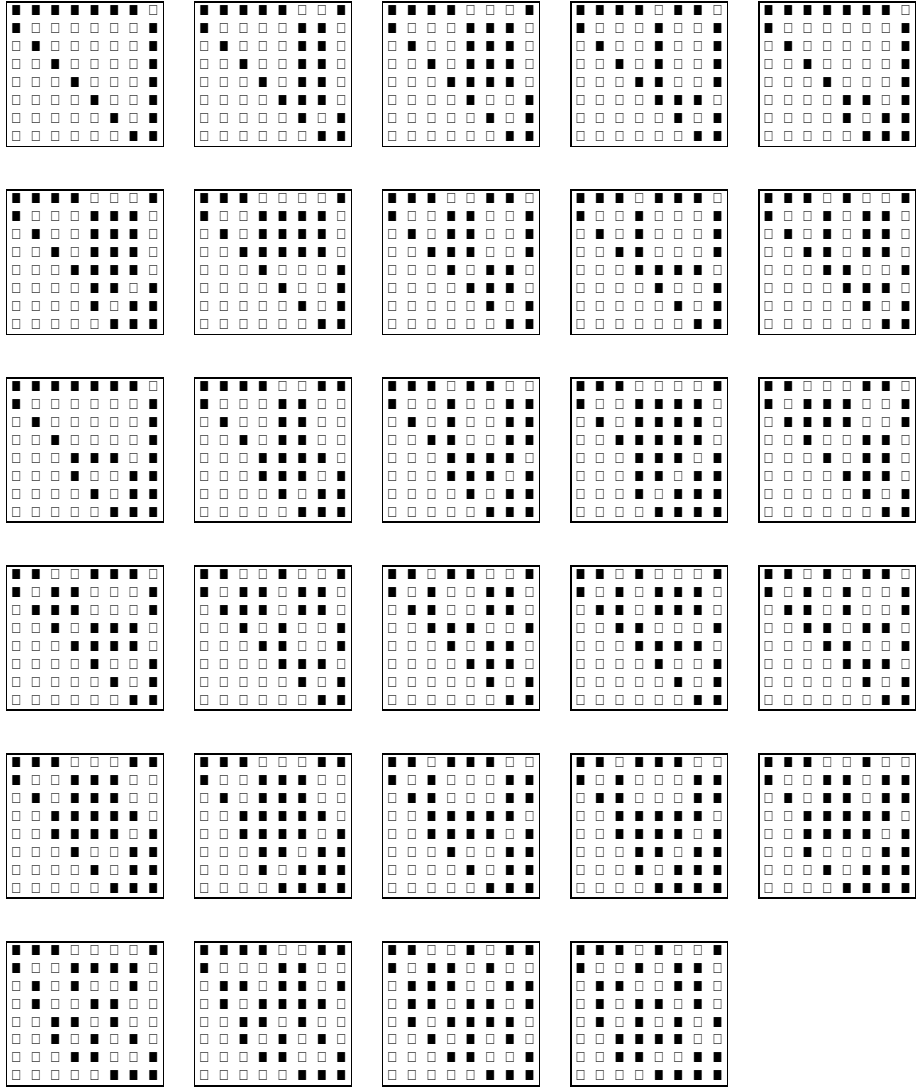

\begin{center}
%
\end{center} 
\caption{\small{Minimal matrix representations of all twenty-nine acute 0/1-$8$-simplices.}}
\label{figureE227}
\end{figure}

\begin{figure}[!]
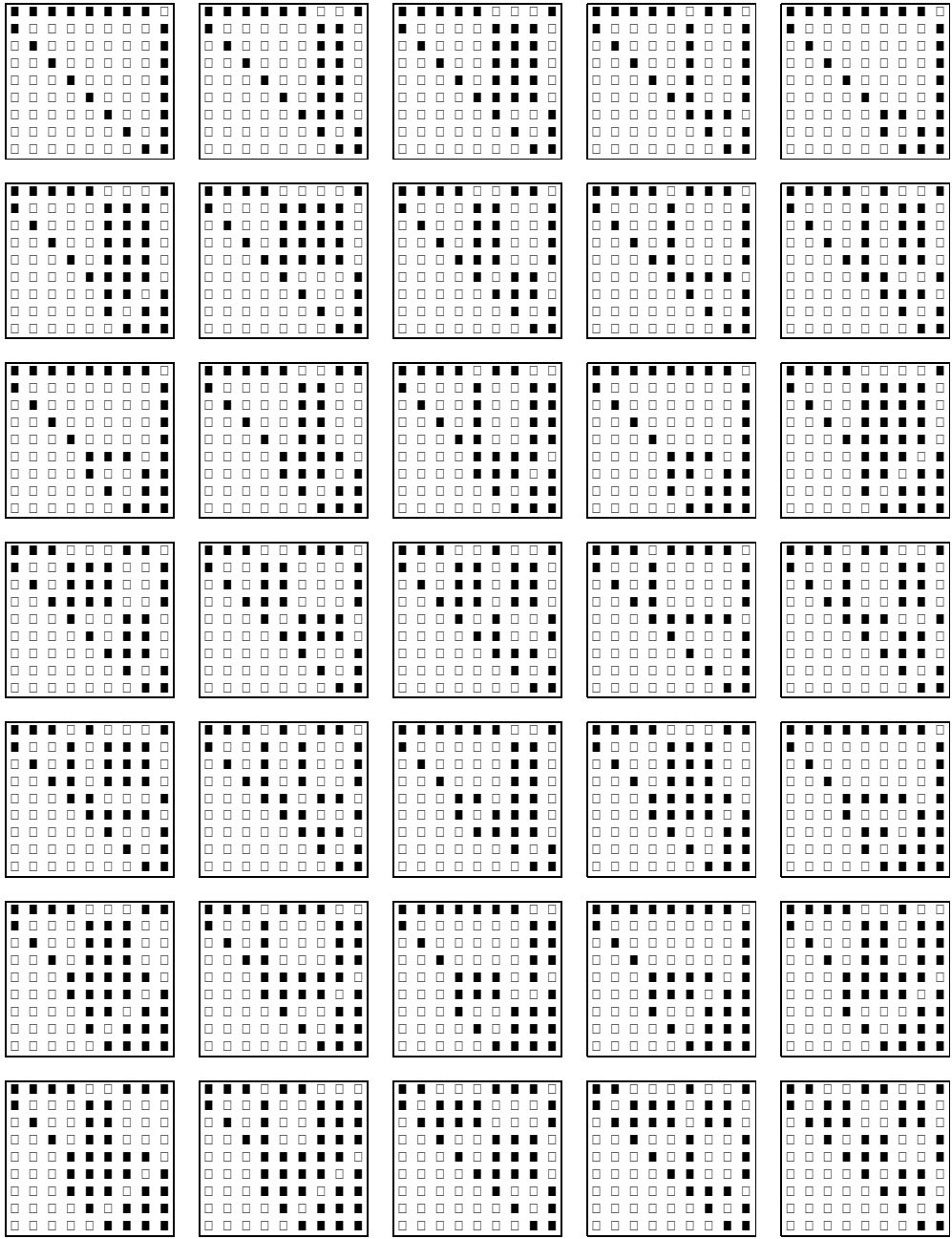

\begin{center}
%
\end{center} 
\caption{\small{Minimal matrix representations of the first~$35$ acute 0/1-$9$-simplices.}}
\label{figureE228}
\end{figure}

\begin{figure}[!]
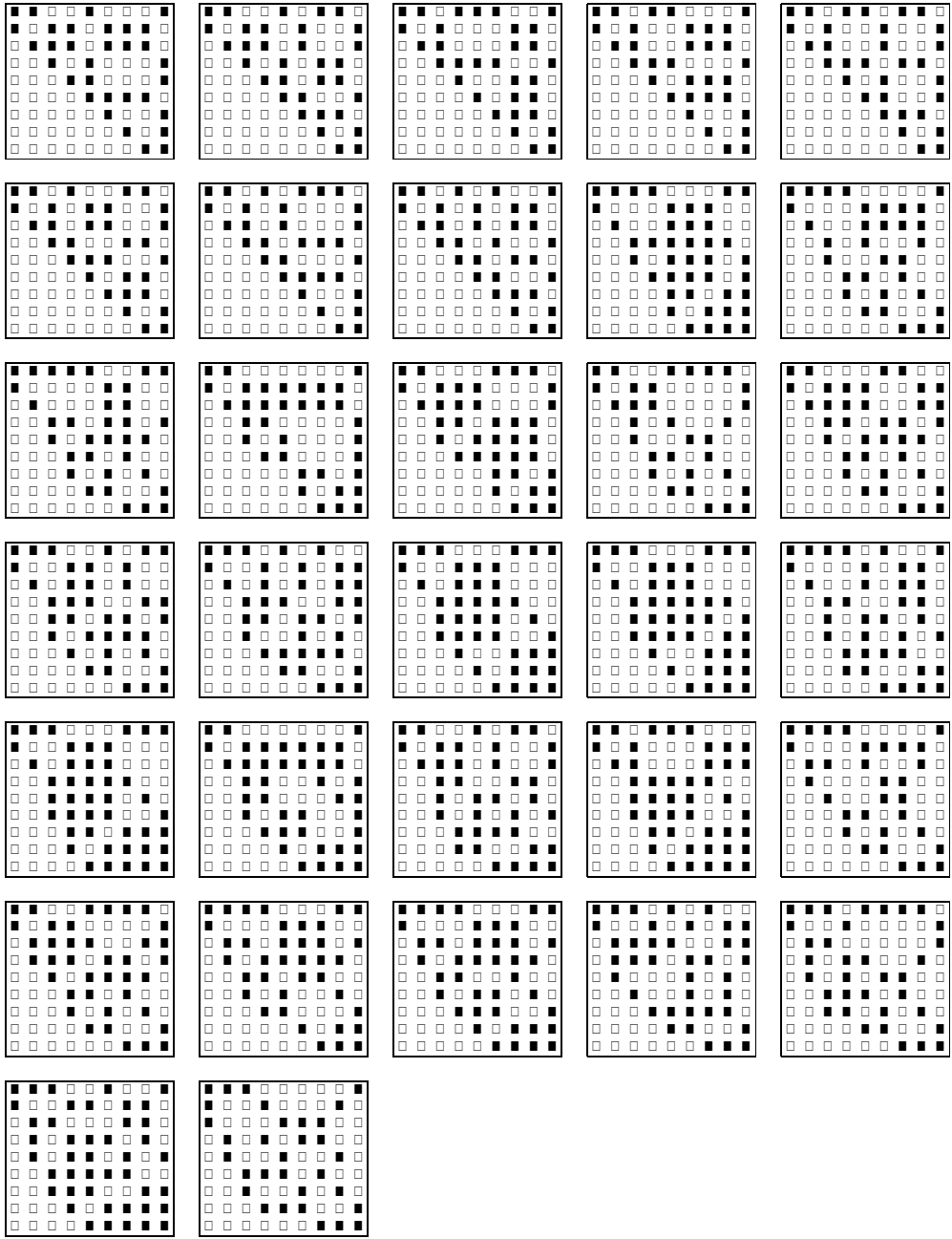

\begin{center}
%
\end{center} 
\caption{\small{Minimal matrix representations of the remaining acute 0/1-$9$-simplices.}}
\label{figureE229}
\end{figure}

%\subsubsection{Minimal matrix representations of acute 0/1-$9$-simplices}
There are~$67$ acute 0/1-$9$-simplices in~$I^9$ modulo the 
action of~$\mathcal{B}_9$.  Their absolute determinants are
\[\det_9=\{8, 13, 16, 17, 19, 20, 21, 22, 23, 24, 25, 26, 27, 28, 29, 30, 31, 32, 34, 35, 45\}\]
\[   \cup \{56, 64, 68, 72, 80, 88, 96\}. \]
In Figure \ref{figureE224}, \ref{figureE225}, \ref{figureE226}, \ref{figureE227}, \ref{figureE228}, \ref{figureE229} we depicte 
the minimal matrix representations of all acute 0/1-$n$-simplices for $3 \leq n \leq 9$. 

\newpage

\begin{figure}[!]
\begin{center}

\end{center} 
\captionof{table}{\small{The cycle indices~$Z_3,Z_4,Z_5,Z_6$ and~$Z_7$ in condensed tabulated form.}}
\label{figureE20}
\end{figure}
\begin{figure}[!]
\begin{center}
%
\end{center} 
\captionof{table}{\small{The cycle index~$Z_8$ in condensed tabulated form.}}
\label{figureE21}
\end{figure}
\begin{figure}[!]
\begin{center}
%
\end{center} 
\captionof{table}{\small{The cycle index~$Z_9$ in condensed tabulated form.}}
\label{figureE22}
\end{figure}
%%%%%%%%%%%
%%%%%%%%%%
%%%%%%%%%%%%%
\begin{figure}[!]
\begin{center}
%
\end{center} 
\captionof{table}{\small{Number of $0/1$-polytopes in $I^n$ with $0\leq k\leq 2^n$ vertices for $n\leq 5$.}}
\label{figureE23}
\end{figure}
\begin{table}
\begin{center}
\small{
\begin{tabular}{|r|r|r|r|r|r|r|r|r|r|}
\hline
$n\backslash k$& 1 & 2 & 3 & 4 & 5 & 6 & 7 & 8 & 9\\
\hline\hline
1&1&1&&&&&&&\\\hline
2&1&2&1&&&&&&\\\hline
3&1&3&3&6&&&&&\\\hline
4&1&4&6&19&27&&&&\\\hline
5&1&5&10&47&131&472&&&\\\hline
6&1&6&16&103&497&3253&19735&&\\\hline
7&1&7&23&203&1606&18435&221778&2773763&\\\hline
8&1&8 & 32 & 373 &4647 & 91028 & 2074059 & 51107344 & 1245930065\\ \hline
\end{tabular}}
\end{center}
\caption{\small{Number of 0/1-simplices in~$I^n$ with~$1\leq k\leq n+1$ vertices for~$n\leq 8$.}}
\label{Etable14}
\end{table}
%%%%%%%%%
%%%%%%%%%%
\begin{table}[!]
\begin{center}
\small{
\begin{tabular}{|r|r|r|r|r|r|r|r|}
\hline
$n\backslash k$& 1 & 2 & 3 & 4 & 5 & 6 & 7\\
\hline\hline
1&2&1&&&&&\\\hline
2&4&6&4&&&&\\\hline
3&8&28&56&70&&&\\\hline
4&16&120&560&1820&4368&&\\\hline
5&32&496&4960&35960&201376&906192&\\\hline
6&64&2016&41664&635376&7624512&74974368&621216192\\\hline
\end{tabular}}
\end{center}
\caption{\small{ Binomial coefficients~$\binom{2^n}{k}$ for comparison with Table \ref{Etable14}.}}
\label{Etable15}
\end{table}
%%%%%%%%%%
%%%%%%%%%%%
\begin{table}[!]
\begin{center}
\small{
\begin{tabular}{|r||r|r|r|r|r|r|r|r|r|r|r|r|r|r|}
\hline
$n$ &1 & 2 & 3 & 4 & 5 & 6 & 7 & 8 & 9 & 10 & 11 \\
 \hline\hline
$a(n)$ &1 & 0 & 1 & 1 & 2 & 6 & 13 & 29 & 67 & 162 & 392\\ 
  \hline
$s(n)$ &1 & 1 & 6 & 27 & 472 & 19735 & 2773763 & 1245930065 & 1.8e12 & 8.7e15& 1.3e20\\
 \hline 
\end{tabular}}
\end{center}
\caption{\small{The number~$a(n)$ of {\em acute} 0/1~$n$-simplices related to their {\em total} number~$s(n)$.}}
\label{Etable16}
\end{table}

\newpage

%%%%%%%%%%%%%%%%%%%%

\end{document}